%% file: jacobian-fin-full.tex
\documentclass[12pt]{amsart}
\usepackage{longtable}
\usepackage{amsfonts,amssymb,amscd,amsmath,epsfig}
\usepackage{latexsym}
\usepackage{mathrsfs}
\usepackage{tocvsec2}

\usepackage[
hyperindex=true,pagebackref=true,
bookmarks=true,
colorlinks=true,linkcolor=blue,citecolor=red]{hyperref}

\begin{document}

\input {macro.tex}

\title [DAHA and plane curve singularities]
{DAHA and plane curve singularities}
\author[Ivan Cherednik]{Ivan Cherednik $^\dag$}
\author[Ian Philipp] {Ian Philipp}

\begin{abstract}
We suggest a relatively simple and totally geometric
conjectural description of uncolored
DAHA superpolynomials 
of arbitrary algebraic knots (conjecturally coinciding with
the reduced stable Khovanov-Rozansky polynomials) via
the flagged Jacobian factors (new objects)
of the corresponding unibranch plane curve singularities. 
This generalizes
the Cherednik-Danilenko conjecture on the Betti numbers of 
Jacobian factors, the Gorsky combinatorial conjectural 
interpretation of superpolynomials of torus knots and that by
Gorsky-Mazin for their constant term. The paper mainly focuses on 
non-torus algebraic knots. A connection with the conjecture due 
to Oblomkov-Rasmussen-Shende is possible, but our approach is 
different. A motivic version of our conjecture is related to 
p-adic orbital A-type integrals for  anisotropic 
centralizers.
\end{abstract}

\thanks{$^\dag$ \today.
\ \ \ Partially supported by NSF grant
DMS--1363138}

\address[I. Cherednik]{Department of Mathematics, UNC
Chapel Hill, North Carolina 27599, USA\\
chered@email.unc.edu}

\address[I. Philipp] {Department of Mathematics, UNC
Chapel Hill, North Carolina 27599, USA\\
iphilipp@live.unc.edu}

 \def\sht{\raisebox{0.4ex}{\hbox{\rm{\tiny sht}}}}
 \def\bysame{{\bf --- }}
 \def\~{{\bf --}}
 \def\rr{{\mathsf r}}
 \def\cc{{\mathsf c}}
 \def\ss{{\mathsf s}}
 \def\mm{{\mathsf m}}
 \def\pp{{\mathsf p}}
 \def\ll{{\mathsf l}}
 \def\aa{{\mathsf a}}
 \def\bb{{\mathsf b}}
 \def\NS{\hbox{\tiny\sf ns}}
 \def\ssum{\hbox{\small$\sum$}}
\newcommand{\comment}[1]{}
\renewcommand{\tilde}{\widetilde}
\renewcommand{\hat}{\widehat}
\renewcommand{\V}{\mathbb{V}}
\renewcommand{\F}{\mathbb{F}}
\newcommand{\dagx}{\hbox{\tiny\mathversion{bold}$\dag$}}
\newcommand{\ddagx}{\hbox{\tiny\mathversion{bold}$\ddag$}}
\newtheorem{conjecture}[theorem]{Conjecture}
\newcommand*\toeq{
\raisebox{-0.15 em}{\,\ensuremath{
\xrightarrow{\raisebox{-0.3 em}{\ensuremath{\sim}}}}\,}
}
\newcommand{\unknot}{\hbox{\tiny\!\raisebox{0.2 em}{$\bigcirc$}}}

\vskip -0.0cm
\maketitle

{\small
\centerline{{\bf MSC} (2010): 14H50, 17B22, 17B45, 20C08, 20F36,}
\centerline{ 22E50, 22E57, 30F10, 33D52, 33D80, 57M25}
}
\smallskip

\vskip -0.0cm
\renewcommand{\baselinestretch}{1.2}
{\small
\tableofcontents   
}
\renewcommand{\baselinestretch}{1.0}

\vfill\eject

\renewcommand{\natural}{\wr}

\setcounter{section}{-1}
\setcounter{equation}{0}
\section{\sc Introduction}

\vfill
We propose a relatively simple and totally computable
conjectural {\em geometric\,} description of  uncolored
DAHA superpolynomials of arbitrary algebraic knots in terms of
{\em flagged Jacobian factors} (new objects)
of the corresponding unibranch 
plane curve singularities, presumably coinciding
with the corresponding {\em stable Khovanov-Rozansky 
polynomials}. This description 
significantly generalizes (a) Cherednik-Danilenko's 
conjecture on the Betti numbers of Jacobian factors
(any unibranch singularities), 
(b) Gorsky's conjectural interpretation of 
superpolynomials of torus knots from \cite{Gor1} etc., 
and (c) that from \cite{GM1,GM2} for their $a$\~constant term.
Our conjecture  is different from the {\em ORS conjecture\,}
from \cite{ORS,GORS} (though some connection is not
impossible).

\vskip 0.2cm
{\sf Motivation.}
Algebraic-geometric theory of topological
invariants of algebraic links has a long history,
starting with the well-known algebraic interpretation
of the {\em Alexander polynomials\,}. This paper provides an 
algebro-geometric description of
stable Khovanov-Rozansky  polynomials via the 
DAHA\~superpolynomials. 
See e.g. \cite{KhR1,KhR2,Kh,Ras,WW}. The geometry of 
flagged Jacobian factors conjecturally provides the 
DAHA\~superpolynomials from \cite{CJ}--\cite{ChD2}.
For instance, this explains the positivity of the latter
in the uncolored case, conjectured in \cite{ChD1}. 

For uncolored {\em torus knots\,}, this positivity 
results from the combinatorial construction
from \cite{Gor1,Gor2}, which conjecturally provides both, 
stable {\em KhR\,}-polynomials and those via DAHA (and is closely
related to rational DAHA). Our conjecture makes the 
positivity entirely geometric (generalizing \cite{GM1,GM2})
for torus and arbitrary algebraic knots.
We expect important implications in the theory 
of plane curve singularities, the theory of $p$\~adic  
{\em orbital integrals\,} and {\em affine Springer fibers\,};
see Conjecture  \ref{CONJSUP} $(iii)$
and Section \ref{sec:CONCL}.

\vskip 0.2cm
{\sf Algebraic knots.}
Torus-type (quasi-homogeneous) plane 
singularities are very special.
Not much is actually known on 
the Jacobian factors of non-torus plane singularities; the
paper \cite{Pi} still remains the main source of examples. 
It was an important development when the DAHA 
approach from \cite{CJ,GoN,CJJ} was extended from torus knots to 
arbitrary algebraic knots in \cite{ChD1} 
and then to any algebraic link in \cite{ChD2}. 
The {\em Newton pairs\,} and
the theory of {\em Puiseux expansion\,}, the key
in the {\em topological\,} 
classification of plane curve singularities,
naturally emerge in the DAHA approach. 

\vfill
\vskip 0.1cm
One of the key advantages of the usage of DAHA is that adding 
colors is relatively direct (via the Macdonald polynomials),
which is well understood for any {\em iterated torus 
link\,} (including all algebraic links). This is well ahead of 
any other approaches (topology included) for such links.
We expect that our present paper can be enhanced by adding
colors via (presumably) the curves suggested in \cite{Ma}. 
The case
of rectangle Young diagrams is exceptional due to 
the conjectured positivity of the corresponding
{\em reduced\,} DAHA 
superpolynomials for algebraic knots \cite{CJ,ChD1}.
The switch from  the rank-one torsion free modules
in the definition of compactified Jacobians 
to arbitrary ranks is expected here (among other modifications),
which is in progress.

\vskip 0.0cm
The passage to arbitrary algebraic knots  and links from
torus knots is important because of multiple reasons. The
generality is an obvious advantage, but not the only one.
All algebraic links (not only torus knots)
are necessary to employ the technique 
of the {\em resolved conifold\,} and similar tools
used in \cite{Ma} to prove the (colored generalization of) the
{\em OS conjecture\,} \cite{ObS} concerning the 
HOMFLYPT polynomials. 
Also, all algebraic {\em links\,} are needed for the theory 
of Hitchin and affine Springer fibers, since 
{\em spectral curves\,} are generally not unibranch. 
Topologically, the class of iterated links is  closed 
with respect to {\em cabling\,}, a major operation in knot theory.

\vskip 0.2cm
{\sf ORS Conjecture.}
Let us briefly comment on 
Conjecture 2 from \cite{ORS};
see Section \ref{sec:CONCL} 
for some further discussion.
It relates the geometry of 
{\em nested Hilbert schemes\,} of arbitrary (germs of)
plane curve singularities to the Khovanov-Rozansky 
{\em unreduced\,} stable polynomials 
of the corresponding links. The main 
component of their conjecture vs. the OS conjecture
is the {\em weight filtration\,}.
The polynomials $\mathscr{P}_{alg}$ there
conjecturally coincide with uncolored ones in the $q,t$-DAHA
theory (upon switching to the standard parameters);
see \cite{ChD1}. They are connected
with the {\em perverse filtartion\,} on the
cohomology of the compactified Jacobians from \cite{MY,MS} 
(see Proposition 4 in \cite{ORS}).

\vskip 0.0cm
Our approach is based on 
{\em admissible\,} {\em flags\,} of submodules $\{\m_i\}$
in the normalization ring;
here dim$_{\C}(\m_{i+1}/\m_{i})=1$ but they are not
full flags and the admissibility is a very restrictive
condition. 
The absence of (nested) Hilbert schemes is 
due to the {\em reduced setting} of our paper 
(continuing \cite{ChD1}); there are other important  
deviations from \cite{ORS}. For instance, we do not need the weight 
filtration and our approach is quite computable. There may be a 
connection with Section 9.1 from \cite{GORS} (a reduced
version of the construction of \cite{ORS}) but this is unclear.
Actually, the weight filtration appears naturally in 
(\ref{conjmotxx}), which follows from a modular variant 
(\ref{conjmot}) of our conjecture, but it is associated with a 
parameter different from that in \cite{ORS}.

\vskip 0.2cm
{\sf Main results.} 
The key is Conjecture  \ref{CONJSUP}; anything else
is about confirmations, examples and connections.
It extends Conjecture 2.4 $(iii)$ from \cite{ChD1}
for Betti numbers of Jacobians factors for unibranch
plane curve singularities (the case $q\!=\!1,a\!=\!0$).
It was 
essentially checked in \cite{Mel} for torus knots.
The Betti numbers for torus knots are due to
\cite{LS} (see also \cite{Pi, GM1}). We focus in this
paper on non-torus knots. 
\vskip 0.0cm

The series of the plane curve
singularities for
{\em Puiseux exponents\,} $(4,2u,v)$
for odd $u,v$ and $v>2u>4$ is the simplest of non-torus
type; the corresponding links are 
$C\!ab(2u+v,2)T(u,2)$. Here we generalize the formulas
from \cite{Pi} for the dimensions of cells
in the corresponding $CW$\~presentation. The most convincing
demonstrations of our Main Conjecture are the examples where
such cells are not all affine. Such examples are well beyond
\cite{Pi} and are actually a new vintage in the theory of
compactified Jacobians as well as our flagged generalization. 
\vskip 0.3cm

{\sf Some perspectives.}
An extension of the 
geometric approach to superpolynomials from this paper
to {\em all\,}
root systems is of obvious
interest, especially due to connections with $p$\~adic orbital 
integrals. 
Paper \cite{ChEl} hints that such a uniform
theory may exist, in spite of the fact that there can be
no rank stabilization for the systems $EFG$. 
Such a theory can be expected to provide {\em refined\,}
generalizations of {\em orbital integrals\,} from the geometric 
Fundamental Lemma; local spectral curves are taken here as 
plane curve singularities. See the end of the paper.  

\vskip 0.0cm
The case $q\!=\!1,a\!=\!0$ 
is directly related to $p$\~adic orbital integrals
of nil-elliptic type $A$. An immediate corollary 
of Conjecture \ref{CONJSUP}  is that such orbital integrals
are {\em topological\,}
invariants of the corresponding plane curve singularities.
This readily follows from \cite{LS} for torus knots, but seems
beyond any existing approaches for non-torus singularities, 
especially in the presence of non-affine cells (see online
version of this paper).
This invariance, the refined orbital integrals, 
the connections with HOMFLYPT homology and an
extension of our paper to arbitrary algebraic links, any
colors and all root systems are natural challenges. 
\vfill
\eject


\comment{
Summation formulas for counting
the cell dimensions and 
Betti numbers of the corresponding Jacobian factors 
for torus knots and this family 
were obtained in \cite{Pi}; for torus knots, they coincide 
with  those in \cite{LS} and match
``statistics" \cite{GM1}. 
Many examples of Betti numbers for torus knots and this series
are given in \cite{Pi}. However, no  {\em general\,} 
conjectural (combinatorial) formulas for Betti numbers
were suggested there; see Theorem 22 and Conjecture 23 in 
\cite{Pi}. The DAHA-Betti polynomials are actually 
quite computable ($q=1$ is the case of trivial center
charge).
Let us mention here \cite{Mel}, where DAHA
superpolynomials for torus knots at $q\!=\!1,a\!=\!0$ were 
connected with the combinatorial formulas from \cite{Gor1, Gor2}.

Importantly, the positivity was conjectured in
\cite{CJ,ChD1} for 
the {\em reduced\,} superpolynomials of algebraic 
knots colored by any {\em rectangle\,} Young diagrams,
i.e. not only in the uncolored case. 
We are going to incorporate such diagrams in our approach. 
Generally, adding colors requires using links, but 
we see some direct ways for such diagrams. 

An extension of the 
geometric approach to superpolynomials from this paper
to {\em all\,}
root systems (including exceptional ones). This of great
interest due to connections with $p$\~adic orbital integrals. 
Let us say something
about this exciting possibility. 
\smallskip

Paper \cite{ChEl} hints that such a uniform
theory may exist, in spite of the fact that there can be
no rank stabilization for the systems $EFG$. Note that this
part of \cite{ChEl} is of clear experimental nature. The rank 
stabilization is the main road to the superpolynomials so far 
(including the Khovanov-Rozansky theory and the DAHA approach). 

Such a uniform theory can be expected to provide
certain {\em refined  orbital integrals\,} from the geometric 
Fundamental Lemma; local spectral curves are taken here as 
plane curve singularities.   
The affine Springer fibers are related to 
our approach, though there are deviations at the level
of our definition of 
{\em flagged Jacobian factors\,}. See the 
end of the paper. This relation 
clearly demonstrates an importance
of extending our construction to any root systems.
\smallskip

To recapitulate, the significance of the general theory of plane 
curve singularities is quite obvious from any viewpoint.
The case of torus knots is special technically;  
superpolynomials and all related considerations 
are much simpler here due to the action of $\C^*$ on 
the corresponding singularities. For instance, the 
cells in the Lusztig-Smelt-Piontkowski-type decompositions
are always of type $\mathbb{A}^{\!N}$ (a general fact for any 
singularities with the action of $\C^*$, not only planar), 
no admissibility conditions from \cite{Pi} occur,
finding the dimensions of the cells becomes much less involved 
and so on. 

We think that the approach we present significantly clarifies
the geometric nature of (uncolored) superpolynomial, including
torus knots; it is actually quite transparent theoretically. 
Also, our methods are mostly elementary by nature; they 
extend those in \cite{Pi}. Potentially we can reach any
algebraic knots and fully compare our geometric approach
with the DAHA-based theory from \cite{ChD1}. 

}

\setcounter{section}{0}
\setcounter{equation}{0}
\section{\sc DAHA superpolynomials}
We will provide here the main facts of
DAHA theory needed for the definition
of the DAHA-Jones polynomials and DAHA
superpolynomials. See \cite{CJ,CJJ,C101}
for details. The construction is totally
uniform for any root systems and weights.

\subsection{\bf Definition of DAHA}
Let $R=\{\al\}   \subset \R^n$ be a root system of type
$A_n,\ldots,\!G_2$
with respect to a euclidean form $(z,z')$ on $\R^n
\ni z,z'$,
$W$ the Weyl group 
generated by the reflections $s_\al$,
$R_{+}$ the set of positive  roots
corresponding to fixed simple 
roots $\al_1,...,\al_n$. 
The form is normalized
by the condition  $(\al,\al)=2$ for 
{\em short\,} roots. 
The weight lattice is
$P=\oplus^n_{i=1}\Z \om_i$, 
where $\{\om_i\}$ are fundamental weights.
The root lattice is $Q=\oplus_{i=1}^n \Z\al_i$.
Replacing $\Z$ by $\Z_{+}=\{\Z\ni m\ge 0\}$, we obtain
$P_+,Q_+$ See  e.g., \cite{Bo} or \cite{C101}. 

Setting 
$\nu_\al\equal (\al,\al)/2$,
the vectors $\ \tal=[\al,\nu_\al j] \in
\R^n\times \R \subset \R^{n+1}$
for $\al \in R, j \in \Z $ form the
{\em twisted affine root system\,}
$\tR \supset R$ ($z\in \R^n$ are identified with $ [z,0]$).
We add $\al_0 \equal [-\vth,1]$ to the simple
roots for the {\em maximal short root\,} $\vth\in R_+$.
The corresponding set
$\tR_+$ of positive roots is 
$R_+\cup \{[\al,\nu_\al j],\ \al\in R, \ j > 0\}$.

The set of the indices of the images of $\al_0$ by all
automorphisms of the affine Dynkin diagram will be denoted by 
$O$ ($O=\{0\} \for E_8,F_4,G_2$). 
Let $O'\equal\{r\in O, r\neq 0\}$; $O'=[1,\ldots,n]$ for $A_n$.
The elements $\om_r$ for $r\in O'$ are  
{\em minuscule weights\,}, defined by the inequalities 
$(\om_r,\al^\vee)\le 1$ for all $\al \in R_+$. We set $\om_0=0$
for the sake of uniformity.
\smallskip

{\sf Affine Weyl groups.}
Given $\tal=[\al,\nu_\al j]\in \tR,  \ b \in P$, let
\begin{align}
&s_{\tal}(\tz)\ =\  \tz-(z,\al^\vee)\tal,\
\ b'(\tz)\ =\ [z,\ze-(z,b)]
\label{ondon}
\end{align}
for $\tz=[z,\ze] \in \R^{n+1}$.
The
{\em affine Weyl group\,} $\tW=\lan s_{\tal}, \tal\in \tR_+\ran$ 
is the semidirect product $W\lsmash Q$ of
its subgroups $W=$ $\lan s_\al,
\al \in R_+\ran$ and $Q$, where $\al$ is identified with
\begin{align*}
& s_{\al}s_{[\al,\,\nu_{\al}]}=\
s_{[-\al,\,\nu_\al]}s_{\al}\for
\al\in R \hbox{\, \, considered in\,} \, \lan s_{\tal}\ran. 
\end{align*}

Using the presentation of $\tW$ as  $W\lsmash Q$,
the {\em extended Weyl group\,} $ \hW$ can be defined as
$W\lsmash P$, where
the corresponding action is 
\begin{align}
&(wb)([z,\ze])\ =\ [w(z),\ze-(z,b)] \for w\in W,\, b\in P.
\label{ondthr}
\end{align}
It is canonically 
isomorphic to $\tW\lsmash \Pi$ for $\Pi\equal P/Q$. 
The latter group consists of $\pi_0=$id\, and the images $\pi_r$
of minuscule $\om_r$ in $P/Q$. 

The group $\Pi$
will be naturally identified with the subgroup of $\hW$ of the
elements of the length zero; the {\em length\, } is defined as 
follows:
\begin{align*}
&l(\hw)=|\la(\hw)| \for \la(\hw)\equal\tR_+\cap \hw^{-1}(-\tR_+).
\end{align*}
One has $\om_r=\pi_r u_r$ for $r\in O'$, where $u_r$ is the
(unique) element $u\in W$ of minimal length such 
that $u(\om_r)\in -P_+$. 
\smallskip

Setting $\hw = \pi_r\tw \in \hW$ for $\pi_r\in \Pi,\, \tw\in \tW,$
the length $l(\hw)$ coincides with the length of any reduced 
decomposition of $\tw$ in terms of the simple reflections
$s_i, 0\le i\le n$ (a standard and important fact). 
\medskip

Let $\mm,$ be the least natural number
such that  $(P,P)=(1/\mm)\Z;$ $\mm=n+1$ for $A_n$.
The double affine Hecke algebra, {\em DAHA\,}, depends
on the parameters
$q, t_\nu\, (\nu\in \{\nu_\al\});\,$ to be exact,
it is defined over the ring
of polynomials in terms of $q^{\pm 1/\mm}$ and
$\{t_\nu^{\pm1/2}\}.$ 

For $\tal=[\al,\nu_\al j] \in \tR$, we set
$t_{\tal}\!=\!t_{\al}\!=\!t_{\nu_\al}$,
$q_{\tal}\!=\!q_{\al}\!=\!q^{\nu_\al}$, 
and introduce $k_{\tal}\!=\!k_\al\!=\!k_{\nu_\al}$ 
from the relation $t_{\nu}=q^{\nu k_\nu}$. For
$i=1,\ldots, n$, let
\begin{align*}
\rho_k\equal \frac{1}{2}\!\sum_{\al>0} k_\al \al=
k_{\sht}\rho_{\sht}\!+\!k_{\lng}\rho_{\lng},\ \,
\rho_\nu=\frac{1}{2}\!\sum_{\nu_\al=\nu} \al=
\!\!\sum_{\nu_{\al_i}=\nu}  \om_i,
\end{align*}
where {\small sht,\, lng\,} are used for short and
long roots. We note that the specialization
$k_{\sht}\!=\!1\!=\!k_{\lng}$ corresponds to quantum groups
and provides the {\em WRT invariants\,} in the construction below;
see \cite{CJ}.
\smallskip

For pairwise commutative $X_1,\ldots,X_n,$
\begin{align}
& X_{\tb}\ \equal\ \prod_{i=1}^nX_i^{l_i} q^{ j}
\iif \tb=[b,j],\ \hw(X_{\tb})\ =\ X_{\hw(\tb)},
\label{Xdex}\\
&\hbox{where\ \,} b=\sum_{i=1}^n l_i \om_i\in P,\ \,\, j \in
(1/ \mm)\,\Z,\ \,\, \hw\in \hW.
\notag \end{align}
For instance, $X_0\equal X_{\al_0}=X_{[-\vth,1]}=qX_\vth^{-1}$.
\medskip

Recall that 
$\om_r=\pi_r u_r$ for $r\in O'$ (see above). Note that
$\pi_r^{-1}$ is $\pi_{\iota(i)}$, where
$\iota$ is the standard involution  of the {\em nonaffine }
Dynkin diagram,
induced by $\al_i\mapsto -w_0(\al_i)$; it is the reflection
of $[1,\ldots,n]$ in type $A_n$.
Finally, we set $m_{ij}=2,3,4,6$
when the number of links between $\al_i$ and $\al_j$ in the affine 
Dynkin diagram is $0,1,2,3$.

\begin{definition}
The double affine Hecke algebra $\HH\ $
is generated by
the elements $\{ T_i,\ 0\le i\le n\}$,
pairwise commutative $\{X_b, \ b\in P\}$ satisfying
(\ref{Xdex})
and the group $\Pi,$ where the following relations are imposed:

(o)\ \  $ (T_i-t_i^{1/2})(T_i+t_i^{-1/2})\ =\
0,\ 0\ \le\ i\ \le\ n$;

(i)\ \ \ $ T_iT_jT_i...\ =\ T_jT_iT_j...,\ m_{ij}$
factors on each side;

(ii)\ \   $ \pi_rT_i\pi_r^{-1}\ =\ T_j \iif
\pi_r(\al_i)=\al_j$;

(iii)\  $T_iX_b \ =\ X_b X_{\al_i}^{-1} T_i^{-1} \iif
(b,\al^\vee_i)=1,\
0 \le i\le  n$;

(iv)\ $T_iX_b\ =\ X_b T_i\ $ if $\ (b,\al^\vee_i)=0
\for 0 \le i\le  n$;

(v)\ \ $\pi_rX_b \pi_r^{-1}\ =\ X_{\pi_r(b)}\ =\
X_{ u^{-1}_r(b)}
 q^{(\om_{\iota(r)},b)},\  r\in O'$.
\label{double}
\end{definition}

Given $\tw \in \tW, r\in O,\ $ the product
\begin{align}
&T_{\pi_r\tw}\equal \pi_r T_{i_l}\cdots T_{i_1},\where
\tw=s_{i_l}\cdots s_{i_1} \for l=l(\tw),
\label{Twx}
\end{align}
does not depend on the choice of the reduced decomposition.
Moreover,
\begin{align}
&T_{\hv}T_{\hw}\ =\ T_{\hv\hw}\  \hbox{ whenever\,}\
 l(\hv\hw)=l(\hv)+l(\hw) \for
\hv,\hw \in \hW. \label{TTx}
\end{align}
In particular, we arrive at the pairwise
commutative elements 
\begin{align}
& Y_{b}\equal
\prod_{i=1}^nY_i^{l_i} \iif
b=\sum_{i=1}^n l_i\om_i\in P,\ 
Y_i\equal T_{\om_i},b\in P.
\label{Ybx}
\end{align}


\subsection{\bf Main features}\label{sect:Aut}
The following maps can be (uniquely) extended to
automorphisms of $\HH\,$, where                    
$q^{1/(2\mm)}$ must be added to the ring of constants
(see \cite{C101}, (3.2.10)-(3.2.15)):
\begin{align}\label{tauplus}
& \tau_+:\  X_b \mapsto X_b, \ T_i\mapsto T_i\, (i>0),\
\ Y_r \mapsto X_rY_r q^{-\frac{(\om_r,\om_r)}{2}}\,,
\\
& \tau_+:\ T_0\mapsto  q^{-1}\,X_\vth T_0^{-1},\
\pi_r \mapsto q^{-\frac{(\om_r,\om_r)}{2}}X_r\pi_r\
(r\in O'),\notag\\
& \label{taumin}
\tau_-:\ Y_b \mapsto \,Y_b, \ T_i\mapsto T_i\, (i\ge 0),\
\ X_r \mapsto Y_r X_r q^\frac{(\om_r,\om_r)}{ 2},\\
&\tau_-(X_{\vth})= 
q T_0 X_\vth^{-1} T_{s_{\vth}}^{-1};\ \
\si\equal \tau_+\tau_-^{-1}\tau_+\, =\,
\tau_-^{-1}\tau_+\tau_-^{-1}.\notag
\end{align}
These automorphisms fix $\ t_\nu,\ q$
and their fractional powers.

The span 
of $\tau_{\pm}$ is the {\em projective $PSL_{2}(\Z)\,$}
(due to Steinberg), which is isomorphic to the braid 
group $B_3$. Let us list the matrices corresponding to the 
automorphisms above upon the natural projection 
onto $SL_2(\Z)$, which is upon the specialization  
$\,t^{\frac{1}{2\mm}}_\nu,q^{\frac{1}{2\mm}}\mapsto 1$. 
The matrix {\tiny 
$\begin{pmatrix} \al & \be \\ \ga & \de\\ \end{pmatrix}$}
will represent the map $X_b\mapsto X_b^\al Y_b^\ga,
Y_b\mapsto X_b^\be Y_b^\de$ for $b\in P$. One has:\ 
$\tau_+\mapsto$ 
{\tiny 
$\begin{pmatrix}1 & 1 \\0 & 1 \\ \end{pmatrix}$},\ 
$\tau_-\mapsto$ 
{\tiny 
$\begin{pmatrix}1 & 0 \\1 & 1 \\ \end{pmatrix}$},\
$\si\mapsto$ 
{\tiny 
$\begin{pmatrix}0 & 1 \\-1 & 0 \\ \end{pmatrix}$}.
\smallskip

We note that there are some simplifications with 
the definition of DAHA and $\tau_{\pm}$ for $A_n$ 
and in part $(i)$ of Theorem \ref{JONITER},    
but they are not significant (the theory is very
much uniform for any root systems). However $A_n$
is obviously needed in part $(ii)$
in this theorem. 
\smallskip 


Following \cite{C101}, 
we use the PBW Theorem to express any $H\in \HH$ in the form 
\,$\sum_{a,w,b} c_{a,w,b}\, X_a T_{w} Y_b$\, for $w\in W$,
$a,b\in P$ (this presentation is unique). Then we substitute:
\begin{align}\label{evfunct}
\{\,\}_{ev}:\ X_a\  \mapsto\  X_a(q^{-\rho_k})\!=
\!q^{-(\rho_k,a)},\ 
Y_b \ \mapsto\  q^{(\rho_k,b)},\ 
T_i \ \mapsto\  t_i^{1/2}. 
\end{align}

The functional $\,\HH\ni H\mapsto \{H\}_{ev}$, 
called {\em coinvariant\,}, acts via the projection 
$H\mapsto H\!\Downarrow \equal H(1)$ of $\HH\,$
onto the {\em polynomial representation \,}$\v$, which is
the $\HH$\~module induced from the one-dimensional
character $T_i(1)=t_i^{-1/2}=Y_i(1)$ for $1\le i\le n$ and
$T_0(1)=t_0^{-1/2}$. Here $t_0=t_{\sht}$; 
see \cite{C101,CJ,CJJ}.
\smallskip

The polynomial representation $\v$
is linearly generated by $X_b (b\in P)$
and the action of $T_i(0\le i\le n)$ there 
is given by 
the {\em Demazure-Lusztig operators\,}:
\begin{align}
&T_i\  = \  t_i^{1/2} s_i\ +\
(t_i^{1/2}-t_i^{-1/2})(X_{\al_i}-1)^{-1}(s_i-1),
\ 0\le i\le n.
\label{Demazx}
\end{align}
The elements $X_b$ become the multiplication operators 
and  $\pi_r (r\in O')$ act via the general formula
$\hw(X_b)=X_{\hw(b)}$ for $\hw\in \hW$. 

{\sf Macdonald polynomials.}
The Macdonald polynomials $P_b(X)$ for 
$b\in P_+$ are uniquely defined
as follows. For $c\in P$, let $c_+$ be a unique element 
such that $c_+\in W(c)\cap P_+$. Given $b\in P_+$ and
assuming that $c\in P$ is such that 
$b\neq c_+\in b-Q_+$,
\begin{align}
&P_b\! -\!\!\!\!\sum_{a\in W(b)}\!\!\! X_{a} 
\in\, \oplus_{c}\Q(q,t_\nu) X_c 
\hbox{\, and\, }
CT\bigl(P_b X_{c^\iota}\,\mu(X;q,t)\bigr)\!=\!0, \label{macdsym}
\\
&\hbox{where\,\, } 
\mu(X;q,t)\equal\!\prod_{\al \in R_+}
\prod_{j=0}^\infty \frac{(1\!-\!X_\al q_\al^{j})
(1\!-\!X_\al^{-1}q_\al^{j+1})
}{
(1\!-\!X_\al t_\al q_\al^{j})
(1\!-\!X_\al^{-1}t_\al^{}q_\al^{j+1})}\,. \notag
\end{align}
Here $CT$ is the constant term; 
$\mu$ is considered
a Laurent series of $X_b$ with 
the coefficients expanded in terms of
positive powers of $q$. The coefficients of
$P_b$ belong to the field $\Q(q,t_\nu)$.
The following evaluation formula (the Macdonald
Evaluation Conjecture) is important
to us:
\begin{align}
\label{macdeval}
&(P_{b}(q^{-\rho_k}))=
q^{-(\rho_k,b)}
\prod_{\al>0}\prod_{j=0}^{(\al^{\!\vee},b)-1}
\Bigl(
\frac{
1- q_\al^{j}t_\al X_\al(q^{\rho_k})
 }{
1- q_\al^{j}X_\al(q^{\rho_k})
}
\Bigr).
\end{align}

\subsection{\bf Algebraic knots}\label{sec:ALG-KNOTS}
Torus knots $T(\rr,\ss)$ are defined for any integers 
$\rr,\ss>0$ such that \,gcd$(\rr,\ss)=1$. One has the 
symmetry $\,T(\rr,\ss)=T(\ss,\rr)$, where we use 
``$=$" for the ambient isotopy equivalence. Also 
$\,T(\rr,1)=\unknot\,$ for  
the {\em unknot\,}, denoted by $\,\unknot$. 

{\em Algebraic knots\,}
$\t(\vec{\rr},\vec{\ss})$
are associated with two sequences of
(strictly) positive integers: 
\begin{align}\label{iterrss} 
\vec{\rr}=\{\rr_1,\ldots \rr_\ell\}, \ 
\vec{\ss}=\{\ss_1,\ldots \ss_\ell\} \hbox{\, such that\, 
gcd}(\rr_i,\ss_i)=1;
\end{align}
$\ell$ will be called the {\em length\,} of $\,\vec\rr,\vec\ss$.
The pairs $\{\rr_i,\ss_i\}$ are 
{\em characteristic\,} or {\em Newton pairs\,}. 

We will need one more sequence: 
\begin{align}\label{Newtonpair}
\aa_1=\ss_1,\,\aa_{i}=\aa_{i-1}\rr_{i-1}\rr_{i}+\ss_{i}\,\ 
(i=2,\ldots,m). 
\end{align}
See  e.g. \cite{EN} and \cite{Pi}. Then, 
\begin{align}\label{Knotsiter}
\t(\vec{\rr},\!\vec{\ss})\equal
C\!ab(\vec{\aa},\!\vec{\rr})(\unknot)=
\bigl(C\!ab(\aa_\ell,\rr_\ell)\cdots C\!ab(\aa_2,\rr_2)\bigr)
\bigl(T(\aa_1,\!\rr_1)\bigr)
\end{align}
in terms of the {\em cabling} defined below.
Note that the first iteration (application of 
$C\!ab$) is for $\{\aa_1,\rr_1\}$ (not for the last pair!).
\smallskip

{\sf Cabling}.
The {\em cabling\,} $C\!ab(\aa,\bb)(K)$ of any oriented
knot $K$ in (oriented) $S^3$ is defined as follows; 
see e.g. \cite{Mo,EN} and references therein. 
We consider a small $2$\~dimensional torus
around $K$ and put there the torus knot $T(\aa,\bb)$
in the direction of $K$, 
which is $C\!ab(\aa,\bb)(K)$ (up to ambient isotopy).

This procedure depends on 
the order of $\aa,\bb$ and the orientation of
$K$. We choose the latter in the standard way (compatible 
with almost all sources, including  the Mathematica package
``KnotTheory"); the parameter $\aa$ gives the number of
turns around $K$. This construction also depends
on the framing of the cable knots; we take the natural
one, associated with the parallel copy of the torus where 
a given cable knot sits (its parallel copy 
has zero linking number with this knot). 

By construction,
$C\!ab(\aa,0)(K)=\unknot$ and  $C\!ab(\aa,1)(K)=K$ for any
knot $K$ and $\aa\neq 0$. See \cite{ChD1} for further discussion
of relations.
The pairs $\{\rr_i,\aa_i\}$ are sometimes called 
{\em topological};  the isotopy equivalence
of algebraic knots generally can be
seen only at the level of $\rr,\aa$\~parameters (not at the level 
of the Newton or Puiseux pairs).
\smallskip

{\sf Newton-Puiseux theory.}
Given a sequence $\rr_i,\ss_i>0$ of Newton (characteristic)
pairs the knot $\t(\vec\rr,\vec\ss)$ is the link 
of the {\em germ of the singularity\,} 
\begin{align}\label{yxcurve}
y = x^{\frac{\ss_1}{\rr_1}}
(c_1 + x^{\frac{\ss_2}{\rr_1\rr_2}}
\bigl(c_2 + \ldots + 
x^{\frac{\ss_\ell}{\rr_1\rr_2\cdots\rr_\ell}}\bigr))
\hbox{\, at\, } 0,
\end{align} 
which is the intersection of the corresponding plane curve 
with a small $3$-dimensional sphere in 
$\C^2$ around $0$. We will always
assume in this paper that this germ is {\em unibranch}.

The inequality  $\ss_1<\rr_1$ is commonly
imposed here (otherwise $x$ and $y$ can be switched).
Formula (\ref{yxcurve}) is the celebrated 
{\em Newton-Puiseux expansion}. See e.g.
\cite{EN}. All algebraic knots can be obtained
in such a way.

{\sf Jacobian factors.}
One can associate with a unibranch $\,\c_{\vec\rr,\vec\ss}\,$ the
{\em Jacobian factor\,} $J(\c_{\vec\rr,\vec\ss})$.
Up to a homeomorphism, it can be introduced as the canonical 
compactification of the {\em generalized Jacobian\,} of an
integral {\em rational\,} planar curve with 
$\,\c_{\vec\rr,\vec\ss}\,$
as its {\em only\,} singularity. It has a purely
local definition, which we will use below. Its dimension is 
the {\em $\de$\~invariant\,} of the singularity
$\,\c_{\vec\rr,\vec\ss},$ also called the
{\em arithmetic genus}.

Calculating the {\em Euler number\,}  
$e(J(\c_{\vec\rr,\vec\ss}))$, the topological Euler
characteristic of $J(\c_{\vec\rr,\vec\ss}),$
and the corresponding {\em Betti numbers\,} 
in terms of $\,\vec\rr,\vec\ss\,$ is a challenging problem. 
For torus knots $T(\rr,\ss)$, one has 
$e(J(\c_{\rr,\ss}))=$ $\frac{1}{\rr+\ss}
\binom{\rr+\ss}{\rr}$ due to \cite{Bea}. This formula 
is related to the perfect modules 
of rational DAHA and the combinatorics of
generalized Catalan numbers; see e.g. \cite{GM1}. 

\smallskip
The Euler numbers of $J(\c_{\vec\rr,\vec\ss})$ were calculated 
in \cite{Pi} (the Main Theorem) for the following triples of 
{\em Puiseux characteristic exponents\,}: 
\begin{align}\label{PiEuler}
(4,2u,v),\, (6,8,v)\, \hbox{\, and\, } 
(6,10,v)  \hbox{\, for odd\, } u,v>0,
\end{align}
where $4\!<\!2u\!<v$, and $8\!<\!v, 10\!<\!v$ respectively.
Here $\de\!=$dim\,$J(\c_{\vec\rr,\vec\ss})$
is $\frac{(\rr-1)(\ss-1)}{2}$ for $T(\rr,\ss)$ and 
$2u+(v-1)/2-1$ for the series $(4,2u,v)$. Generally,
$\de$ equals the cardinality
$|\N\setminus \Gamma|$, where $\Gamma$ is the {\em valuation
semigroup\,} associated with $\c_{\vec\rr,\vec\ss}$;
see \cite{Pi} and \cite{Za}. The Euler numbers of
the Jacobian factors can be also calculated via the 
HOMFLYPT polynomials of 
the corresponding links (see below) due to \cite{ObS,Ma}.
\smallskip

Concerning the Betti numbers for the torus knots and
the series $(4,2u,v)$, the odd (co)homology 
of $\,J(\c_{\vec\rr,\vec\ss})\,$ vanishes.
The formulas for the corresponding even Betti numbers 
$\,h^{(2k)}=\hbox{rk\,}(H^{2k}(J(\c_{\vec\rr,\vec\ss})))$
were calculated explicitly for many 
values of $k$ in \cite{Pi}, where $0\le k\le \de$. Not 
much was and is known/expected beyond these two series. 

\subsection{\bf DAHA-Jones theory}
The following results and conjectures are mainly 
from \cite{ChD1}; see also  
Theorem 1.2 from \cite{CJJ} and \cite{CJ,GoN}. 

The construction is given directly in terms of
the parameters $\{\vec{\rr},\vec{\ss}\}$, though it
actually depends only on the corresponding topological 
parameters $\{\vec{\aa},\vec{\rr}\}$. Recall that 
torus knots $T(\rr,\ss)$ are naturally represented   
by $\ga_{\rr,\ss}\in PSL_{2}(\Z)$ with the 
first column $(\rr,\ss)^{tr}$ ($tr$ denotes transposition),
where $\,\rr,\ss>0$ and we 
assume that \,gcd$(\rr,\ss)=1$. 
Let $\hat{\ga}_{\rr,\ss}$ be
any pullback of $\ga_{\rr,\ss}$ to the projective $PSL_2(\Z)$.
\smallskip

For a polynomial $F$ in terms of 
fractional powers of $q$ and $t_\nu$, 
the {\em tilde-normalization}
$\tilde{F}$ will be the result of the division of $F$
by the lowest $q,t_\nu$\~monomial, assuming that it 
is well defined. We put $q^\bullet t^\bullet$ for
a monomial factor (possibly fractional)
in terms of $q,t_\nu$. See \cite{ChD1} for the
following theorem. We will also apply this definition
to the superpolynomials, where the lowest $q,t$\~monomial
is picked from the constant $a$\~term.

\begin{theorem}\label{JONITER}
Let $R$ be a reduced irreducible root system.
Recall that $H\mapsto H\!\!\Downarrow\, \equal H(1)$, where the
action of $H\in \HH$ in $\v$ is used.

(i) Given two strictly
positive sequences $\vec\rr,\vec\ss$ of length $\ell$ as 
in (\ref{iterrss}), we lift $(\rr_i,\ss_i)^{tr}$ to
$\ga_i$ and then to $\hat{\ga_i}$ (acting in $\HH$) as above.
For a weight $b\in P_+\,$, the {\sf DAHA-Jones polynomial} is
\begin{align}\label{jones-dit}
&J\!D_{\,\vec\rr,\,\vec\ss\,}^{R}(b\,;\,q,t)\ =\ 
J\!D_{\,\vec\rr,\,\vec\ss\,}(b\,;\,q,t)\, \equal\\ 
\Bigl\{\hat{\ga_{1}}&\Bigr(
\cdots\Bigl(\hat{\ga}_{\ell-1}
\Bigl(\bigl(\hat{\ga}_\ell(P_b)/P_b(q^{-\rho_k})\bigr)\!\Downarrow
\Bigr)\!\Downarrow\Bigr) \cdots\Bigr)\Bigr\}_{ev}.\notag
\end{align}
It does not depend
on the particular choice of the lifts \,$\ga_i$
and $\hat{\ga}_i\in GL_{\,2}^{\wedge}(\Z)$.
The {\sf tilde-normalization} 
$\tilde{J\!D}_{\vec\rr,\vec\ss}\,(b\,;\,q,t)$ is well defined
and is a polynomial in terms of $\,q,t_\nu$ with the
constant term $1$. 

(ii) Let us switch to the root system $A_n$ for $sl_{n+1}$,
setting $t\!=\!t_{\sht}\!=\!q^k$ and  considering $P_+\ni b=$ 
$\sum_{i=1}^n b_i \om_i$ as 
(dominant) weights for any $A_m$ (for $sl_{m+1}$)
with $m\ge n-1$, where we assume that $\om_{n}=0$ upon the 
restriction to $A_{n-1}$.

Then given $\t(\vec\rr,\vec\ss)$ as above, there exists a 
{\sf DAHA\~superpolynomial\,} 
$\h_{\,\vec\rr,\,\vec\ss}\,(b\,;\,q,t,a)$
 from $\Z[q,t^{\pm 1},a]$ satisfying the relations
\begin{align}\label{jones-sup}
\h_{\,\vec\rr,\,\vec\ss}\,(b;q,t,a\!=\!-t^{m+1})\!=\!
\tilde{J\!D}_{\,\vec\rr,\,\vec\ss\,}^{A_{m}}(b;q,t)
\hbox{\, for any } m\!\ge\! n\!-\!1;
\end{align} 
then its $a$\~constant term 
is automatically tilde-normalized.
\sq
\end{theorem}

\comment{
\begin{theorem}\label{STABILIZ}
We switch to the root system $A_n$ for $sl_{n+1}$,
setting $t=t_{\sht}=q^k$. 
Let us consider $P_+\ni b=$ $\sum_{i=1}^n b_i \om_i$ as 
a (dominant) weight for any $A_m$ (for $sl_{m+1}$)
with $m\ge n-1$,
where we set $\om_{n}=0$ upon the restriction to $A_{n-1}$.

(i) {\sf Stabilization.} Given $\t(\vec\rr,\vec\ss)$, there exists 
$\h_{\,\vec\rr,\,\vec\ss}\,(b\,;\,q,t,a)$,
a polynomial from $\Z[q,t^{\pm 1},a]$, satisfying the relations
\begin{align}\label{jones-sup}
\h_{\,\vec\rr,\,\vec\ss}\,(b;q,t,a\!=\!-t^{m+1})\!=\!
\tilde{J\!D}_{\,\vec\rr,\,\vec\ss\,}^{A_{m}}(b;q,t)
\hbox{\, for any } m\!\ge\! n\!-\!1.
\end{align} 

(ii) {\sf Specialization.}
Making here $q=1$, one has for any weights and
iterated knots: 
\begin{align}\label{jones-seval}
\h_{\,\vec\rr,\,\vec\ss}\,(b\,;\,q\!=\!1,t,a)\!=\!
\prod_{i=1}^n\h_{\,\vec\rr,\,\vec\ss}\,
(\om_i\,;\,q\!=\!1,t,a)^{b_i}
\hbox{\, for\, } b\!=\!\sum_{i=1}^n b_i \om_i.
\end{align}

(iii) {\sf The $a$-degree.}
Assuming  that $\rr_1>\ss_1$,                                  
deg${}_a \h_{\,\vec\rr,\,\vec\ss}\,(b\,;\,q,t,a)\,$
(the $a$-degree) is no greater than                  
$\ss_1\rr_2\cdots\rr_\ell-1$         
times \, ord $(\la_b)$, the 
number of boxes in the Young diagram $\la_{b}$
associated with $b\in P_+$. This is the multiplicity of the 
corresponding singularity minus $1$ \cite{Rego} in the uncolored 
case. (Conjecturally the equality always holds.)

(iv) {\sf Super-duality.} 
Up to a power of $q$ and $t$,
\begin{align}\label{iter-duality}
\h_{\,\vec\rr,\,\vec\ss}\,(\la\,;q,t,a)=q^{\bullet}t^{\bullet}
\h_{\,\vec\rr,\,\vec\ss}\,(\la^{tr}\,;t^{-1},q^{-1},a),
\end{align}
where we switch from (arbitrary) dominant
weights $\,b\,$ to the corresponding Young diagrams $\la=\la_b$,
and $\la^{tr}$ in the right-hand side is the transposition of 
$\la$. 
\end{theorem} 
\vskip -1.25cm \sq
\vskip 0.2cm

Concerning the stabilization and duality, we mainly follow  
\cite{GoN,ChD2}, switching from $P_{\la}$
to the {\em modified Macdonald polynomials\,}. See also
\cite{GS,CJ,CJJ}; there is a sketch of the proof in the 
latter based on the theory at roots of unity (potentially 
not only in type $A$).  
Note that the DAHA of type $gl_{n+1}$ and
that for $sl_{n+1}$ (of type $A_n$) result in {\em coinciding
DAHA-Jones polynomials\,} up to $q^\bullet t^\bullet$, which is 
not difficult to justify. The compatibility
of DAHA and the Macdonald polynomials of types $gl_{n+1}$ and 
$sl_{n+1}$ is used here.
 
}

{\sf Topological connection.}
Let us briefly discuss the conjectural relation
of DAHA\~super\-polynomials to 
{\em stable Khovanov-Rozansky polynomials\,}
denoted by  $K\!h\!R_{\hbox{\tiny stab}}$.
See \cite{KhR1,KhR2,Kh,Ras}. We consider only the
reduced setting (actually not quite developed
topologically). 

The passage to the Khovanov-Rozansky polynomials
for $sl_N$ for sufficiently large $N$   
is the substitution
$a\mapsto\! t^{N}\sqrt{q/t}$. Note the relation 
to the {\em Heegaard-Floer homology\,} for $N=0$.
Equivalently, this passage is
$a_{st}\mapsto q_{st}^{N}$ in the {\em standard topological
parameters} (also used in the {\em ORS Conjecture\,}),
which are related to the DAHA parameters as follows:
\begin{align}\label{qtareli}
&t=q^2_{st},\  q=(q_{st}t_{st})^2,\  a=a_{st}^2 t_{st},\notag\\
&q^2_{st}=t,\  t_{st}=\sqrt{q/t},\  a_{st}^2=a\sqrt{t/q}.
\end{align}
\smallskip

For the DAHA\~superpolynomials from
Theorem \ref{jones-sup}, 
\begin{align}\label{khrconj}
\h_{\,\vec\rr,\,\vec\ss}\,(\square\,;\,q,t,a)_{st}=
\tilde{K\!h\!R}_{\hbox{\tiny stab}}(q_{st},t_{st},a_{st})\
\hbox{\, where\ \,} \square=\omega_1,
\end{align}
and $\tilde{K\!h\!R}_{\hbox{\tiny stab}}$ is {\em reduced\,}
$K\!h\!R_{\hbox{\tiny stab}}$ divided  
by the smallest power of $a_{st}$ and then by
$q_{st}^\bullet t_{st}^\bullet\,$ such that
$\,\tilde{K\!h\!R}_{\hbox{\tiny stab}}(a_{st}\!=\!0)
\in \Z_+[q_{st},t_{st}]$ with the constant term $1$.
Here $\{\cdot\}_{st}$ means the switch from the DAHA
parameters to the {\em standard topological parameters}.

\smallskip
Also, the polynomials $K\!h\!R_{\hbox{\tiny stab}}$
are expected to coincide with the
(reduced) physics {\em superpolynomials}  
based on the {\em BPS states\,} \cite{DGR,AS,FGS,GGS}   
and those obtained in terms of {\em rational DAHA\,}
\cite{GORS,GoN} for torus knots. The latter approach is
developed so far only for torus knots and 
in the uncolored case; there is some progress for
symmetric powers of the fundamental representation  
(see \cite{GGS}). We will not touch the connections
with rational DAHA in this paper. Concerning physics
origins, let us mention that using the Macdonald 
polynomials at roots of unity $q$ (for
$t\!=\!q^k, k\in \Z_+$) instead of Schur functions in the 
usual construction of {\em knot opertators\,} was suggested 
in \cite{AS}.
\medskip

{\sf Betti numbers 
$h^{(i)}\!\!=\!\hbox{rk\,}H^{i}(J(\c_{\vec\rr,\vec\ss}))$.}
We are very much motivated by the 
DAHA approach to these numbers. 
Technically, we generalize the interpretation
of Gorsky's superpolynomials for torus knots at $a\!=\!0$
from \cite{GM1,GM2} and the following conjecture 
from \cite{ChD1}:
  
\comment{
Following Section \ref{sec:ALG-KNOTS}, let
$\,\c_{\vec\rr,\vec\ss}$ be the (unibranch) germ
at $0$ of the plain curve singularity corresponding
to $\{\rr_i,\ss_i>0\}$ and $\,h^{(k)}\,
=\hbox{rk\,}(H^{k}(J(\c_{\vec\rr,\vec\ss})))$
the Betti numbers of 
the Jacobian factor $\,J(\c_{\vec\rr,\vec\ss})$. Then
Conjecture 2 from \cite{ORS} 
and (conjectural) relation (\ref{khrconj}) formally 
result in
the coincidence of the tilde-normalization of {\em reduced\,}
$\mathscr{P}_{\! \hbox{\tiny alg}}$ there
with our
$\h_{\,\vec\rr,\,\vec\ss}\,(\square\,;\,q,t,a)$.

See (4.37) from \cite{ChD1} and Section 4
from \cite{ORS} for the relation between reduced and
unreduced (the main in \cite{ORS}) setting. The reduced
setting generally relies on Jacobian factors vs. Hilbert
schemes; see below and Section 9.1 from \cite{GORS}.
}

\begin{align}\label{bettit}
\h_{\,\vec\rr,\,\vec\ss}\,(\square\,;\,q=1,t,a=0)=
\hbox{\small$\sum$}_{i=0\,}^{2\de} h^{(i)}t^{i/2} \hbox{\, for\, }
\de\!=\!\hbox{dim}\,J(\c_{\vec\rr,\vec\ss}).
\end{align}
It implies that $h^{\hbox{\tiny odd}}=0$ (the van Straten- Warmt 
conjecture). Relation (\ref{bettit}) will be generalized below 
to the whole superpolynomial 
$\h_{\,\vec\rr,\,\vec\ss}\,(\square\,;\,q,t,a)$, which is the main
result of our paper.
\smallskip

\setcounter{equation}{0}
\section{\sc Geometric superpolynomials}
\subsection{\bf Modules of semigroups} 
Let $\r$ be the complete local ring of the {\em unibranch\,}
germ of the plane curve singularity, embedded into the 
normalization ring $\o=\C[[z]]$. The {\em conductor\,}
of $\r$ is the smallest $\cc$ such that $z^{\cc}\o\subset \r$;
actually it is the ideal $(z^{\cc})\subset \o$, but we will 
call $\cc$ the conductor in this paper. We set $\k\equal\C((z))$.

The corresponding {\em semigroup\,} $\Ga_{\r}$ is formed by the
orders of the smallest powers, i.e. {\em valuations\,} 
$\nu(x)$ (minimal $z$\~degrees) $x\in\r$. The 
$\de$\~invariant (the arithmetic genus) is
then $\de_{\r}=|\Z_+\setminus \Ga_{\r}|$. We will call 
$\Z_+\setminus \Ga_{\r}$ the {\em set of gaps} and  
denote it $G_{\r}$; thus $\de=|G_{\r}|$. Also, $\cc=2\de$.
 
Compactified Jacobians for projective curves
 are generally defined as the varieties of coherent
torsion free sheaves of rank one and fixed  degree 
up to isomorphisms. The {\em Jacobian factor} $J_{\r}$,
we are going to define, is a 
local version of the compactified Jacobian.
It is (as a set) formed by all finitely generated 
$\r$\~submodules $\m\subset \k=\C((z))$, 
of (any) prescribed degree, also called
$\r$\~{\em lattices}.

We define {\em $\o$\~degree\,} deg$_{\o}\m$ with respect 
to $\o$; 
it is dim$_{\C}(\o/\m)$ if $\m\subset\o$. For
arbitrary submodules $\m\subset \k$:
\begin{align}\label{deg-def}
\hbox{deg}_\o(\m)=\hbox{dim}_{\C}(\o/(\o\cap\m))-
\hbox{dim}_{\C}(\m/(\o\cap\m)).
\end{align}
This definition is a natural counterpart 
of the degree of a divisor
at a given point (here at $z=0$) in the smooth situation.
Actually we will mainly need below 
$-$deg$_{\r}(\m)=\de-\hbox{deg}_\o(\m)$. 
\smallskip

The valuations $\nu(x)$ of the elements $x\in\m\subset \k$ form
a $\Ga_{\r}$-module $\De_{\m}$; the modules for semigroups
$\Ga$ (with $0$) are subsets $\De\subset \Z_+$
such that $\Ga+\De\subset\De$. 
Unless stated otherwise, we assume that $\m\subset \o$ and
that it contains the element 
$1+\sum_{i=1}^{\cc\!-\!1}\la_0^i z^i$
(of valuation $0$), such an embedding can be achieved by the
division by $z^{m}$ for $m=\min(\De_{\m})$. Here the 
upper limit $\cc\!-\!1$
is sufficient in the sum due to the definition of the 
conductor. The notation $\De_\circ$ is used for such a 
normalization in \cite{Pi},
we call it the {\em standard normalization}. See 
there for these  and related definitions and facts. 

For a {\em standard\,} $\m$, 
we will use the notation $D_{\m}$ or $D[\m]$ for 
$G_{\r}\cap \De_{\m}$ 
and call it the {\em set of added gaps\,}
or simply the {\em $D$\~set\,}.
The square brackets will be used
for the list of its elements. For instance,  
$D_{\r}\!=\!\varnothing\!=\![\,]$ corresponds to 
the (trivial) $\Ga$\~module
$\De_{\r}=\Ga$, and  
$D_{\o}\!=\!G_{\r}\!=\![1,\ldots, \cc\!-\!1]$ for $\m\!=\!\o$
(recall that $\cc\!-\!1$ is the last gap in $G$). 

Due to the normalization we impose, one has:
\begin{align}\label{deg-fomula}
\hbox{deg}_\o(\m)=\de-|D_{\m}|,\ \,
 -\hbox{deg}_{\r}(\m)=|D_{\m}|
\hbox{\, for standard\, } \m .
\end{align} 

Not all $\Ga$\~modules
$\De$ can be realized as $\De_{\m}$ for non-torus
singularities. Recall that torus knots $T(\rr,\ss)$
are associated with the rings
$\r=\C[[z^{\rr},z^{\ss}]]$. The simplest example
of a non-torus singularity is $\r=\C[[z^4,z^6+z^7]]$
with $\Ga=\lan 4,6,13\ran$ and $\cc=16$.
Then the sets of added gaps $D=[2,15]$ and $D=[2,11,15]\,$ do 
not come from any modules $\m$; following \cite{Pi}, we call
them {\em non-admissible\,}. All other $D$ are admissible
in this case (i.e. can be obtained as $D_{\m}$).
\smallskip

\subsection{\bf 
\texorpdfstring{{\mathversion{bold}$J_{\r}$}}{J}
as a projective variety}
Let $J_{\r}[D]\!\subset\! J_{\r}$ be the set of modules $\m$ with 
$D_{\m}\!=\!D$, $J_{\r}[d]\!=\!\cup_{\,|D|= d\,} J_r[D]$ and
$\bar{J}_{\r}[d]\!=\!\cup_{\,d'\ge d\,} J_r[d']$ for $d\ge 0$. I.e.
the latter is the set of all standard submodules of $\o$\~degree
$\de\!-\!d$ or smaller.
The set $J_{\r}[0]=J_{\r}[\varnothing]$ is {\em the big cell\,}; 
it is formed by all invertible modules $\m$ (with one generator
and of $\o$\~degree $\de$). 
Also, $J_{\r}[\de]=\bar{J}_{\r}[\de]=\{\m=\o\}$,   
$J_{\r}[>\!\de]=\emptyset$. Finally, we note that the isomorphisms 
of a particular submodule $\m$  are those induced by the action 
of the group of units $\r^*$.

\smallskip

We can give $\bar{J}_{\r}[d]$ the structure of a projective 
subvariety in $J_{\r}=\bar{J}_{\r}[0]$ following \cite{GP}; 
actually they are projective schemes over $\Z$. 
Generally, the $\o$\~degree of $z^a\m$
for an $\r$\~module $\m\subset \o$ and $a\in \Z_+$ is 
deg$_{\o}\m+a$.
Indeed, deg$_{\o}(z^a\m)=|\Z_+\setminus \nu(z^a\m)|=a+
|\Z_+\setminus \nu(\m)|$.

Given $0\le d\le \de$,
let $\m^\circ$ be a standard submodule of $\o$\~degree 
$\de-d^\circ$ 
for  $d^\circ\ge d$ and $a=d^\circ-d$. 
Then the submodule $\m=z^a\m^\circ$ is of $\o$\~degree 
$\de\!-\!d^\circ\!+\!d^\circ\!-\!d=\de\!-\!d$
and one has:
$$
z^{2\de}\o\subset \m\subset \o \hbox{\,\, and\,\, } 
\text{dim}_{\C}(\o/\m)=\de-d.
$$ 
Vice versa, $\m^\circ=z^{-a}\m$ are standard for any  
such $\r$\~submodules $\m$, where $a$ is 
the smallest evaluation $\nu(z)$ among all $z\in \m$. 
Such $\m$ for $d=0$ are called {\em $\de$\~normalized
modules\,} in \cite{Pi}; our standard submodules are called
there {\em $0$\~normalized}.
\smallskip

Now consider the Grassmanian $Gr(\o/z^{2\de}\o,\de+d)$. 
A subspace of dimension $\de+d$ is an $\r$-module
if and only if it is invariant under the action of 
$\r/z^{2\de}\o$ by multiplication. 
Thus $\bar{J}_{\r}[d]$ is the set of fixed points in 
$Gr(\o/z^{2\de}\o,\de+d)$
under the group action of $(\r/z^{2\de}\o)^{\star}$. 
To obtain the structure we desire, extend the action of 
$(\r/z^{2\de}\o)^{\star}$ 
to $\bigwedge^{\de+d} \o/z^{2\de}\o$ and consider the image under 
the Pl\"ucker embedding. The condition of being a fixed point 
under the action $(\r/z^{2\de}\o)^{\star}$ defines a linear 
subspace of the projective space 
$\P(\bigwedge^{\de+d} \o/z^{2\de}\o)$, 
so $J_{\r}$ is the intersection of the image of the Pl\"ucker 
embedding of $Gr(\o/z^{2\de}\o,\de+d)$ and this linear subspace.
Finally, we go back to our standard modules using the 
identification maps above.
\smallskip

There is an alternative {\em intrinsic\,} 
interpretation of the projective structure of $J_{\r}$.
Let us consider 
$\r$\~submodules $\tilde{\m}\subset \k\otimes_{\C} \C[[u]]$,
finitely generated over $\C[[u]]$, i.e. a families of
submodules $\m$ (formally) analytically depending on a 
parameter $u$. Let $\tilde{\m}$
contain a principle ideal $\C[[u,z]]z^p$ for 
some $p$ (arbitrarily large).

We define the {\em enhanced (flat) limit\,} $\tilde{\m}_0$ of
$\tilde{\m}$ as the linear subspace of $\k$ of {\em all\,}
linear combinations of the vectors in $\tilde{\m}$
divided by the smallest possible power of $u$ and
then evaluated at $u=0$. This is one of the key 
definitions in the theory of sheaves and bundles over curves.

Such a limit obviously contains the straightforward
specialization $\m_0\equal\tilde{\m}(u=0)$. The
space $\tilde{\m}_0\subset \o$ has a
natural structure of a $\r$\~submodule by construction. 
Therefore it becomes a standard module (an element of $J_{\r}$) 
upon the division by a proper power $z^n (n\ge 0)$. 

As in the definition of the standard normalization, we assume 
here that $\tilde{\m}$ contains an element with zero 
$z$\~valuation; its $z$\~constant term can be any
nonzero element (possibly non-invertible) from $\C[[u]]$. 
By the way, one can check that the enhanced limit will remain 
the same if the principle ideal $\C[[u,z]]z^{\cc}$ is
added to $\tilde{\m}$, i.e. $p$ above can be assumed to be
no greater than the conductor $\cc$. This results from a 
standard theory of flat limits. Without general theory here, 
this fact can be justified by enlarging $\tilde{\m}$
with all linear combinations of the vectors in $\tilde{\m}$
divided by the corresponding minimal power of $u$  and
verifying that this procedure will eventually add
$\C[[u,z]]z^{\cc}$ to $\m$.  The limit $\tilde{\m}_0$
will remain unchanged under such an ``enhancement" of
$\tilde{\m}$.

Finally, the {\em boundary\,} of any family of modules 
considered in $J_{\r}$
is the collection of enhanced limits for
all one-parametric sub-families $\tilde{\m}$.
Using one $u$ is sufficient in this approach (another
general fact). 



\begin{proposition} \label{DEGJAC}
Let $\m'$ be the enhanced limit $\tilde{\m}_0$ of a 
$u$\~family $\tilde{\m}$ 
of modules invertible over $\C((u))$, where we assume that 
$0\in \De_{\tilde{\m}}$ as above. Setting $\tilde{\m}_0\subset 
(z^{d'})=\o z^{d'}$ for the maximal such $d'$,
one has that $d'\ge d$ for the number of added
gaps for the standard module $\m$ corresponding to $\m'$ and
$d'=d$ for sufficiently general such $\tilde{\m}$.
\end{proposition}  

{\it Proof.} For generic $\tilde{\m}$,
the limit $\tilde{\m}_0$ is of degree $0$ relative to 
$\r\in \o$; we use that $J_{\r}$ is a closed
subvariety in $Gr(\o/z^{2\de}\o,\de)$. Therefore
$z^{-d}\m'$ is standard for $d=|D_{\m}|$
and $d'\ge d$ for any (non-generic) $\tilde{\m}$. \sq


{\sf Two examples of enhanced limits.}
To illustrate Proposition \ref{DEGJAC},
let us consider $\r\!=\!\lan z^4,z^6+z^9,z^7\ran$ 
with $\Ga\!=\!\{0,4,6,7,8,10,11,\cdots\}$, $G\!=\![1,2,3,5,9]$,
$\de\!=\!|G|\!=\!5$, and the conductor $\cc=10$. 
We take the following family of modules invertible
over $\C((u))$:
$$
\tilde{\m}=\r[[u,z]](u^4+u (z+z^3)+ u^3 z^2) +\C[[u,z]]z^{10},
$$
which is an invertible $\r$\~module over 
$\C((u))$ (but not over $\C[[u]]$
due to  $u^4$).  
Then its {\em enhanced\,} limit $\tilde{\m}_0$
is a module with the corresponding set of added gaps 
$D=[2,3,5,9]$,
so $d=|D|=4$. However it is contained in $(z^5)$, which is smaller
than  $(z^4)$ guaranteed by the proposition.

Taking here $\tilde{\m}=\r[[u,z]](u(1+z^5)+z^9) +
\C[[u,z]]z^{10}$, the limit  $\tilde{\m}_0$ has the same
$D=[2,3,5,9]$. Now it strictly belongs to $(z^4)$, 
so this is the
case of a general position for such $D$. 
\medskip

The cardinalities of the $D$\~sets will be  
associated below with the parameter $q$. If 
$J_{\r}[D]$ are all affine spaces, their dimensions
give the powers of $t$. The third parameter $a$
of our construction will be due to the flag-length
of the {\em flagged Jacobian factors\,} defined as follows.

\subsection{\bf Flagged Jacobian factors}
\begin{definition}\label{JACDEF}
For a $D$-set $D$ and the sequence  
$\vec{g}=\{g_i,\, 1\leq i \leq m\}$ in $G\setminus D$ 
such that $g_i<g_{i+1}$, the
{\sf $D$\~flag} is
$$ 
\d\!=\!\{D_0\!=\!D,D_1\!=\!D\!\cup g_1,
D_2\!=\!D\!\cup\{g_1,g_2\},\ldots,
D_m\!=\!D\!\cup\{g_1,\ldots,g_m\}\},
$$
provided that all $\De_i=\Ga\cup D_i\,$ are 
standard $\,\Ga$\~modules. Then the
corresponding  
{\sf \,flagged Jacobian factor\,} $J_{\r}^{m}$
is defined as the union of the following varieties of flags of 
standard submodules $\m\subset\o$:
\begin{align}\label{flagjac}
&J_{\r}^m[\d]\equal \bigl\{\mathscr{M}=
\bigl\{\!\m_0\!\subset\! \m_1\!\subset\!\cdots\!\subset\!
\m_m\!\bigr\} \hbox{\, where\, } D_i=D_{\m_i}\bigr\}
\end{align}
for all $0\le i\le m$;
we will sometimes omit $m$ here. If at least one
such flag $\mathscr{M}$ exists  we call the 
corresponding $\,\d$ and $\,\{\De_i\}\,$
{\sf admissible}.

Considering $D$\~flags $\d$ of
length $m$, we set \, $J_{\r}^m[d>\de]=\emptyset$,
\begin{align}\label{flagjacd}
J_{\r}^m[d]\,\equal\,\bigcup_{\d, |D_0|=d}\!\!J_{\r}^m[\d]
\subset J_{\r}^{m},\ \  \overline{J}_{\r}^m[d]=
\bigcup_{d'\ge d}\!J_{\r}^m[d'],
\end{align}
where $d=0,1,\ldots,\de.$
For $m=1$, the $D$\~flag $\d$ is obtained from $D$ by adding
one gap $g$ and we use the notation $J_{\r}^1[D,g]$. 
When $m=0$, we put $J_{\r}[D]$, which are from \cite{Pi}.\sq
\end{definition}

\comment{
{\sf Hilbert schemes.}
Let us briefly discuss the switch to Hilbert schemes in this
and similar definitions. Generally, 
Hilbert schemes can be defined for flags of 
ideals $\i\subset \r$ instead of  modules $\m$.
Such ideals are special torsion free rank-one modules and
any standard $\m$ can be made inside $\r$ (the multiplication
by $z^{\cc}$ will do this). Obviously we cannot assume any longer
the standard normalization $0\in \De$ and 
use the $D$\~flags; so we  
need to go back to $\{\De_i\}$ (omitting the standard 
normalization). For any flag $\{\De_0\subset\De_1\subset
\ldots \subset \De_n\}$, let
\begin{align}\label{clusthilb}
&\hbox{Hilb}_{\r}^n(\{\De_i\})\! =\!
\bigl\{\i_0\subset\i_1\!\subset\!\cdots\!\subset\!\i_n\bigr\},
\ 0\!\le\! i\!\le\! n,\ \De_{\i_i}=\De_i.
\end{align}
Here $\De_i\subset \De_{i+1}$, but we do not generally assume 
that $|\De_{i+1}\setminus \De_{i}|=1$. However, the corresponding
flags $\{\De_i^\circ\}$ for the standard normalizations
$\De_i^\circ$ of $\De_i$ must have extensions to $D$\~flags,
as in the definition of flagged Jacobian factors. For instance,
not all pairs $\mathfrak{m}\i_n\subset \i_0\subset\i_n$ 
can be embedded into such a flag, but if $\De_0^\circ$ and
$\De_n^\circ$ come from a $D$\~flag, then $\i_0,\i_n$
have such an extension. 
See Proposition \ref{NESTED} $(iii)$; this case is
of special interest due 
to \cite{ObS,ORS,GORS} and a (possible) passage from 
algebraic knots to links
to be addressed elsewhere.



The maps 
from $\hbox{Hilb}_{\r}^n(\{\De_i\})$ to the corresponding
$J_{\r}^n(\{\De^\circ_i\})$
for small $n$  are not surjective (as well as for 
classical Jacobians), where $\De^\circ_i$ is
our standard normalization of $\De_i$. For sufficiently large
degrees, the fibers of the corresponding maps are not too
difficult to describe following Proposition \ref{NESTED}. 
For instance, the
series $\sum_{n=0}^\infty q^{n}e_{top}(\hbox{Hilb}_{\r}^n)$ for
the topological Euler characteristic $e_{top}$ can be transformed
to a finite sum, similar to the Gopakumar-Vafa formula
(see \cite{PaT}). 

The ORS conjecture states that the geometry of 
{\em nested\,} Hilbert schemes  
is related to the unreduced stable Khovanov-Rozansky 
polynomials. See also Theorem 9.5 from \cite{GORS},
a reformulation of the {\em ORS conjecture\,} in 
the reduced setting (the unibranch case). They define 
certain {\em isotypic components} 
in cohomology of special partial flag varieties with
respect to the Springer action of the corresponding symmetric
group $\S_N$, 
which eventually provide the coefficients of
$a$ in their superpolynomial. We use the {\em flagged}
Jacobian factors defined directly via $\r$ and do not need $N$, 
the rank of the corresponding matrices.

One can expect that the homology of our $J^m_{\r}[\d]$ 
upon some transformations are related to the isotypic 
components from \cite{GORS}.  However the connection of 
Theorem 9.5 (a reduced 
variant of the ORS conjecture) with our Main Conjecture is 
unclear at the moment.
}


Let us begin with some general properties of
flags $\mathscr{M}$.  We note that our usage of Nakayama's Lemma 
in Proposition \ref{NESTED} below
is actually similar to that in \cite{ORS}, Section 2.1.

 For an arbitrary module $\m$, we 
set $\m^{(i)}=\m\cap
z^i\o$, which is obviously an $\r$\~module, and
$\bar{\m}=\m/\mathfrak{m}\m$ for the maximal ideal
$\mathfrak{m}\subset \r$. Accordingly,  $\bar{\m}^{(i)}$
is the image of $\m^{(i)}$ in $\bar{\m}$. Obviously,
\{dim\,$\m^{(g)}/\m^{(g\!+\!1)}=1\Leftrightarrow
g\in \De_\m$\}\,, and\, 
dim\,$\bar{\m}^{(g)}/\bar{\m}^{(g\!+\!1)}\le 1$.

\begin{proposition}\label{NESTED}
(i) For a $D$\~flag $\d=
\bigl\{D_0\!\subset D_1\!=\!D_0\cup\{g_1\}\subset\!
\cdots\!\subset D_m\!=\!D_0\cup\{g_1,\ldots,g_m\}\bigr\}$
from Definition \ref{JACDEF},
all $D_0\cup \{g_i\}$ for $1\!\le\! i\!\le\! m$ 
are $D$\~sets. Let $\mathscr{M}=\{\m_i\}$ be a flag 
of (standard) modules corresponding 
to a $D$\~flag (admissible).
Picking an arbitrary $m_i\in \m_i$ such that
$\upsilon(m_i)=g_i$, the space $\m_0\oplus\C m_i$
is an $\r$\~module. 

Then for $1\le i\le m$,
$
\m_i=\m_0\oplus \C m_1\oplus\ldots\oplus\C m_i, 
$
$\mathfrak{m}\m_m\subset \m_0$ and 
dim\,$\bar{\m}^{(g_i)}/\bar{\m}^{(g_i\!+\!1)}=1$.
The elements $m_i$ modulo $\m_0$ are
uniquely determined up to proportionality, i.e. depend
only on the flag $\mathscr{M}$. 
\smallskip

(ii) Vice versa, for the $D$\~flag $\d$ as above, let us assume 
the existence of a module $\m_{top}$ such that $D[\m_{top}]=D_m$ 
and 
dim\,$\bar{\m}_{top}^{(g_i)}/\bar{\m}_{top}^{(g_i\!+\!1)}\!=\!1$ 
for all $i$. We do 
not impose the admissibility of the whole $\d$. 
Then $\m_m=\m_{top}$ can be extended to a flag 
$\mathscr{M}$ corresponding to $\d$ (so it is admissible)
and all such flags can be described as follows. 

(a) For any subspace $\bar{\m}_0\subset \bar{\m}_{top}$  
such that\, dim\,$\bar{\m}_{top}/\bar{\m}_0=m$,
dim\,$(\bar{\m}_0\!+\!\bar{\m}_{top}^{(g_i)})/
(\bar{\m}_0\!+\!\bar{\m}_{top}^{(g_i\!+\!1)})=1$ for
$1\le i\le m$, $\bar{\m_0}\not\subset
\bar{\m}_{top}^{(1)}$, let $\mathfrak{m}\m_{top}\subset 
\m_0\subset \m_{top}$
be a unique lift of $\bar{\m}_0$ 
(the Nakayama Lemma).


(b) Then 
 for any pullbacks $\bar{m}_i$ of the generators of the
latter quotients to 
$(\bar{\m}_0+\bar{\m}_{top}^{(g_i)})/\bar{\m}_0$, 
we take $\mathscr{M}=\{\m_i\}$, where
$\m_i$ is a similar lift of the space 
$\bar{\m}_i=\bar{\m}_0\oplus \C \bar{m}_1\oplus\ldots\oplus
\C \bar{m}_i$ to an $\r$\~submodule of $\m_{top}$
containing $\mathfrak{m}\m_{top}$. These imply that \,
$D[\m_0]=D_0,\,\m_{i=m}=\m_{top}$.
\smallskip

(iii) Continuing (ii),
all flags $\mathscr{M}$ in (i) are uniquely described by the 
following data: ($\al$) the space $\bar{\m}_0\subset 
\bar{\m}_{top}$ as above,\ ($\be$) the set of 
elements 
$\bar{m}_i\in (\bar{\m}_0+\bar{\m}_{top}^{(g_i)})/\bar{\m}_0$  
considered up to proportionality. I.e. for a fixed
$\m_{top}$ and  $\bar{\m_0}$ as in (ii), such flags of modules
are naturally parametrized by unipotent complex matrices 
of size $m\times m$, which is due to the action of the
Borel subgroup preserving the full flag  
$\{(\bar{\m}_0+\bar{\m}_{top}^{(g_i)})/\bar{\m}_0\}$ in
 $\bar{\m}_{top}/\bar{\m}_0$. 
\end{proposition}

{\it Proof.} Part $(iii)$ formally follows from $(i,ii)$,
including the equivalence of the admissibility of a $D$\~flag
$\d$ and the existence of the pair $\{\m_0,\m_{m}\}$ for
$D_0,D_m$ such that $\mathfrak{m}\m_{m}\subset 
\m_0\subset \m_{m}$. Let us justify $(i,ii)$. Below $\{g_i\}$
means a single element $g_i$ considered as a set.

First of all, $D_0\cup \{g_i\}$ correspond to 
certain $\Gamma$\~modules for any $1\le i\le m$.
Indeed, adding $g_i$ to $D_{i-1}$ cannot add anything
new to $D_i$ (since the corresponding
$\Delta_i$ is assumed a $\Ga$\~module).
The same holds for $D_0\cup \{g_i\}$, since extra
elements in the minimal $\Gamma$\~module containing 
the latter can be only greater than $g_i$, i.e. can be
only in $D_i\setminus D_{i-1}$

This reasoning is equally applicable to the 
flags $\mathscr{M}$ for admissible $\d$. 
We claim that $D_0\cup \{g_i\}$ is admissible for each $g_i$. 
Indeed, adding  
$\r m_i\subset \m_i$ to  
$\m_0$ for an element $m_i$ with the valuation
$\upsilon(m_i)=g_i$ cannot create any 
new valuations vs. $\De_0$ but $g_i$, since this does not create 
them  when going from $\m_{i-1}$ to  $\m_i$. Thus 
$\m_0+\r m_i=\m_0+\C m_i$ is an $\r$\~module corresponding 
to $D_0\cup \{g_i\}$.

Next, $\psi m_i$ must belong to $\m_0$
for any $\psi\in \mathfrak{m}$. We use that
for any $N>0$ there exists an element $m'\in \m_0$ such
that  $\upsilon(\psi m_i-m')>N$; thus this difference can 
be assumed in $\m_0$.  

Similarly, $\m_i$ is the linear span of $\m_{i-1}$ and $m_i$
for the elements $m_i$ introduced above,
since for every $m\in \m_i$, there exists $m'\in
\m_{i-1}+\r m_i$ such that $\upsilon(m-m')$ is greater
than any given number. Thus $m-m'$ can be assumed
in $\m_{i-1}$.  We obtain that adding the elements 
$\{m_i,1\le i\le j\}$ to $\m_0$ generates $\m_j$.
Combining this claim for $j=m$ with 
 $\mathfrak{m}\, m_i\in\m_0$ checked above, we conclude that
$\mathfrak{m}\m_m\subset \m_0$.

Part $(ii)$ uses the same arguments and 
the Nakayama Lemma. We lift $\bar{m}_i$ to arbitrary
elements in $\m_{top}^{(g_i)}$. The corresponding
$\r$\~modules does not depend on such choices.
Note that 
the condition $\bar{\m}_m=\bar{\m}_{top}$ provides that
all $\m_i$ are standard ($\De_i\ni 0$). 
We need
only to check that these modules
have the required $D$\~sets,
which results from the following lemma by induction
with respect to the length $m$ of $\mathscr{M}$.

\begin{lemma} For a pair of modules $\m'\subset\m$, 
where $\m'$ is not necessarily under the standard
normalization, let $g\in D[\m]$, $\mathfrak{m}\m\subset \m'$, 
dim\,$(\bar{\m}'\!+\!\bar{\m}^{(g)})/
(\bar{\m}'\!+\!\bar{\m}^{(g\!+\!1)})=1$ and
$\bar{\m}=\bar{\m}'\oplus \C\bar{m}$ for a
pullback $\bar{m}$ of the generator of the latter
quotient. Then we claim that $\m'$ is a standard
module and  $D[\m']=D[\m]\setminus\{g\}$.
\end{lemma}
{\it Proof.} The module
$\m$ is linearly generated by $\m'$ and
any lift of  $\bar{m}$ to $\m$ (the Nakayama Lemma).
Since
dim\,$(\bar{\m}'\!+\!\bar{\m}^{(g)})/
(\bar{\m}'\!+\!\bar{\m}^{(g\!+\!1)})$ is $1$,
the gap $g$ cannot be among the valuations of $\m'$
and this is the only missing gap there vs. $\m$,
which gives the required.  \sq 
\smallskip

Part $(iii)$ follows from $(ii)$. 
The inequalities dim\,$\bar{\m}_{top}^{(g_i)}/
\bar{\m}_{top}^{(g_i\!+\!1)}>0$ are obviously
necessary and sufficient to ensure the existence
of the required $\bar{\m}_0$. The rest is a combination
of $(i)$ and $(ii)$.\sq
\smallskip

Thus given a $D$\~flag $\d$ as in $(i)$, and
$\m_{top}$ such that $D[\m_{top}]\!=\!D_m\!=\!
D_0\cup\{g_1,\ldots,g_m\}$, the inequalities 
dim\,$\bar{\m}_{top}^{(g_i)}/
\bar{\m}_{top}^{(g_i\!+\!1)}=1$ are   
equivalent to the admissibility of $\d$. 
We note that in all examples we calculated, the admissibility of
a $D$\~flag is equivalent to that of
every $D_i$. We cannot justify this in general. 

This proposition is of independent interest and can 
be used to count the dimensions of $J_{\r}[\d]$ 
and to check (in some cases) that the latter are 
biregular to affine spaces.

\comment{
The most general approach to
establishing that $J_{\r}[\d]$ are affine is as follows.
Let us consider our modules $\m$ up to 
isomorphisms (we will use $\cong$ for this).
In the standard realization as submodules of $\o$, one has:
$
\m\cong\m' \Longleftrightarrow \m'=\phi\m \for
\phi\in \o^*.
$
Equivalently, $\m'\l=\m$ for an invertible
$\r$\~module $\l$ (submodule of $\o$ of rank one under
the standard normalization), where
the multiplication in $\o$ is used.
}

\smallskip
We note that there is a natural (biregular)
action of the algebraic group $J_{\r}[\varnothing]=\{\l\}$ in 
$J_{\r}$. Without going into details, Theorem 1 from \cite{Pi}, 
can be extended to flagged Jacobian factors as follows.
If the GIT quotient $J^m_{\r}[\d]\/\!/J_{\r}[\varnothing]$ is 
biregular to an affine space and with the stabilizers of points
in $J^m_{\r}[\d]$ of the same dimension, then the latter
space is biregular to some $\mathbb{A}^{\!N}$. This is not always 
the case (see online version of this paper),
but it is not impossible that $J^m_{\r}[d]$ can be always paved 
by affine spaces.

\comment{ 
\begin{proposition}\label{Affinespace}
The $J_{\r}[\varnothing]$\~orbit of
any flag $\mathscr{M}$ is biregular 
to an affine space $\mathbb{A}^{\!N}$. 
Let $\d$ be a $D$\~flag. If the GIT quotient
$J^m_{\r}[\d]\/\!/J_{\r}[\varnothing]$ is biregular 
to an affine space and with the stabilizers of points
in $J^m_{\r}[\d]$ of the same dimension, then the latter
space is biregular
to some $\mathbb{A}^{\!N}$. 
\end{proposition}
{\it Proof.} The first claim is essentially from \cite{Pi} (the 
Chevalley -– Rosenlicht theorem is used). The second one
results from 
the first and the Quillen-Suslin theorem
(Serre's problem). \sq
\smallskip
When calculating $J^m_{\r}[\d]/\!/J_{\r}[\varnothing]$
we can assume that $1\in \m_0$, which condition  
``eliminates" the action of $J_{\r}[\varnothing]$ up to
the units of the ring of automorphisms of a flag $\m$.  
This is not always the case (see Appendix \ref{app:noaff}),
but we think that $J^m_{\r}[d]$ can be always paved by
affine spaces.
}

\subsection{\bf The Main Conjecture}
The flagged Jacobian factor $J_{\r}^{m}$
has a natural structure of a projective variety. 
Accordingly,
$J^m_{\r}[\d]$ and $\overline{J}^m_{\r}[d]$ are its subvarieties;
the latter is closed. Proposition \ref{NESTED}
gives that $J^m_{\r}$ is a certain {\em subvariety\,} in a proper 
{\em parahoric Springer fibers\,} 
in the terminology from Section 2.2.9 from
\cite{Yu}. They are formed by {\em partial\,} periodic flags of 
$\r$\~invariant lattices. We do not use them in the 
present work.  By $H_i$, we will mean singular (relative)
homology with the
$\C$\~coefficients.

\begin{conjecture} \label{CONJSUP}
Let $\r$ be the ring of a unibranch plane curve
singularity $\c_{\vec\rr,\vec\ss}\,$ 
from Section \ref{sec:ALG-KNOTS},
$\h_{\,\vec\rr,\,\vec\ss}\,(\square;q,t,a)$
be the DAHA uncolored superpolynomial
from (\ref{jones-sup}). Recall that
$\de=|G_{\r}|$ is the arithmetic genus of $\r$;
we assume that $\rr_1>\ss_1$.
\smallskip

(i) We conjecture that the relative homology 
$H_{2i+1}(\overline{J}^m_{\r}[d],\overline{J}^m_{\r}[d\!+\!1])$
vanishes for all $i,d\ge 0$ and   
\begin{align}\label{conjcoh}
\h_{\,\vec\rr,\,\vec\ss}\,(\square;q,t,a)=
\sum_{d,i,m}\hbox{rk}\bigl(H_{2i}
(\overline{J}^m_{\r}[d],\overline{J}^m_{\r}[d\!+\!1])\bigr)
q^{d+m} t^{\delta-i} a^m,
\end{align}
where $0\le d,i \le \de$ and the range of $\,m\,$
is from $0$ to
$\,\ss_1\rr_2\cdots \rr_{\ell}-1$,
which is the number of (admissible) $D$ such that 
dim$\left(J_{\r}[D]\right)=\de\!-\!1.$ The right-hand 
side of (\ref{conjcoh}) will be called $\h^{h\!om}(q,t,a)$. 
Also, we conjecture that $J^m_{\r}[d]$ are paved 
by affine spaces.  
\smallskip

(ii) If all varieties $J_{\r}^m[\d]$ are 
affine
spaces $\mathbb{A}^{\!i}$ then their total number 
for any fixed $i,d$ gives the corresponding rank 
rk$(H_{2i})$ from (\ref{conjcoh}) and this relation reads:
\begin{align}\label{conjaff}
\h_{\,\vec\rr,\,\vec\ss}\,(\square;q,t,a)=
\sum_{\d} q^{d+m\,} t^{\de-\hbox{\tiny dim}
(J^m_{\r}[\d])\,} a^m,
\end{align}
where the summation is over all admissible $D$\~flags
$\d, d\!=\!|D_0|$. Here the cells $J_{\r}^m[\d]$ contribute to 
$H_{2i}(\overline{J}_{\r}^m[d])\!\subset\! H_{2i}(J_{\r}^m)$ 
for $i\!=\!\hbox{\small dim}(J^m_{\r}[\d])$,  and
then to $H_{2i}(\overline{J}_{\r}^m[d],
\overline{J}_{\r}^m[d\!+\!1])$. Formula (\ref{conjaff}) can
be readily extended to any affine cell decompositions of
$J^m_{\r}[d]$ (not only those via $J^m_{\r}[\d]$).
\smallskip

(iii) 
Let $1/t=p^{\ell}$ for prime $p$ and $\ell\in \N$, 
$\mathbb{F}=\mathbb{F}_{1/t}$ the field with $p^{\ell}$ elements, 
$|X(\mathbb{F})|$ the number of $\mathbb{F}$\~points
of a scheme $X$ defined over $\mathbb{F}$.
One can assume that $\r$ is defined
over $\Z$ and consider $J^m_{\r}[d]$ as schemes
over $\mathbb{F}$. We conjecture that apart from finitely 
many $p$,
\begin{align}\label{conjmot}
\h_{\,\vec\rr,\,\vec\ss}\,(\square;q,t,a)=
t^{\delta}\sum_{d,m}|J^m_{\r}[d](\mathbb{F})|
q^{d+m} a^m \equal \h^{mod}(q,t,a).
\end{align}
If  $J^m_{\r}[d]$ are (a)\, paved by 
affine spaces over  $\mathbb{F}$, and (b)\, (non)admissible $\d$ 
remain such over this field, then (\ref{conjmot})
is equivalent to (\ref{conjcoh}).
\sq 
\end{conjecture}

{\sf Relative homology.}
Vanishing $H_{2i+1}(\overline{J}^m_{\r}[d],
\overline{J}^m_{\r}[d\!+\!1])$ generalizes 
the {\em van Straten-Warmt conjecture\,}, 
which claims that odd Betti numbers 
of $J_{\r}$ vanish. This assumption and the exact
sequences for relative (singular) homology imply that
the natural maps $H_{2i+1}(\overline{J}^m_{\r}[d\!+\!1])
\to H_{2i+1}(\overline{J}^m_{\r}[d])$ are surjective, which
readily gives that $H_{2i+1}(\overline{J}^m_{\r}[d])=\{0\}$
for any $i,d$. Using this, the following
sequence is exact: 
 \begin{align}\label{relcoho}
0\!\to\! H_{2i}(\overline{J}^m_{\r}[d\!+\!1])
\!\to\! H_{2i}(\overline{J}^m_{\r}[d])\!\to\!  
H_{2i}(\overline{J}^m_{\r}[d],\overline{J}^m_{\r}[d\!+\!1])
\!\to\! 0.
\end{align}


As an application, (\ref{conjcoh}) upon the 
substitution $q=1$ becomes: 
$$
\h_{\,\vec\rr,\,\vec\ss}\,(\square;q=1,t,a)=
\sum_{i,m}\hbox{\em rk\,}H_{2i}(J^m_{\r}) t^{\de-i} a^m,
$$
where the special case $a=0$ is Conjecture 2.4, $(iii)$
from \cite{ChD1}.  We note that 
$H_{2i}(\overline{J}^m_{\r}[d],
\overline{J}^m_{\r}[d\!+\!1])$ are potentially connected with
compactly supported cohomology 
$H^{2k}_c(J^m_{\r}[d])$
for $k\!=\!j$ or $k\!+\!j\!=\!\hbox{dim}\,J^m_{\r}[d]$, though 
Poincar\`e duality fails for
singular varieties. Compare with (\ref{conjmotxx}).

\smallskip

{\sf Motivic approach.\,}
Relation (\ref{conjmot}) readily follows from Part $(ii)$. 
We think that it can generally hold. In examples, it suffices to 
take any $p^\ell=1/t$ such that $\Gamma$ remains the
same over $\mathbb{F}_{1/t}$, but we do not know how far
this goes. Say, (\ref{conjmot}) holds for 
$\r=\mathbb{F}_2[[z^4+z^5,z^6]]$ under $a=0$, but not for 
$\mathbb{F}_3$, which changes $\Gamma$. 
Though $\r=\mathbb{F}_3[[z^4,z^6+z^7]]$ is fine and we note that 
$\nu(\mathbb{F}_2[[z^4+z^5,z^6]])=
\nu(\mathbb{F}_3[[z^4,z^6+z^7]])$. 
These rings over $\C$ correspond to coinciding 
DAHA-superpolynomials. See (\ref{conjmotxx}) for
a reformulation of (\ref{conjmot}) in terms of the weight
filtration.

When $a=0,q=1$, Part $(iii)$
is closely related to Theorems 0.1,0.2 from \cite{GKM} 
in the case of $A_n$ for affine Springer fibers. 
The latter are not always paved by affine
spaces for other types \cite{KL}. The $A_n$\~case is 
exceptional; the positivity of the orbital 
integrals and vanishing odd rational homology are widely 
expected to hold for anisotropic centralizers (i.e.
in the nil-elliptic case). This matches well the conjectured 
positivity of uncolored DAHA superpolynomials.

Knowing the spectral curve is sufficient
here. However, it is far from obvious beyond the torus case
(quasi-homogeneous plane curve singularities) that only
the  {\em topological\,} type of the singularity matters here.
This follows from of our conjecture. The analytic equivalence of
plane curve singularities is generally very different from the 
topological perspective. See
\cite{La,ChU} and Section \ref{sec:ORBITAL}
for the identification of affine Springer fibers 
of anisotropic type with the compactified Jacobian
of rational curves such that $\r$ is the local ring 
at its {\em unique\,} singular point. 

In full generality, with the nonzero
parameters $q$ and $a$  in (\ref{conjmot}),  
we count points in $J^m_{\r}[d](\mathbb{F})$
with $q$\~weights $q^{|D_m|}$, where 
$|D|\!=\!\de\!-\!\hbox{deg}_\o\m$ 
for the corresponding {\em standard\,} $\m$. This seems a new turn
in the theory of affine Springer fibers and related 
$p$\~adic {\em orbital
integrals}. Possible (expected) adding colors and the
multi-branch generalization of our construction
 make this even more interesting. Generally, we think that
using the powerful modern theory of {\em topological\,}
invariants of plane curve singularities can be expected to impact
the theory of orbital integrals (at least in type $A$) and related 
part of the Fundamental Lemma;  
see Sections \ref{sec:ORBITAL},\ref{sec:MOTIV} 
for further discussion. 
\medskip

{\sf The case $\r=\lan z^4,z^6+z^7\ran$.}
We will stick to {\em admissible
pairs\,} $\{\m=\m_0\subset\m_1\}$, where 
$D_1=D_0\cup \{g_1\}$ (corresponding to the coefficient
of $a^1$ in the Main Conjecture). Recall that $D_i=
\De_i\cap G$, $G=\Z_+\setminus \Ga$.
Using directly the definition of
$J_{\r}[D_0\!\subset\!D_1]=J_{\r}^1[D_0,g_1]$,
we obtain the following lemma.

\begin{lemma}\label{ADD1GAP}
Provided the admissibility of the $\Ga$\~modules
$D_0$ and $D_1=D_0\cup \{g_1\}$ for $g_1\not\in \Ga\cup D_0$,
one has: 
\begin{align}\label{add-oneg}
&\hbox{dim} J_{\r}[D_0\!\subset\! D_1]\le 
\hbox{dim} J_{\r}[D_1]\,+\,|\,\{g\in \De_0\mid g<g_1\}\,|\,,\\
&\hbox{dim} J_{\r}[D_0\!\subset\!D_1]=
\hbox{dim} J_{\r}[D_1]+1 \hbox{\,\, if\,\, }
g_1 < g, \,\forall g\!\in\!D_0,\notag\\ 
&\hbox{dim} J_{\r}[D_0\!\subset\!D_1]=
\hbox{dim} J_{\r}[D_0] \hbox{\,\, if\ \,}
\{g\!\in\! G\!\mid g\!>\!g_1\}\!\subset\! D_0.\label{add-onem}
\end{align}

In the case of $\r=\lan z^4,z^6+z^7\ran$, 
the pairs $\{D_0,g_1\}$ that are not included in the
latter two formulas are all governed by the first one
where ``$\le$" is replaced by ``$=$" there. Thus these 
formulas 
provide all dimensions of the varieties of admissible 
pairs $D_0\!\subset\! D_1$. Using this for such $\r$, 
one obtains the formula 
\begin{align*}
\sum_{D_0\subset\!D_1} q^{|D_0|+1\,} t^{\de-\hbox{\tiny dim}
(J^m_{\r}[D_0\subset\!D_1])}=q + q^2 \bigl(1 + t\bigr) + 
q^3 \bigl(1 + 2 t + t^2\bigr)\\ 
+ q^4 \bigl(3 t + 2 t^2 + t^3\bigr) + 
    q^5 \bigl(t + 4 t^2 + 2 t^3 + t^4\bigr) + q^6 \bigl(t^2 
+ 4 t^3 + 2 t^4 + t^5\bigr)\\ 
+ q^7 \bigl(t^3 + 3 t^4 + 2 t^5 + t^6\bigr) + q^8 \bigl(t^5 
+ t^6 + t^7\bigr),
\end{align*}
which matches the coefficient of $a^1$ of the uncolored DAHA 
superpolynomial
for $\r=\lan z^4,z^6+z^7\ran$ from \cite{ChD1}. 
\end{lemma} 

{\it Proof}. To justify the first inequality, we begin
with any parametric family
of modules $\m_0$ corresponding to
$D_0$, assuming that they can be extended
to $\m_{g_1}$ by adding certain $m_{g_1}$ with
$\upsilon(m_{g_1})=g_1$ and that they are different 
(as submodules of $\o$) for different values of
parameters.  The dimension of the resulting
family of modules $\m_1$ is no greater than
$\hbox{dim} J_{\r}[D_1]$.
Given  $\m_1\ni m_{g_1}$, its different submodules $\m_0$
can be only obtained by adding the terms from $\C m_{g_1}$ to
the generators $m_{g}\in \m_0$ with $\upsilon(m_g)=g$,
where $g\in D_0$ such that $g<g_1$. This gives the required
inequality.

If $g_1<D_0$, then only $m_0=1+z(\cdot)$ can be altered 
by $\C m_{g_1}$ and there is a one-dimensional family
of (pairwise distinct) such submodules $\m_0$ inside fixed
$\m_1$, which gives the second formula
in (\ref{add-oneg}). Here we use that the relations for
the coefficients of the generators $m\in \m$ necessary for
the equality $D_{\m}=D$ allow such a deformation. 

Generally, if there several such $m_g$ (not
just $g=0$) with  $g$ before $g_1$, then 
it is not true that the increase of dimension 
for $D[\m_0]$ will coincide with the number of such $g$.
It really occurs with  $5$ such 
admissible pairs $\{D_0,D_1\}$ for 
$\r=\lan z^4,z^6+z^{13}\ran$:

\begin{align}\label{4pairs}
&D_0=[2, 9, 15], \ \ \ \ \  D_1=[2, 9, 11, 15], 
&dim_0=7, dim_1=6;\\
&D_0=[2, 7, 11, 15], \ D_1=
[2, 7, 9, 11, 15], &dim_0=6, dim_1=5;\notag\\
&D_0=[2, 9, 11, 15], \ D_1=[2, 7, 9, 11, 15], 
&dim_0=6, dim_1=5;\notag\\
&[2, 3, 7, 9, 11, 15],\ \ \ \ \ [2, 3, 5, 7, 9, 11, 15],
&dim_0=6, dim_1=5;\notag\\
&[1, 2, 5, 7, 9, 11, 15], [1, 2, 3, 5, 7, 9, 11, 15],
 &dim_0=3, dim_1=0, \notag
\end{align}
where $dim_i=\hbox{dim}(J_{\r}[D_i])$. However in all these 
cases, $\hbox{dim}(J_{\r}[D_0\!\subset\! D_1])$ coincides with
$dim_0$, which is (\ref{add-oneg}) with strict equality there.  

Finally, formula (\ref{add-onem})
holds because adding {\em any} element
$m_{g_1}\in \o$ with such a valuation $g_1$ to $\m_0$ results
in the required $\m_1$ in this case. 
{\em All\,} cases from (\ref{4pairs})
are of the latter type. 
Thus any admissible pair for  $\r=\lan z^4,z^6+z^{7}\ran$
satisfies either (\ref{add-oneg}) with the equality there or 
(\ref{add-onem}). Then one can
use the formulas for dim\,$J_{\r}^{m=0}[D]$ from \cite{Pi}.

\setcounter{equation}{0}
\section{\sc The family 
\texorpdfstring{ 
\mathversion{normal}$\C[[z^4,z^{2u}+z^v]]$}
{z\textasciicircum 4 }}
\subsection{\bf Dimensions of 
\texorpdfstring{\mathversion{normal}$J_{\r}^m[\d]$}
{J[D]}}
It is shown in \cite{Pi} that
the varieties $J_{\r}[D]$ are isomorphic to affine spaces for
the family $\r=\C[[z^4,z^{2u}+z^v]]$, where $(uv,2)=1$, 
$v>2u$ and their dimensions are computed. 
We will generalize these claims to admissible $D$\~flags
and the corresponding varieties 
$J_{\r}^m[\mathcal{D}]$. 

To write 
our formula for the dimensions we will need a few definitions. 
Let 
\[
\mu_{D,g}\equal\ \text{dim}\left(J_{\r}^1[D,g]\right)-
\text{dim}\left(J_{\r}[D]\right)
\]
 be the \emph{dimension change}. Also, we write 
 \[
 \gamma_{\Delta}(\ell)\equal\ | [\ell,\infty)\setminus \Delta|
 \]
 for the \emph{gap counting function}. Recall that
$D=\De\setminus \Ga$ is called a $D$\~set corresponding to
a {\em standard\,} $\Ga$\~module $\De$ (i.e. $0\in \Delta$).

\begin{theorem}\label{dimprop}
Let $\r=\C[[z^4,z^{2u}+z^v]]$,  where $(uv,2)=1$, $v>2u$,
and $\d$ be an admissible $D$-flag (with $D_i$ corresponding
to admissible $\De_i$). Then $J_{\r}^m[\d]$ 
(with the notation as above) is isomorphic to an affine 
space $\mathbb{A}^{\!N}$ with
\[
N=\text{\emph{dim}}\left(J_{\r}[D_0]\right)
+\sum^{m-1}_{i=0}\mu_{D_i,g_{i+1}}, \hbox{ where:\,\, }
\]
\[
\mu_{D,g}=\left\{\begin{array} {lrl}
      \gamma_{\Delta\cup \{g\}}(g)-\gamma_{\Delta\cup \{g\}}
(g+4) &\hbox{for\,\,} 
      g\equiv 1\text{ or }3&\hspace{-.5cm}\mod 4, \\[.5em]
     \begin{minipage}{6cm}
     $\gamma_{\Delta\cup \{g\}}(g)-\gamma_{\Delta\cup \{g\}}(g+4)
      \\ -(\gamma_{\Delta\cup \{g\}}(g+n)-
\gamma_{\Delta\cup \{g\}}(g+4+n))$
      \end{minipage} &\hbox{for\,\,} g\equiv 2 & 
\hspace{-.5cm}\mod 4, \\
   \end{array}\right.
\]
\vskip 0.1cm
\noindent
$n=n(D)$ is the smallest odd number in 
$(\Delta \cup {v})\cap [2u,\infty)$. Here $D=D_i,g=g_{i+1}$
or $D$ is any admissible $D$\~set  
such that $D\cup\{g\}$ is admissible.
\end{theorem}

\subsection{\bf Basic Definitions} 
Before we proceed with the proof, we will briefly 
summarize and adjust the approach taken in \cite{Pi} to prove 
that the $J_{\r}[D]$ are affine; our argument 
relies heavily on his method. Let $\Delta$ be a $\Ga$\~module 
and begin by choosing $a_0,a_1,a_2$ and $a_3$ such that 
$a_i=\min\{k\in \Delta\hspace{1mm}|\hspace{1mm}k\equiv i\mod 4\}$.
Consider the following elements in $\o$: 
\begin{align}
\label{gen}
m_0=1+\sum_{k\in \N\setminus \Delta}\lambda_k^{0}z^k, \ 
m_1=z^{a_1}+\!\sum_{k\in [a_1,\infty)\setminus
\Delta}\lambda_{k-a_1}^{1}z^k, \\
m_2=z^{a_2}+\!\sum_{k\in [a_2,\infty)\setminus
\Delta}\lambda_{k-a_2}^{2}z^k, \ 
m_3=z^{a_3}+\!\sum_{k\in [a_3,\infty)\setminus
\Delta}\lambda_{k-a_3}^{3}z^k,\notag
\end{align}
where the $\lambda$\~coefficients are treated as variables. 
The valuation $\Gamma$-module of the module $\m$ 
generated by $\{m_i\}$ will then contain $\Delta$, since any 
element of $\Delta$ has the form $a_i+4n$ for 
some $n\in \N$ and because $z^4\in \r$. Thus 
$\{a_i\}$ form a \emph{basis} for $\upsilon(\m)$
in a natural sense.

An important component of
Piontkowski's method is the  
observation that $\upsilon(\m)=\Delta$ if and only if the 
relations 
among the elements $m_i$ do not produce elements with 
valuation not in $\Delta$. Thus the {\em syzygies} of the 
set $\{m_i\}$ as well as the syzygies of the set of
their leading terms are of importance. Lemma 11 of 
\cite{Pi} uses this basic idea to give an equivalent condition for 
$\upsilon(\m)=\Delta$; it will be provided. We need 
the notion of 
\emph{initial vector}, which is also from \cite{Pi}, 
to state the aforementioned lemma:

\begin{definition}\label{DEFINITX}
For $\vec{r}=(r_0,\ldots,r_3)\in \r^4$, let
 $\sigma=\min\{\upsilon(r_j)+a_j\}$. The 
{\sf initial vector} $in(\vec{r})$ is as follows:  
$in(\vec{r})=(\ze_j)$ 
with $\ze_j$ equal to the monomial of lowest degree in 
$\,r_j\,$ if $\,\upsilon(r_j)+a_j=\sigma\,$ and $\,0\,$ otherwise.
\end{definition}


\begin{lemma}\label{PIO}
Let $\m$ be an $\r$-module generated by 
$\{m_0,m_1,m_2,m_3\}$,  
$V$ an $\r$\~submodule in  $\oplus_{j=0}^3 \r$   
such that the initial vectors  
$\{\text{in}(\vec{r})\mid \vec{r}\in V\}$ of $V$ linearly generate 
the {\sf syzygies} 
of the set $(z^{a_j})$ of $\C[\Delta]$. Here 
$\C[\Delta]$ 
is the vector space generated
by the elements $\{z^k\}$ for $k\in \Delta$. 
We use  
$\sigma=\min\{\upsilon(r_jm_j)\}=\min\{\upsilon(r_j)+a_j\}$
from Definition \ref{DEFINITX}.

Then  $\upsilon(\m)=\Delta$ if and only if for each 
$\vec{r}=(r_j)\in V$ the initial terms $\,in\,$  
$\sum_{j=0}^3 r_jm_j$ cancel, i.e., 
$\upsilon(\sum_{j=0}^3 r_jm_j)>\sigma$ and for every $j$,
there exists $p_j\in \r$ such that 
$\upsilon(p_jm_j)>\sigma$ and $\sum_{j=0}^3r_jm_j=
\sum_{j=0}^3p_jm_j$. 
If such $p_j$ exist, then the element
$\sum_{j=0}^3 p_jm_j\in \m$ 
is called a {\sf higher order expression} for 
$\sum_{j=0}^3 r_jm_j\in \m$, which is generally not
uniquely determined by the latter element.\sq
\end{lemma}

We will need the following
\emph{reduction procedure} from \cite{Pi} for 
series $y\in \o$. Let $y_0=y$. Then define inductively, 
$y_{i+1}=y_i$ if $i\not\in \Delta$ and we also set 
$s_i=0$ in this case.
If $i\in \Delta$, then find the monomial, $c_iz^i$, with power 
$i$ in 
$y_i$ and the element
$s_i\in \r$ such that $s_im_{j_i}=c_iz^i+...$ for 
one of the   
generators $m_{j_i}$ (it can be non-unique). Then we 
let $y_{i+1}=y_i-s_im_{j_i}$. The sequence of
elements $y_i\in \o$ converges to an element $y_{\infty}$, 
which has the form $\sum_{k\in \N\setminus\Delta}d_k z^k$ 
for some coefficients $d_k$.
 
The key facts from \cite{Pi} about the reduction procedure are 
that (a) $y_{\infty}=0$ if and only if $y\in \m$ and (b) 
if $y_{\infty}=0$
then $\sum s_i m_{j_i}$ is a higher order expression for $y$.

\begin{definition}
For any element $y\in \o$, $y^{\dagger}$ will be the 
result of the reduction procedure
applied to $y$. One has:
$((y)^{\dagger})^{\dagger}=y^{\dagger}$.
\end{definition}

The reduction procedure depends on the choices above. 
We can standardize the procedure by always taking the elements
$s_i$ involved to be 
of the form $(z^4)^k$ for some non-negative integers $k$, which 
then makes the reduction procedure unique. To see that picking
such $s_i$ is possible, observe that when eliminating 
$c_iz^i$ for $i\in \Delta$, we can choose a unique
generator $m_j$ such that 
$a_j\equiv i\mod 4$. We will call such a procedure {\em 
standard}, and we will always assume such a standardization 
in what will follow. 
 
One has 
$x^{\dagger}=x+f$\, for 
some $f\in \m$. If for two elements $x,y\in \o$, we know that
$x^{\dagger}=x+f_1$ and $y^{\dagger}=y+f_2$, then 
the standard reduced form $(x+y)^{\dagger}$ of $x+y$\, is\,
$x+y+f_1+f_2$. Generally,
if $(x+y)^{\dagger}=x+y+f$; the uniqueness
guarantees  $f=f_1+f_2$, and we obtain
that 
$x^{\dagger}+y^{\dagger}=(x+y)^{\dagger}$ for the
standardization we will always impose.
Thus $\dagger$ is a $\C$-linear projection.

Piontkowski proves that $J_{\r}[D]$ are affine by using the 
previous lemma and the reduction procedure. In $\ref{PIO}$, it 
is sufficient to 
consider the degrees of syzygies of $(z^{a_j})$ less than 
$\max\{a_j:j=0,...,3\}-3$ 
and there are only finitely many of such syzygies. Thus it is  
sufficient to check that only finitely many linear combinations 
of the $\{m_i\}$
prescribed by the syzygies have 
higher order expressions to ensure that $\Delta=\upsilon(\m)$. 
To obtain the higher order expressions of these
elements the reduction procedure is used and the
resulting elements of $\o$ must vanish since these elements 
are in $\m$. 
The coefficients of the 
higher order expressions are polynomials in terms of  
the coefficients of $\{m_i\}$, which have the form
\[
\lambda_k^j-\lambda_k^i+\text{polynomial in }\lambda^*_{\ell}
\text{, }\ell<k.
\]
These expressions must vanish due to the properties of
the reduction procedure. Since they vanish and are
linear in the parameters 
$\lambda^j_k,\lambda^i_k$ we can express $\lambda^j_k$ in 
terms of $\lambda^*_{\ell}$ for $\ell<k$ and $\lambda^i_k$. 

Applying the same process to $\lambda_{\ell}^*$ 
such that $\ell< k$, we can eventually show that
$J_{\r}[D]$ is a graph of a 
regular function on an affine space $\mathbb{A}^n$. Since the 
graph of a regular function is always isomorphic to the domain of 
the function, we have that $J_{\r}[D]\cong \mathbb{A}^n$.
We will use a similar technique to prove 
Theorem \ref{dimprop}.

\subsection{\bf The cells
\texorpdfstring{{\mathversion{bold}$J^m_{\r}[\d]$}}{ for flags }
are affine}
Let $\d$ be an admissible $D$-flag with 
$\vec{g}=\{g_i,1\leq i\leq m\}$. 
Define  $\e$ to be the admissible $D$-flag of length $m-1$ 
such that 
$E_0=D_0$ and $E_{j}=D_0\cup \{g_i,1\leq i\leq j\}$ where 
$1\leq j\leq m-1$, i.e. $\e$ is a truncation of $\d$. 

Let $\mathscr{M}$ be an element of $J_{\r}^{m-1}[\e]$ and define  
\[
h_{g_m}=z^{g_m}+\sum_{k\in [g_m,\infty)\setminus \Delta_m}
\lambda^{h_{g_m}}_{k-g_m}z^k.
\]
We can extend $\mathscr{M}$ to an element of $J_{\r}^m[\d]$ by 
adjoining the module 
$\m_m= \m _{m-1}\oplus \r h_{g_m}$
to the end of 
$\mathscr{M}$ (of course there are restrictions on 
the $\lambda$ coefficients which are addressed below).  The other 
way around, every element of $J_{\r}^m[\d]$ can be obtained
by this procedure from flags of length $m-1$. 

To ensure that $\upsilon(\m_m)=\Delta_m$ and that 
$\m_{m-1}\subset \m_m$ 
we need to use Lemma \ref{PIO} above and Lemma
\ref{containment} below. The hypotheses for these lemmas only 
concern 
$\m_{m-1}$ and $g_m$. For this reason, we can reduce the argument
to the case of $m=1$.

From now on we set $D=D_0$, so it corresponds to an 
admissible $\Gamma$-module $\Delta=\Delta_0$. Accordingly, set 
$g=g_1$ and $D_1=D\cup\{g\}$ where $D_1$ will be assumed 
admissible. Also, we let $\Delta'=\Delta\cup\{g\}$.
 
Recall that we let the module generated by the  
$\{m_i\}$ as in (\ref{gen}) be denoted $\m$. Then we consider 
adding 
\[
h\equal z^g+\sum_{k\in [g,\infty)\setminus \Delta}
\lambda_{k-g}^hz^k 
\]
to the set of generators $\{m_i\}$ and  set
$\m'=\langle \m,h\rangle$. Note that $h$ will replace $m_{\ell}$ where
$g\equiv \ell\text{ mod }4$.
 
The modules $\m$ and $\m'$ must satisfy $\upsilon(\m)=\Delta$ and 
$\upsilon(\m')=\Delta'$. 
By (\ref{PIO}), this holds for $\m$ if and only if 
the following elements have higher order expressions: 

\begin{align}\label{REL}
T^1\equal(z^{2u}+z^v)m_0-z^{4\alpha_2}m_2,\  
T^2\equal z^{4(u-\alpha_2)}m_0-(z^{2u}+z^v)m_2,\notag \\
T^3\equal \left(z^{2u+v}+\tfrac{1}{2}z^{2v}\right)m_0
-z^{4\alpha_1}m_j,
\end{align}
where $\alpha_i$ is the unique integer such that 
$a_i=\beta_i(2u+v)+\gamma_i(2u)-4\alpha_i$ for
$\beta_i,\gamma_i\in\{0,1\}$ and $j\equiv 2u+v\mod 4$. To obtain 
the higher order
expressions, the reduction procedure is applied to $T^i$. As a 
result of the 
reduction procedure polynomial relations among 
the $\lambda$ variables are obtained (as discussed before the 
beginning of the proof). For $\m'$ we use a similar approach. 
When we consider the syzygies
(\ref{REL}) in the context of $\m'$ we will denote the
resulting $T$  by 
$(T^i)'$. See pages 14-17 of \cite{Pi} concerning the 
existence of the higher order expressions. 

Given a module $\m$ with $\upsilon(\m)=\Delta$,
not every module $\n$ with $\upsilon(\n)=\Delta'$ is its
extension. In order to understand when $\m\subset \n$, let
\begin{align}
F\equal m_i-(z^4)h\in \n,\notag
\end{align}
where $i\equiv g\mod 4$. Note that $\upsilon(F)>a_i$, 
which gives a new type of syzygy
(not from \cite{Pi}); there is one such syzygy
for each pair $\m\subset \n$. 

Let us take $i$ such that $g\equiv i \hbox{\, mod\,} 4$. 
Then replacing  $a_i$ 
with $g$ results in a basis for $\Delta'$; let us check this.
If there were $\ell$ such that
$\ell \equiv a_i\mod 4$, $g<\ell<a_i$ and $\ell\not\in\Delta$, 
then the presence of $g$ in $\Delta'$ implies that 
$\ell\in \Delta'$ since $\Delta'$ is a $\Gamma$-module. This 
contradicts to $\Delta'=\Delta\cup\{g\}$. Thus $g+4=a_i$,
which gives the required. 
\smallskip

The following lemma provides necessary and sufficient
conditions for $\m\subset\n$ when $\upsilon(\n)=\Delta'$ and 
$\upsilon(\m)=\Delta$. It is important to distinguish performing the 
reduction procedure with respect to $\m$ or $\n$ and we will 
use $\dag_1$ to denote reduction with respect to the generators 
(\ref{gen}) of  $\m$ and $\dag_2$ for reduction with respect to 
(\ref{gen}) except with the changes necessary for $\n$.

\begin{lemma}\label{containment}
Suppose $\m$ and $\n$ are standard $\r$-modules  
in $\o$ such that $\upsilon(\m)=\Delta$ and
$\upsilon(\n)=\Delta\cup\{g\}$. Furthermore, suppose that $\n$ 
contains the 
generators $m_j$ of $\m$ from (\ref{gen}) satisfying
$g\not\equiv j\mod 4$. 
Then $\m\subset \n$ if and only if $F^{\dag_2}=0$
\end{lemma}
{\it Proof.}
Since $g\in \upsilon(\n)$ we see that $h\in \n$ for some choice of 
values for the variables $\lambda^h_k$ and that $h$ is a 
normalized 
generator of $\n$. For some $i$, $g\equiv a_i\mod 4$ (i.e.
$g$ replaces $a_i$ in the basis for $\upsilon(\m)$). 
Observe that $M\subset \n$ if and only if $m_i\in \n$ which will 
happen if and only if $(m_i)^{\dag_2}=0$. $F$ is the first step 
in the reduction of $m_i$ and so $(m_i)^{\dag_2}=0$ if and only 
if $F^{\dag_2}=0$.
\sq

Let us now return to considering  $\m\subset\m'$ with 
$\upsilon(\m)=\Delta$ and $\upsilon(\m')=\Delta'$. We only need 
to consider the equations resulting from $F^{\dag_2}$ since 
the equations resulting from $\upsilon(\m)=\Delta$ and 
$\upsilon(\m')=\Delta'$ are already solved 
for in \cite{Pi}. We may write $F^{\dag_2}=\sum_{k=1}^{\infty}
\widetilde{c}_kz^{a_i+k}$. Recall that the only powers of $z$ present in 
$F^{\dag_2}$ are those greater than $a_i$ that are not in $\Delta$. 
By (\ref{containment}), we have $F^{\dag_2}=0$, which implies 
$\widetilde{c}_k=0$. 
Similar to the analysis of $(T^j)^{\dag_1}$ from the discussion
before the proof, the 
$\widetilde{c}_k$ are in the form
\[
\widetilde{c}_k=\lambda^i_k-\lambda^h_k+
(\text{a polynomial in terms of }\lambda_p^{\bullet} 
\text{ for } p<k).
\]

For a given $k$, we can then express $\lambda^h_k$ in terms of  
$\lambda_p^{\bullet}$ for $p<k$. This gives that $J^1_{\r}[D,g]$ 
is an affine space because it is a graph of a regular function 
on an affine space. Since the $\gamma_{\Delta'}(g+4)$ equations
$\widetilde{c}_k=0$ are solvable, we see that $\mu_{D,g}\leq 
\gamma_{\Delta'}(g)- \gamma_{\Delta'}(g+4)$. The exact value of 
$\mu_{D,g}$  depends on the congruence class of $g$ modulo 4. 
We will now obtain the formulas for the dimensions.

\subsection{\bf Calculating dimensions}
We are going now to justify the dimension formulas
in Theorem \ref{dimprop}. Recall that our approach extends 
the formulas and techniques used in \cite{Pi} to the case
of flags of modules. 
 
First assume $g$ is odd. Since $2u+v$ is odd, 
it is either congruent to $1$ or $3\mod 4$. 
If $g\not\equiv 2u+v\mod 4$, then $T^i=(T^i)'$ for all $i$ and 
hence we do not need to impose any further relations on the 
coefficients. Thus
$\mu_{D,g}=\gamma_{\Delta'}(g)-\gamma_{\Delta'}(g+4)$
in this case.
\smallskip

When $g\equiv 2u+v\mod 4$, we have $(T^3)'=
z^{2u+v}g_0-z^{4\alpha_h}h$ 
and $T^i=(T^i)'$ for $i\not=3$. At the bottoms of 
page 15 and 16 of  \cite{Pi} it is shown that a higher order 
expression exists for $T^3$ when the smallest odd number 
$n\in \Delta\cap [2u,\infty)$ 
is less than or equal to $v$. When $n>v$ it is also shown that 
we can use the higher 
order expressions of $T^1$ and $T^2$ to obtain a higher 
order expression for 
$T^3$ (without imposing any new relations among the 
$\lambda$\~parameters). Hence 
$\mu_{D,g}=\gamma_{\Delta'}(g)-\gamma_{\Delta'}(g+4)$. 
This finishes the proof when $g\equiv 1$ or $3\mod 4$.

In the last case, $g\equiv 2\hbox{ mod }4$\,  implies that
\begin{align}\label{REL2}
\begin{split}
(T^1)'&=(z^{2u}+z^{v})m_0-z^{4\alpha_h}h, \\
(T^2)'&=z^{4(u-\alpha_h)}m_0-(z^{2u}+z^{v})h,
\end{split}
\end{align}
and $T^3=(T^3)'$. Following \cite{Pi}, 
$\gamma_{\Delta'}(2u)$ 
equations for the $\la${'\small s}
result from 
the coefficients of $((T^1)')^{\dag_2}$ and 
$\gamma_{\Delta'}
(g+n)$ distinct equations result from 
the coefficients of  $((T^2)')^{\dag_2}$. We claim 
that the  $\gamma_{\Delta'}(2u)$ equations from 
$((T^1)')^{\dag_2}$ 
are equal to those from $(T^1)^{\dag_1}$ and all of the
$\gamma_{\Delta'}(a_2+n)$ 
equations from $(T^2)^{\dag_1}$ can be obtained from the 
coefficients of $((T^2)')^{\dag_2}$. To prove this we introduce 
the following definition and lemma.

Let $P$ be any power series in $z$ whose 
coefficients are polynomials in 
the $\lambda$ variables. We let 
$\mathfrak{I}(P)$ be the ideal generated by 
the coefficients of $P$ 
in the polynomial ring over the $\lambda$ variables. We have the 
following basic result concerning $\mathfrak{I}$ and $\dag_2$.

\begin{lemma}\label{ideals}
$\mathfrak{I}((rP)^{\dag_2})\subset \mathfrak{I}(P^{\dag_2})$ 
where $r$ is a polynomial in $\r$.
\end{lemma}
{\it Proof.}
Since $\dag_2$ is linear it is sufficient to prove the lemma 
when $r$ is a monomial. The reduction procedure for $rP$ is 
exactly the same as for $P$ until the first $k$ such that 
$k\not\in \Delta \ni k+\upsilon(r)$. Beyond this range, a 
multiple of the 
coefficient of 
$z^k$ in $P^{\dag_2}$ may be added to the 
remaining
coefficients of $(rP)^{\dag_2}$. Because of such $k$,
all coefficients of 
$(rP)^{\dag_2}$
will be coefficients of $P^{\dag_2}$ 
plus some multiples of the previous coefficients
of $P^{\dag_2}$. Thus we have 
$\mathfrak{I}[(rP)^{\dag_2}]\subset 
\mathfrak{I}[P^{\dag_2}]$.
\sq

On page 14 of \cite{Pi} we see that the valuations of 
$T^1, T^2, (T^1)'$ and $(T^2)'$ are greater than $2u$
and therefore greater than $g$. Let us use this.
\begin{lemma}\label{dag}
If $\m\subset \m'$, then 
$(T^1)^{\dag_1}\!=\!(T^1)^{\dag_2}$ and 
$(T^2)^{\dag_1}\!=\!(T^2)^{\dag_2}$.
\end{lemma}
{\it Proof.}
Observe that 
$(T^1)^{\dag_1}-(T^1)^{\dag_2}=r_0m_0+r_1m_1+r_2F+r_3m_3$ 
where the $r_i\in R$.
When we apply $\dag_2$ to the left hand side we get
\[
((T^1)^{\dag_1}-(T^1)^{\dag_2})^{\dag_2}=
(T^1)^{\dag_1}-(T^1)^{\dag_2}
\]
since $(T^1)^{\dag_1}$ has valuations greater
than $g$ which implies the left hand side is an 
eigenvector of $\dag_2$. Applying 
$\dag_2$ to the right hand side we have 
\[
(r_0m_0+r_1m_1+r_2F+r_3m_3)^{\dag_2}=(r_2F)^{\dag_2}
\]
since $m_0,m_1,m_3$ are in $\m'$. Hence we have
\[
(T^1)^{\dag_1}-(T^1)^{\dag_2}=(r_2F)^{\dag_2}.
\]
Note that $r_2$ may be a power series but only a 
truncation of it determines $(r_2F)^{\dag_2}$ because terms 
with valuation greater than the conductor 
will all eventually be eliminated. Therefore by 
Lemma \ref{ideals} we have 
$\mathfrak{I}((r_2F)^{\dag_2})\subset \mathfrak{I}(F^{\dag_2})$. 
Now, $\m\subset \m'$ which means $F^{\dag_2}=0$ by 
Lemma \ref{containment} and so we have $(r_2F)^{\dag_2}=0$. 
The proof for $T^2$ is identical because $T^2$ has valuation 
greater than $2u$.
\sq

Now we prove that $(T^1)^{\dag_1}=((T^1)')^{\dag_2}$ 
when $\m\subset \m'$. Notice that $T^1-(T^1)'=-z^{4\alpha_2}F$, 
where $\alpha_2$, $\alpha_h$ are defined in (\ref{REL}), 
(\ref{REL2}) respectively and we have used 
that $\alpha_h=\alpha_2+1$. Thus
\[
(T^1)^{\dag_2}-((T^1)')^{\dag_2}=(-z^{4\alpha_2}F)^{\dag_2}.
\]
By Lemmas \ref{containment} and \ref{ideals} we have 
$(-z^{4\alpha_2}F)^{\dag_2}=0$ so by the previous Lemma 
\ref{dag} we have 
$(T^1)^{\dag_1}=((T^1)')^{\dag_2}$.

To prove the equations from $(T^2)^{\dag_1}$ are 
redundant we first observe that $T^2-z^4(T^2)'=-(z^{2u}+z^v)F$. 
This implies that
\[
(T^2)^{\dag_2}-(z^4(T^2)')^{\dag_2}=(-(z^{2u}+z^v)F)^{\dag_2}.
\]
Again, by Lemmas \ref{containment} and \ref{ideals} we see that 
$(-(z^{2u}+z^v)F)^{\dag_2}=0$. By Lemma \ref{dag}, we have 
$(T^2)^{\dag_2}= (T^2)^{\dag_1}$ which means 
$(z^4(T^2)')^{\dag_2}=(T^2)^{\dag_1}$. Finally, 
$\mathfrak{I}(((z^4(T^2)')^{\dag_2})\subset 
\mathfrak{I}(((T^2)')^{\dag_2})$ 
by Lemma \ref{ideals} which readily implies that the equations from 
$(T^2)^{\dag_1}$ are redundant.

Thus we have shown that the $\gamma_{\Delta'}(2u)+
\gamma_{\Delta'}(a_2+n)$ 
equations from $(T^1)^{\dag_1}$ and $(T^2)^{\dag_1}$ are 
actually redundant. So, as required:
\begin{align*}
&\mu_{D,g}\!=\!\gamma_{\Delta'}(g)\!-\!\gamma_{\Delta'}(a_2)\!-\!
\gamma_{\Delta'}(2u)\!-\!\gamma_{\Delta'}(g\!+\!n)\!+\! 
\gamma_{\Delta'}(2u)\!+\!\gamma_{\Delta'}(a_2\!+\!n) \\
&=\gamma_{\Delta'}(g)-\gamma_{\Delta'}(a_2)-
(\gamma_{\Delta'}(g\!+\!n)-\gamma_{\Delta'}(a_2\!+\!n)) \\
&=\gamma_{\Delta\cup \{g\}}(g)-\gamma_{\Delta\cup \{g\}}(g\!+\!4)
-(\gamma_{\Delta\cup \{g\}}(g\!+\!n)-\gamma_{\Delta\cup \{g\}}
(g\!+\!4\!+\!n)).
\end{align*}

\setcounter{equation}{0}
\section{\sc Numerical support}
We provide here dimensions of cells for some
basic examples and the corresponding non-admissible 
$D$\~flags $\d$. It is important to know whether these
dimensions, non-admissible $\d$ and admissible
ones with non-affine $J^m_r[\d]$ are {\em topological\,}
properties of the singularity. We found no counterexamples,
but formally this can be wrong even if our Main Conjecture 
holds. Zariski proved in Chapter IV, Section 3 of \cite{Za} that
$\Ga=\lan 4,6,13+2v\ran$ for $v\ge 0$ uniquely determines
the corresponding {\em analytic\,} singularity
(i.e. that each equisingularity class is one point),
but generally such questions can be (very) involved. Also,
these examples seem absolutely necessary for (restarting) 
the theory of Jacobian factors beyond torus knots, which is
of obvious importance for topology and geometry of
plane curve singularities and for orbital integrals.

\subsection{\bf Two simplest cables}
We begin with  formula (\ref{conjaff}) for the knots 
$C\!ab(13,2)T(3,2)$ and $C\!ab(15,2)T(3,2)$.
Here one can use the general Theorem \ref{dimprop} 
or the special Lemma \ref{ADD1GAP}. The
latter can be extended to any $m$ for these
two cases (with minor adjustments).

Recall that 
$\r\!=\!\lan z^4,z^6\!+\!z^7\ran, \Gamma\!=\!
\langle 4,6,13\rangle$ and  
$\r\!=\!\lan z^4,z^6\!+\!z^9\ran, \Gamma\!=
\!\langle 4,6,15\rangle$ 
in these cases.
The first DAHA superpolynomial is as follows:
$$
\vec\rr=\{3,2\},\,
\vec\ss=\{2,1\},\,
\t= C\!ab(13,2)T(3,2);\ 
\h_{\,\vec\rr,\,\vec\ss}\,(\square\,;\,\,q,t,a)=
$$
\renewcommand{\baselinestretch}{0.5} 
{\small
\(
1 + q t + q^8 t^8 + q^2 \bigl(t + t^2\bigr) 
+ a^3 \bigl(q^6 + q^7 t + q^8 t^2\bigr) + 
 q^3 \bigl(t + t^2 + t^3\bigr) + q^4 \bigl(2 t^2 + t^3 + t^4\bigr) 
+  q^5 \bigl(2 t^3 + t^4 + t^5\bigr) + q^6 \bigl(2 t^4 + t^5 
+ t^6\bigr) + 
 q^7 \bigl(t^5 + t^6 + t^7\bigr) + 
 a^2 \bigl(q^3 + q^4 \bigl(1 + t\bigr) + q^5 \bigl(1 + 2 t 
+ t^2\bigr) + 
    q^6 \bigl(2 t + 2 t^2 + t^3\bigr) + q^7 \bigl(2 t^2 + 2 t^3 
+ t^4\bigr) + 
    q^8 \bigl(t^3 + t^4 + t^5\bigr)\bigr) + 
 a \bigl(q + q^2 \bigl(1 + t\bigr) + q^3 \bigl(1 + 2 t + t^2\bigr) 
+ q^4 \bigl(3 t + 2 t^2 + t^3\bigr) + 
    q^5 \bigl(t + 4 t^2 + 2 t^3 + t^4\bigr) + q^6 \bigl(t^2 
+ 4 t^3 + 2 t^4 + t^5\bigr) + 
    q^7 \bigl(t^3 + 3 t^4 + 2 t^5 + t^6\bigr) + q^8 \bigl(t^5 
+ t^6 + t^7\bigr)\bigr).
\)
}
\renewcommand{\baselinestretch}{1.2} 

Here and further we use \cite{ChD1}.
Let us list the necessary information to verify 
(\ref{conjaff}). For the greatest possible
$m=3$ there are only $3$ 
admissible $D$\~sets
$D_0$ that can occur in such a (long) flag. Namely,
these flags and the dimensions dim\,$J^3_{\r}[\d]$ are:

\begin{align}\label{13-3} 
&D_0=[9,11,15],\, \vec{g}= (2,5,7),\, \hbox{dim}=8 
&\rightsquigarrow q^6t^0a^3,\\
&D_0=[7,9,11,15],\, \vec{g}=(2,3,5),\, \hbox{dim}=7
&\rightsquigarrow q^7t^1a^3,\notag\\
&D_0=[5,7,9,11,15],\, \vec{g}=(1,2,3),\, \hbox{dim}=6
&\rightsquigarrow q^8t^2a^3;\notag
\end{align}
we show their contributions to the 
corresponding superpolynomial.


{\footnotesize
\begin{table}[ht!]
\[
\centering
\begin{tabular}{|l|l|}
 \hline 
\hbox{$D$-sets} & $dim$\\
\hline

$\varnothing$ & 8\\
15 & 7\\
11,15 & 6\\
7,11,15 & 6\\
9,15 & 7\\
9,11,15 & 5\\
7,9,11,15 & 4\\
3,7,9,11,15 & 4\\
5,9,11,15 & 5\\
5,7,9,11,15 & 3\\
3,5,7,9,11,15 & 2\\
1,5,7,9,11,15 & 4\\
\hline
\end{tabular}
\hspace{0.1cm}
\begin{tabular}{|l|l|}
\hline
\hbox{$D$-sets} & $dim$\\
\hline
1,3,5,7,9,11,15 & 2\\
2,7,11,15 & 6\\
2,9,15 & 7\\
2,9,11,15 & 6\\
2,7,9,11,15 & 5\\
2,3,7,9,11,15 & 4\\
2,5,9,11,15 & 5\\
2,5,7,9,11,15 & 3\\
2,3,5,7,9,11,15 & 1\\
1,2,5,7,9,11,15 & 3\\
1,2,3,5,7,9,11,15 & 0\\
&\\
\hline   
\end{tabular}
\]
\caption{Dimensions for $\Ga=\lan 4,6,13\ran, m=0$}
\label{Table4-6-13-0}
\end{table}
}

{\small
\begin{table}[ht!]
\[
\centering
\begin{tabular}{|l|l|l|}
\hline 
\hbox{$D$-sets} & $g$ & $dim$\\
\hline
$\varnothing$ & 15 & 8\\
15 & 9 & 8\\
15 & 11 & 7\\
11,15 & 7 & 7\\
11,15 & 9 & 6\\
7,11,15 & 2 & 7\\
7,11,15 & 9 & 6\\
9,15 & 2 & 8\\
9,15 & 11 & 7\\
9,11,15 & 2 & 7\\
9,11,15 & 5 & 6\\
9,11,15 & 7 & 5\\
7,9,11,15 & 2 & 6\\
7,9,11,15 & 3 & 5\\
7,9,11,15 & 5 & 4\\
3,7,9,11,15 & 2 & 5\\
3,7,9,11,15 & 5 & 4\\
5,9,11,15 & 2 & 6\\
5,9,11,15 & 7 & 5\\
5,7,9,11,15 & 1 & 5\\
\hline
\end{tabular}
\hspace{0.1cm}
\begin{tabular}{|l|l|l|}
\hline
\hbox{$D$-sets} & $g$ & $dim$\\
\hline
5,7,9,11,15 & 2 & 4\\
5,7,9,11,15 & 3 & 3\\
3,5,7,9,11,15 & 1 & 3\\
3,5,7,9,11,15 & 2 & 2\\
1,5,7,9,11,15 & 2 & 5\\
1,5,7,9,11,15 & 3 & 4\\
1,3,5,7,9,11,15 & 2 & 2\\
2,7,11,15 & 9 & 6\\
2,9,15 & 11 & 7\\
2,9,11,15 & 5 & 7\\
2,9,11,15 & 7 & 6\\
2,7,9,11,15 & 3 & 6\\
2,7,9,11,15 & 5 & 5\\
2,3,7,9,11,15 & 5 & 4\\
2,5,9,11,15 & 7 & 5\\
2,5,7,9,11,15 & 1 & 4\\
2,5,7,9,11,15 & 3 & 3\\
2,3,5,7,9,11,15 & 1 & 1\\
1,2,5,7,9,11,15 & 3 & 3\\
&&\\
\hline   
\end{tabular}
\]
\caption{Dimensions for $\Ga=\lan 4,6,13\ran, m=1$}
\label{Table4-6-13-1}
\end{table}
}
Tables \ref{Table4-6-13-0}, \ref{Table4-6-13-1},
\ref{Table4-6-13-2}
show $D_0$, the corresponding $\vec{g}$
and the dimensions of $J^m_{\r}[\d]$
for all admissible flags as $m=0,1,2$.

{\small
\begin{table}[ht!]
\[
\centering
\begin{tabular}{|l|l|l|}
\hline 
\hbox{$D$-sets} & $\vec{g}$ & $dim$\\
\hline
15 & 9,11 & 8\\
11,15 & 7,9 & 7\\
7,11,15 & 2,9 & 7\\
9,15 & 2,11 & 8\\
9,11,15 & 2,5 & 8\\
9,11,15 & 2,7 & 7\\
9,11,15 & 5,7 & 6\\
7,9,11,15 & 2,3 & 7\\
7,9,11,15 & 2,5 & 6\\
7,9,11,15 & 3,5 & 5\\
\hline
\end{tabular}
\hspace{0.1cm}
\begin{tabular}{|l|l|l|}
\hline
\hbox{$D$-sets} & $\vec{g}$ & $dim$\\
\hline
3,7,9,11,15 & 2,5 & 5\\
5,9,11,15 & 2,7 & 6\\
5,7,9,11,15 & 1,2 & 6\\
5,7,9,11,15 & 1,3 & 5\\
5,7,9,11,15 & 2,3 & 4\\
3,5,7,9,11,15 & 1,2 & 3\\
1,5,7,9,11,15 & 2,3 & 5\\
2,9,11,15 & 5,7 & 7\\
2,7,9,11,15 & 3,5 & 6\\
2,5,7,9,11,15 & 1,3 & 4\\
\hline   
\end{tabular}
\]
\caption{Dimensions for $\Ga=\lan 4,6,13\ran, m=2$}
\label{Table4-6-13-2}
\end{table}
}

\comment{
{\small
\begin{table}[ht!]
$$
\centering
\begin{tabular}{|l|l|l|}
\hline 
\hbox{$D$-sets} & $\vec{g}$ & $dim$\\
\hline
9,11,15 & 2,5,7 & 8\\
7,9,11,15 & 2,3,5 & 7\\
5,7,9,11,15 & 1,2,3 & 6\\
\hline   
\end{tabular}
$$
\caption{Dimensions for $\Ga=\lan 4,6,13\ran, m=3$}
\label{Table4-6-13-3}
\end{table}
}
}

There are no admissible extensions of degree $4$, so we have
$$
\sum_{D} q^{|D|+m\,} t^{\de-\hbox{\tiny dim}
(J^m_{\r}[D,D'])\,} a^m=
$$
\renewcommand{\baselinestretch}{0.5} 
{\small
\(
1 + q t + q^8 t^8 + q^2 \bigl(t + t^2\bigr) 
+ a^3 \bigl(q^6 + q^7 t + q^8 t^2\bigr) + 
 q^3 \bigl(t + t^2 + t^3\bigr) + q^4 \bigl(2 t^2 + t^3 
+ t^4\bigr) + 
 q^5 \bigl(2 t^3 + t^4 + t^5\bigr) + q^6 \bigl(2 t^4 + t^5 
+ t^6\bigr) + 
 q^7 \bigl(t^5 + t^6 + t^7\bigr) + 
 a^2 \bigl(q^3 + q^4 \bigl(1 + t\bigr) + q^5 \bigl(1 + 2 t 
+ t^2\bigr) + 
    q^6 \bigl(2 t + 2 t^2 + t^3\bigr) + q^7 \bigl(2 t^2 + 2 t^3 
+ t^4\bigr) + 
    q^8 \bigl(t^3 + t^4 + t^5\bigr)\bigr) + 
 a \bigl(q + q^2 \bigl(1 + t\bigr) + q^3 \bigl(1 + 2 t + t^2\bigr) 
+ q^4 \bigl(3 t + 2 t^2 + t^3\bigr) + 
    q^5 \bigl(t + 4 t^2 + 2 t^3 + t^4\bigr) + q^6 \bigl(t^2 
+ 4 t^3 + 2 t^4 + t^5\bigr) + 
    q^7 \bigl(t^3 + 3 t^4 + 2 t^5 + t^6\bigr) + q^8 \bigl(t^5 
+ t^6 + t^7\bigr)\bigr),
\)
}
\renewcommand{\baselinestretch}{1.2} 
\smallskip

\noindent which coincides with $\h_{\,\vec\rr,\,\vec\ss}\,
(\square\,;\,\,q,t,a)$ from Section 3.1 of \cite{ChD1}.
\smallskip

For $\r=\C[[ t^4, t^6+t^9]]$ corresponding to the
$\Ga=\lan 4,6,15\ran$ and
 cable $C\!ab(15,2)T(3,2)$ the situation is very similar. We
checked that:
\begin{align*}
&\sum_{m=0}^\infty \,
\sum_{\{D_0=D,\ldots,D_m\}} q^{|D|+m\,} 
t^{\de-\hbox{\tiny dim}
(J^m_{\r}[D])\,} a^m=
\h_{\,\{3,2\},\{2,3\}}(\square\,;\,\,q,t,a)=
\end{align*}

\renewcommand{\baselinestretch}{0.5} 
{\small
\(
1 + q t + q^9 t^9 + q^2 \bigl(t + t^2\bigr) + q^3 \bigl(t + t^2 
+ t^3\bigr) + 
 a^3 \bigl(q^6 + q^7 t + q^8 t^2 + q^9 t^3\bigr) 
+ q^4 \bigl(2 t^2 + t^3 + t^4\bigr) + 
 q^5 \bigl(2 t^3 + t^4 + t^5\bigr) + q^6 \bigl(2 t^4 + t^5 
+ t^6\bigr) + 
 q^7 \bigl(2 t^5 + t^6 + t^7\bigr) + q^8 \bigl(t^6 + t^7 
+ t^8\bigr) + 
 a^2 \bigl(q^3 + q^4 \bigl(1 + t\bigr) + q^5 \bigl(1 + 2 t 
+ t^2\bigr) + 
    q^6 \bigl(2 t + 2 t^2 + t^3\bigr) + q^7 \bigl(2 t^2 + 2 t^3 
+ t^4\bigr) + 
    q^8 \bigl(2 t^3 + 2 t^4 + t^5\bigr) + q^9 \bigl(t^4 + t^5 
+ t^6\bigr)\bigr) + 
 a \bigl(q + q^2 \bigl(1 + t\bigr) + q^3 \bigl(1 + 2 t + t^2\bigr)
 + q^4 \bigl(3 t + 2 t^2 + t^3\bigr) + 
    q^5 \bigl(t + 4 t^2 + 2 t^3 + t^4\bigr) + q^6 \bigl(t^2 
+ 4 t^3 + 2 t^4 + t^5\bigr) + 
    q^7 \bigl(t^3 + 4 t^4 + 2 t^5 + t^6\bigr) + 
    q^8 \bigl(t^4 + 3 t^5 + 2 t^6 + t^7\bigr) + q^9 \bigl(t^6 
+ t^7 + t^8\bigr)\bigr).
\)
}
\renewcommand{\baselinestretch}{1.2} 

\smallskip
\noindent
See the same section of \cite{ChD1}. Note that formulas
(3.1) and (3.2) there for the DAHA-Betti polynomials are
obtained from the superpolynomials as $a=0,q=1$. 

\subsection{\bf The case of 
\texorpdfstring{\mathversion{bold}$\Ga\!=\!\lan4,14,31\ran$}
{<4,14,31>}}
We checked the Main Conjecture in
many cases for the series 
$(4,2u,v)$ of {\em Puiseux exponents}. 
This example is of importance because quite a few features of
Theorem \ref{dimprop} (and our proof) cannot be seen for the
sub-family $(4,6,v)$.  
For $\Ga=\lan4,14,31\ran$, 
we list all admissible $D$\~flags $\d$ and the  
dimensions dim\,$J^m_{\r}[\d]$\, for $m=3$ (the greatest 
possible value). The calculation of dimensions
at the maximal $m$ as a matter of fact includes a lot of 
information about the dimensions for previous lengths $m$,
so this is a good test of our conjecture. However we restrict
ourselves with $a^3$ here due to practical
reasons.  The total number of admissible $\d$ is 1071, but there
are ``only" 85 such $D$\~flags of the top length for $m=3$.
The dimensions are in Tables \ref{Table4-14-31-1},
\ref{Table4-14-31-2} (its continuation) in the same format 
as above. 
Compare them with the coefficients of $a^3$ of the corresponding
superpolynomial from \cite{ChD1}, which we provide:

\vskip -0.5cm
$$\r=\C[[ z^4,z^{14}+z^{17}]],\ \t=C\!ab(31,2)T(7,2),\ 
\h_{\{7,2\},\{2,3\}}(\square;q,t,a)=$$

\renewcommand{\baselinestretch}{0.5} 
{\small
\(
1+q t+q^2 t+q^3 t+q^2 t^2+q^3 t^2+2 q^4 t^2+q^5 t^2+q^6 t^2+q^3 t^3
+q^4 t^3+2 q^5 t^3+2 q^6 t^3+2 q^7 t^3+q^8 t^3
+q^9 t^3+q^4 t^4+q^5 t^4+2 q^6 t^4+2 q^7 t^4+3 q^8 t^4+2 q^9 t^4
+2 q^{10} t^4+q^5 t^5+q^6 t^5+2 q^7 t^5+2 q^8 t^5+3 q^9 t^5
+3 q^{10} t^5+3 q^{11} t^5+q^6 t^6+q^7 t^6+2 q^8 t^6+2 q^9 t^6
+3 q^{10} t^6+3 q^{11} t^6+4 q^{12} t^6+q^7 t^7+q^8 t^7+2 q^9 t^7
+2 q^{10} t^7+3 q^{11} t^7+3 q^{12} t^7+4 q^{13} t^7+q^8 t^8
+q^9 t^8+2 q^{10} t^8+2 q^{11} t^8+3 q^{12} t^8+3 q^{13} t^8
+4 q^{14} t^8+q^9 t^9+q^{10} t^9+2 q^{11} t^9+2 q^{12} t^9
+3 q^{13} t^9+3 q^{14} t^9+4 q^{15} t^9+q^{10} t^{10}+q^{11} t^{10}
+2 q^{12} t^{10}+2 q^{13} t^{10}+3 q^{14} t^{10}+3 q^{15} t^{10}
+3 q^{16} t^{10}+q^{11} t^{11}+q^{12} t^{11}+2 q^{13} t^{11}
+2 q^{14} t^{11}+3 q^{15} t^{11}+3 q^{16} t^{11}+2 q^{17} t^{11}
+q^{12} t^{12}+q^{13} t^{12}+2 q^{14} t^{12}+2 q^{15} t^{12}
+3 q^{16} t^{12}+2 q^{17} t^{12}+q^{18} t^{12}+q^{13} t^{13}
+q^{14} t^{13}+2 q^{15} t^{13}+2 q^{16} t^{13}+3 q^{17} t^{13}
+q^{18} t^{13}+q^{14} t^{14}+q^{15} t^{14}+2 q^{16} t^{14}
+2 q^{17} t^{14}+2 q^{18} t^{14}+q^{15} t^{15}+q^{16} t^{15}
+2 q^{17} t^{15}+2 q^{18} t^{15}+q^{19} t^{15}+q^{16} t^{16}
+q^{17} t^{16}+2 q^{18} t^{16}+q^{19} t^{16}+q^{17} t^{17}
+q^{18} t^{17}+2 q^{19} t^{17}+q^{18} t^{18}+q^{19} t^{18}
+q^{20} t^{18}+q^{19} t^{19}+q^{20} t^{19}+q^{20} t^{20}
+q^{21} t^{21}
\)

\noindent
\(
+a^3 \bigl(q^6+q^7 t+q^8 t+q^9 t+q^8 t^2+q^9 t^2
+2 q^{10} t^2+q^{11} t^2+q^{12} t^2+q^9 t^3+q^{10} t^3+2 q^{11} t^3
+2 q^{12} t^3+2 q^{13} t^3+q^{10} t^4+q^{11} t^4+2 q^{12} t^4
+2 q^{13} t^4+3 q^{14} t^4+q^{11} t^5+q^{12} t^5+2 q^{13} t^5
+2 q^{14} t^5+3 q^{15} t^5+q^{12} t^6+q^{13} t^6+2 q^{14} t^6
+2 q^{15} t^6+3 q^{16} t^6+q^{13} t^7+q^{14} t^7+2 q^{15} t^7
+2 q^{16} t^7+3 q^{17} t^7+q^{14} t^8+q^{15} t^8+2 q^{16} t^8
+2 q^{17} t^8+2 q^{18} t^8+q^{15} t^9+q^{16} t^9+2 q^{17} t^9
+2 q^{18} t^9+q^{19} t^9+q^{16} t^{10}+q^{17} t^{10}
+2 q^{18} t^{10}+q^{19} t^{10}+q^{17} t^{11}+q^{18} t^{11}
+2 q^{19} t^{11}+q^{18} t^{12}+q^{19} t^{12}+q^{20} t^{12}
+q^{19} t^{13}+q^{20} t^{13}+q^{20} t^{14}+q^{21} t^{15}\bigr)
\)

\noindent
\(
+a^2 \bigl(q^3+q^4+q^5+q^4 t+2 q^5 t+3 q^6 t
+2 q^7 t+q^8 t+q^5 t^2+2 q^6 t^2+4 q^7 t^2+4 q^8 t^2+4 q^9 t^2
+2 q^{10} t^2+q^{11} t^2+q^6 t^3+2 q^7 t^3+4 q^8 t^3+5 q^9 t^3
+6 q^{10} t^3+5 q^{11} t^3+3 q^{12} t^3+q^7 t^4+2 q^8 t^4
+4 q^9 t^4+5 q^{10} t^4+7 q^{11} t^4+7 q^{12} t^4+5 q^{13} t^4
+q^8 t^5+2 q^9 t^5+4 q^{10} t^5+5 q^{11} t^5+7 q^{12} t^5
+8 q^{13} t^5+6 q^{14} t^5+q^9 t^6+2 q^{10} t^6+4 q^{11} t^6
+5 q^{12} t^6+7 q^{13} t^6+8 q^{14} t^6+6 q^{15} t^6+q^{10} t^7
+2 q^{11} t^7+4 q^{12} t^7+5 q^{13} t^7+7 q^{14} t^7+8 q^{15} t^7
+6 q^{16} t^7+q^{11} t^8+2 q^{12} t^8+4 q^{13} t^8+5 q^{14} t^8
+7 q^{15} t^8+8 q^{16} t^8+5 q^{17} t^8+q^{12} t^9+2 q^{13} t^9
+4 q^{14} t^9+5 q^{15} t^9+7 q^{16} t^9+7 q^{17} t^9+3 q^{18} t^9
+q^{13} t^{10}+2 q^{14} t^{10}+4 q^{15} t^{10}+5 q^{16} t^{10}
+7 q^{17} t^{10}+5 q^{18} t^{10}+q^{19} t^{10}+q^{14} t^{11}
+2 q^{15} t^{11}+4 q^{16} t^{11}+5 q^{17} t^{11}+6 q^{18} t^{11}
+2 q^{19} t^{11}+q^{15} t^{12}+2 q^{16} t^{12}+4 q^{17} t^{12}
+5 q^{18} t^{12}+4 q^{19} t^{12}+q^{16} t^{13}+2 q^{17} t^{13}
+4 q^{18} t^{13}+4 q^{19} t^{13}+q^{20} t^{13}+q^{17} t^{14}
+2 q^{18} t^{14}+4 q^{19} t^{14}+2 q^{20} t^{14}+q^{18} t^{15}
+2 q^{19} t^{15}+3 q^{20} t^{15}+q^{19} t^{16}+2 q^{20} t^{16}
+q^{21} t^{16}+q^{20} t^{17}+q^{21} t^{17}+q^{21} t^{18}\bigr)
\)

\noindent
\(
+a \bigl(q+q^2+q^3+q^2 t+2 q^3 t+3 q^4 t+2 q^5 t+q^6 t+q^3 t^2
+2 q^4 t^2+4 q^5 t^2+4 q^6 t^2+4 q^7 t^2+2 q^8 t^2+q^9 t^2+q^4 t^3
+2 q^5 t^3+4 q^6 t^3+5 q^7 t^3+6 q^8 t^3+5 q^9 t^3+4 q^{10} t^3
+q^{11} t^3+q^5 t^4+2 q^6 t^4+4 q^7 t^4+5 q^8 t^4+7 q^9 t^4
+7 q^{10} t^4+7 q^{11} t^4+2 q^{12} t^4+q^6 t^5+2 q^7 t^5
+4 q^8 t^5
+5 q^9 t^5+7 q^{10} t^5+8 q^{11} t^5+9 q^{12} t^5+3 q^{13} t^5
+q^7 t^6+2 q^8 t^6+4 q^9 t^6+5 q^{10} t^6+7 q^{11} t^6
+8 q^{12} t^6
+10 q^{13} t^6+3 q^{14} t^6+q^8 t^7+2 q^9 t^7+4 q^{10} t^7
+5 q^{11} t^7+7 q^{12} t^7+8 q^{13} t^7+10 q^{14} t^7+3 q^{15} t^7
+q^9 t^8+2 q^{10} t^8+4 q^{11} t^8+5 q^{12} t^8+7 q^{13} t^8
+8 q^{14} t^8+10 q^{15} t^8+3 q^{16} t^8+q^{10} t^9+2 q^{11} t^9
+4 q^{12} t^9+5 q^{13} t^9+7 q^{14} t^9+8 q^{15} t^9+9 q^{16} t^9
+2 q^{17} t^9+q^{11} t^{10}+2 q^{12} t^{10}+4 q^{13} t^{10}
+5 q^{14} t^{10}+7 q^{15} t^{10}+8 q^{16} t^{10}+7 q^{17} t^{10}
+q^{18} t^{10}+q^{12} t^{11}+2 q^{13} t^{11}+4 q^{14} t^{11}
+5 q^{15} t^{11}+7 q^{16} t^{11}+7 q^{17} t^{11}
+4 q^{18} t^{11}+q^{13} t^{12}+2 q^{14} t^{12}
+4 q^{15} t^{12}+5 q^{16} t^{12}+7 q^{17} t^{12}
+5 q^{18} t^{12}+q^{19} t^{12}+q^{14} t^{13}
+2 q^{15} t^{13}+4 q^{16} t^{13}+5 q^{17} t^{13}
+6 q^{18} t^{13}+2 q^{19} t^{13}+q^{15} t^{14}
+2 q^{16} t^{14}+4 q^{17} t^{14}+5 q^{18} t^{14}
+4 q^{19} t^{14}+q^{16} t^{15}+2 q^{17} t^{15}
+4 q^{18} t^{15}+4 q^{19} t^{15}+q^{20} t^{15}
+q^{17} t^{16}+2 q^{18} t^{16}+4 q^{19} t^{16}
+2 q^{20} t^{16}+q^{18} t^{17}+2 q^{19} t^{17}
+3 q^{20} t^{17}+q^{19} t^{18}+2 q^{20} t^{18}
+q^{21} t^{18}+q^{20} t^{19}+q^{21} t^{19}+q^{21} t^{20}\bigr).
\)
}
\renewcommand{\baselinestretch}{1.2} 
\smallskip

\renewcommand{\baselinestretch}{1.1} 
{\tiny
\begin{table}[ht!]
$$
\centering
\begin{tabular}{|l|l|l|}
\hline 
\hbox{$D$-sets} & $\vec{g}$ & $dim$\\
\hline
27,37,41 & 10,23,33 & 21\\
27,33,37,41 & 10,23,29 & 20\\
27,29,33,37,41 & 10,23,25 & 19\\
25,27,29,33,37,41 & 10,21,23 & 18\\
21,25,27,29,33,37,41 & 10,17,23 & 18\\
17,21,25,27,29,33,37,41 & 10,13,23 & 18\\
23,27,33,37,41 & 10,19,29 & 20\\
23,27,29,33,37,41 & 10,19,25 & 19\\
23,25,27,29,33,37,41 & 10,19,21 & 17\\
21,23,25,27,29,33,37,41 & 10,17,19 & 16\\
17,21,23,25,27,29,33,37,41 & 10,13,19 & 16\\
13,17,21,23,25,27,29,33,37,41 & 9,10,19 & 16\\
19,23,27,29,33,37,41 & 10,15,25 & 19\\
19,23,25,27,29,33,37,41 & 10,15,21 & 17\\
19,21,23,25,27,29,33,37,41 & 10,15,17 & 15\\
17,19,21,23,25,27,29,33,37,41 & 10,13,15 & 14\\
13,17,19,21,23,25,27,29,33,37,41 & 9,10,15 & 14\\
9,13,17,19,21,23,25,27,29,33,37,41 & 5,10,15 & 15\\
15,19,23,25,27,29,33,37,41 & 10,11,21 & 17\\
15,19,21,23,25,27,29,33,37,41 & 10,11,17 & 15\\
15,17,19,21,23,25,27,29,33,37,41 & 10,11,13 & 13\\
13,15,17,19,21,23,25,27,29,33,37,41 & 9,10,11 & 12\\
9,13,15,17,19,21,23,25,27,29,33,37,41 & 5,10,11 & 13\\
5,9,13,15,17,19,21,23,25,27,29,33,37,41 & 1,10,11 & 14\\
11,15,19,21,23,25,27,29,33,37,41 & 7,10,17 & 15\\
11,15,17,19,21,23,25,27,29,33,37,41 & 7,10,13 & 13\\
11,13,15,17,19,21,23,25,27,29,33,37,41 & 7,9,10 & 11\\
9,11,13,15,17,19,21,23,25,27,29,33,37,41 & 5,7,10 & 11\\
5,9,11,13,15,17,19,21,23,25,27,29,33,37,41 & 1,7,10 & 12\\
7,11,15,17,19,21,23,25,27,29,33,37,41 & 3,10,13 & 14\\
7,11,13,15,17,19,21,23,25,27,29,33,37,41 & 3,9,10 & 12\\
7,9,11,13,15,17,19,21,23,25,27,29,33,37,41 & 3,5,10 & 11\\
5,7,9,11,13,15,17,19,21,23,25,27,29,33,37,41 & 1,3,10 & 11\\
10,17,21,25,27,29,33,37,41 & 6,13,23 & 18\\
10,23,27,33,37,41 & 6,19,29 & 20\\
10,23,27,29,33,37,41 & 6,19,25 & 19\\
10,23,25,27,29,33,37,41 & 6,19,21 & 18\\
10,21,23,25,27,29,33,37,41 & 6,17,19 & 17\\
10,17,21,23,25,27,29,33,37,41 & 6,13,19 & 17\\
10,13,17,21,23,25,27,29,33,37,41 & 6,9,19 & 16\\
10,19,23,27,29,33,37,41 & 6,15,25 & 19\\
10,19,23,25,27,29,33,37,41 & 6,15,21 & 18\\
10,19,21,23,25,27,29,33,37,41 & 6,15,17 & 16\\
\hline
\end{tabular}
$$
\caption{Dimensions for $\Ga=\lan 4,14,31\ran, m=3\, 
\hbox{(I)}$}
\label{Table4-14-31-1}
\end{table}
}
\renewcommand{\baselinestretch}{1.2} 

\renewcommand{\baselinestretch}{1.1} 
{\tiny
\begin{table}[ht!]
$$
\begin{tabular}{|l|l|l|}
\hline
\hbox{$D$-sets} & $\vec{g}$ & $dim$\\
\hline
10,17,19,21,23,25,27,29,33,37,41 & 6,13,15 & 15\\
10,13,17,19,21,23,25,27,29,33,37,41 & 6,9,15 & 14\\
9,10,13,17,19,21,23,25,27,29,33,37,41 & 5,6,15 & 14\\
10,15,19,23,25,27,29,33,37,41 & 6,11,21 & 17\\
10,15,19,21,23,25,27,29,33,37,41 & 6,11,17 & 16\\
10,15,17,19,21,23,25,27,29,33,37,41 & 6,11,13 & 14\\
10,13,15,17,19,21,23,25,27,29,33,37,41 & 6,9,11 & 12\\
9,10,13,15,17,19,21,23,25,27,29,33,37,41 & 5,6,11 & 12\\
5,9,10,13,15,17,19,21,23,25,27,29,33,37,41 & 1,6,11 & 13\\
10,11,15,19,21,23,25,27,29,33,37,41 & 6,7,17 & 15\\
10,11,15,17,19,21,23,25,27,29,33,37,41 & 6,7,13 & 13\\
10,11,13,15,17,19,21,23,25,27,29,33,37,41 & 6,7,9 & 10\\
9,10,11,13,15,17,19,21,23,25,27,29,33,37,41 & 5,6,7 & 9\\
5,9,10,11,13,15,17,19,21,23,25,27,29,33,37,41 & 1,6,7 & 10\\
7,10,11,15,17,19,21,23,25,27,29,33,37,41 & 3,6,13 & 13\\
7,10,11,13,15,17,19,21,23,25,27,29,33,37,41 & 3,6,9 & 10\\
7,9,10,11,13,15,17,19,21,23,25,27,29,33,37,41 & 3,5,6 & 8\\
5,7,9,10,11,13,15,17,19,21,23,25,27,29,33,37,41 & 1,3,6 & 8\\
6,10,17,21,25,27,29,33,37,41 & 2,13,23 & 18\\
6,10,17,21,23,25,27,29,33,37,41 & 2,13,19 & 17\\
6,10,13,17,21,23,25,27,29,33,37,41 & 2,9,19 & 16\\
6,10,19,23,27,29,33,37,41 & 2,15,25 & 19\\
6,10,19,23,25,27,29,33,37,41 & 2,15,21 & 18\\
6,10,19,21,23,25,27,29,33,37,41 & 2,15,17 & 17\\
6,10,17,19,21,23,25,27,29,33,37,41 & 2,13,15 & 16\\
6,10,13,17,19,21,23,25,27,29,33,37,41 & 2,9,15 & 15\\
6,9,10,13,17,19,21,23,25,27,29,33,37,41 & 2,5,15 & 14\\
6,10,15,19,23,25,27,29,33,37,41 & 2,11,21 & 17\\
6,10,15,19,21,23,25,27,29,33,37,41 & 2,11,17 & 16\\
6,10,15,17,19,21,23,25,27,29,33,37,41 & 2,11,13 & 15\\
6,10,13,15,17,19,21,23,25,27,29,33,37,41 & 2,9,11 & 13\\
6,9,10,13,15,17,19,21,23,25,27,29,33,37,41 & 2,5,11 & 12\\
5,6,9,10,13,15,17,19,21,23,25,27,29,33,37,41 & 1,2,11 & 12\\
6,10,11,15,19,21,23,25,27,29,33,37,41 & 2,7,17 & 15\\
6,10,11,15,17,19,21,23,25,27,29,33,37,41 & 2,7,13 & 14\\
6,10,11,13,15,17,19,21,23,25,27,29,33,37,41 & 2,7,9 & 11\\
6,9,10,11,13,15,17,19,21,23,25,27,29,33,37,41 & 2,5,7 & 9\\
5,6,9,10,11,13,15,17,19,21,23,25,27,29,33,37,41 & 1,2,7 & 9\\
6,7,10,11,15,17,19,21,23,25,27,29,33,37,41 & 2,3,13 & 13\\
6,7,10,11,13,15,17,19,21,23,25,27,29,33,37,41 & 2,3,9 & 10\\
6,7,9,10,11,13,15,17,19,21,23,25,27,29,33,37,41 & 2,3,5 & 7\\
5,6,7,9,10,11,13,15,17,19,21,23,25,27,29,33,37,41 & 1,2,3 & 6\\
\hline   
\end{tabular}
$$
\caption{Dimensions for $\Ga=\lan 4,14,31\ran, m=3\ \hbox{(II)}$}
\label{Table4-14-31-2}
\end{table}
}
\subsection{\bf The series 
\texorpdfstring{{\mathversion{bold}$(6,8,v)$}}{(6,8,v)}}
The series with  {\em Puiseux exponents\,} $(4,2u,v)$ above
corresponds to links $C\!ab(2u+v,2)T(u,2)$, which are
somewhat special; torus knots and cables for $(2p+1,2)$
are known to have some special symmetries.

Let us consider the series 
$(6,8,v)$ for odd $v\ge 9$ from (\ref{PiEuler}) with the link 
$C\!ab(25+(v-9),2)T(4,3)$. The corresponding ring, semigroup, 
$\de$\~invariant, and the Euler number of $J_{\r}$ are
$\r=\C[[z^6,z^{8}+z^v]], \Ga=\lan 6,8,25+(v-9)\ran, 
\de =  18+(v-9)/2,\hspace{1mm} e(J_{\r})=227+25(v-9)/2$.

Paper \cite{Pi} provides the Euler number, the total number
of $\Gamma$\~modules, which is $273+25(v-9)/2$, and the
number of {\em non-admissible} ones. The latter is
$46$ for any $v$. The difference 
$273+25(v-9)/2-46$ is exactly the Euler number, since 
all admissible cells in the decomposition from \cite{Pi} are
diffeomorphic to $\mathbb{A}^{\!N}$ in the case under
consideration and the homology can be readily calculated.

For our conjecture, we need to know the set of
all $\Gamma$\~modules $\De$,
those that are non-admissible, and the dimensions
dim\!$(J_{\r}[\d])$. It is
not too difficult to find all non-admissible
modules following \cite{Pi}, but they are not 
provided in his paper.
There are ``generic" non-admissible modules and $3$ 
exceptional ones, which we are going to describe.

We examine the elements 
\[\label{tijpq}
T_{ij}^{pq}=\phi_i m_p-\phi_j m_q\in \m
\]
where $D=D_{\m}$ is the set of gaps of $\m$, 
$i<j\in \Gamma$, $p>q\in \De_{\m}=\upsilon(\m)$,
$\upsilon(\phi_i)=i,\, \upsilon(m_p)=p$,
where we use the valuation $\upsilon: \o\to \Z_+$.
Recall that modules $\m$ are submodules in $\o=\C[[z]]$
with an element of valuation $0$. 

We will assume 
that the leading $z$\~monomial in $\phi_i,m_k$ has the 
coefficient $1$. 
For instance, $m_0=1+\sum_{p>0}\la_0^p z^p$.
The choice of these elements is of course non-unique
(higher terms can be added to them).

\begin{proposition}
All non-admissible $\m$ for $\r=\C[[t^6,t^8+t^v]]$
can be described as follows. In the differences
from (\ref{tijpq}), let:
$$
(a)\  q=0,p=2,4,10; \ 
(b)\  q=2,p=4;\  (c)\ p>q \in  \{0,2,4\}.
$$
Then, let us impose the following relations for $i,j$
in (\ref{tijpq}):
\begin{align}\label{no-admcon}
i+p+1\not\in \De_{\m} \hbox{\, and\,\, } i+p=j+q
\for i,j\in \Gamma.
\end{align} 
Here the latter results in the inequalities
$\upsilon(T_{ij}^{pq})\ge 1+p+1$, which can
be only strict due to the former since 
$T_{ij}^{pq}\in\m$.
 
The non-admissibility of $\m$ of type (a) or (b)
is if and only if there exist
no $m_p,m_q\in \m$ in (\ref{tijpq}) satisfying
(\ref{no-admcon}) 
for {\sf all possible} choices of  $i,j$ there.
In the case of (c), the absence of
$m_0,m_2,m_4\in \m$   must be checked 
for each of the 3 choices
of $p,q$ there altogether and every possible $i,j$
satisfying  (\ref{tijpq}).
\sq
\end{proposition}

Let us list the $D$-sets (the sets of gaps) for all 
$46$ non-admissible modules $\m$.
In Table \ref{Table6-8-25}, $(v-9)$ must be
added to all gaps in the second half of the first
column (clearly visible). The second column
contains the smallest $p>0$ in $D$
that ensures the non-admissibility
in the case of $(a)$, the letters $b$ or $c$ stand for 
the remaining three exceptional cases.
The third column contains the first 
$g=g_{min}=i+p+1\not\in D$
such that its absence in $D$ contradicts the absence
of the previous gaps $g'i'+p+1<g_{min}$ in $D$, which is for a
given pair from $(a,b,c)$. Such a gap is provided only
for the pair $(0,2)$ in the case of $(c)$. 

The non-admissibility can occur only if there are
at least two possible pairs $(i,j)$ satisfying
(\ref{no-admcon}).  The corresponding conditions are simple 
equalities for the differences $\la^1_p-\la^1_q$, where
$m_p=z^p+\la^{1}_p z^{p+1}+\cdots$. If the difference
 $\la^1_p-\la^1_q$ takes different values for different
pairs $(i,j)$, then the module cannot be admissible. 
In the case of
$c$, such a contradiction can be reached only if 
{\em three} such sequences of equalities are considered together,
i.e. for $\la_0^1-\la_2^1$, $\la_0^1-\la_4^1$,
$\la_2^1-\la_4^1$.   
Generally, higher order $z$\~expansions can be necessary
here, i.e. $\la_p^{i>1}$ may occur;
however this is not the case with $(6,8,v)$.

{\tiny
\begin{table}[ht!]
$$
\centering
\begin{tabular}{|ll|l|l|}
  \hline 
\hbox{$D$-sets} & under $\{\cdot\}+(v-9)$ & p & g\\
\hline    
 10,  &  35 &   10 &  19  \\                         
 10,  &  27, 35 &   10 &  19  \\                     
 10,  &  29, 35 &   10 &  19  \\                     
 2, 10,  &  27, 35 &   2 &  19  \\                   
 4, 10,  &  29, 35 &   4 &  17  \\                   
 10,  &  23, 29, 35 &   10 &  19  \\                 
 10,  &  27, 29, 35 &   10 &  19  \\                 
 2, 10,  &  27, 29, 35 &   2 &  19  \\               
 4, 10,  &  23, 29, 35 &   4 &  17  \\               
 4, 10,  &  27, 29, 35 &   4 &  17  \\               
 10,  &  21, 27, 29, 35 &   10 &  19  \\             
 10,  &  23, 27, 29, 35 &   10 &  19  \\             
 2, 4, 10,  &  27, 29, 35 &   2 &  19  \\            
 2, 10,  &  21, 27, 29, 35 &   2 &  19  \\           
 2, 10,  &  23, 27, 29, 35 &   2 &  19  \\           
 4, 10,  &  19, 27, 29, 35 &   4 &  17  \\           
 4, 10,  &  21, 27, 29, 35 &   4 &  17  \\           
 4, 10,  &  23, 27, 29, 35 &   4 &  17 \\             
 10,  &  21, 23, 27, 29, 35 &   10 &  19 \\           
 2, 4, 10,  &  19, 27, 29, 35 &   4 &  17 \\          
 2, 4, 10,  &  21, 27, 29, 35 &   2 &  19 \\          
 2, 4, 10,  &  23, 27, 29, 35 &   2 &  19 \\          
 2, 10,  &  17, 23, 27, 29, 35 &   2 &  19 \\         
 \hline   
\end{tabular}     
\hspace {0.1cm}
\begin{tabular}{|ll|l|l|}
  \hline         
\hbox{$D$-sets} & under $\{\cdot\}+(v-9)$ & p & g\\  
\hline
2, 10,  &  21, 23, 27, 29, 35 &   2 &  19 \\             
4, 10,  &  19, 21, 27, 29, 35 &   4 &  17 \\             
4, 10,  &  19, 23, 27, 29, 35 &   4 &  17 \\             
4, 10,  &  21, 23, 27, 29, 35 &   4 &  17 \\             
10,  &  15, 21, 23, 27, 29, 35 &   10 &  19 \\           
2, 4, 10,  &  17, 23, 27, 29, 35 &   2 &  19 \\          
2, 4, 10,  &  19, 21, 27, 29, 35 &   4 &  17 \\          
2, 4, 10,  &  19, 23, 27, 29, 35 &   4 &  17 \\          
2, 4, 10,  &  21, 23, 27, 29, 35 &   2 &  19 \\          
2, 10,  &  15, 21, 23, 27, 29, 35 &   2 &  19 \\         
2, 10,  &  17, 21, 23, 27, 29, 35 &   2 &  19 \\         
4, 10,  &  15, 21, 23, 27, 29, 35 &   4 &  17 \\         
4, 10,  &  19, 21, 23, 27, 29, 35 &   4 &  17 \\         
2, 4, 10,  &  15, 21, 23, 27, 29, 35 &   2 &  19 \\      
2, 4, 10,  &  17, 19, 23, 27, 29, 35 & {\bf b} &  21 \\      
2, 4, 10,  &  17, 21, 23, 27, 29, 35 &   2 &  19 \\      
2, 4, 10,  &  19, 21, 23, 27, 29, 35 &   4 &  17 \\      
2, 10,  &  15, 17, 21, 23, 27, 29, 35 &   2 &  19 \\     
4, 10,  &  15, 19, 21, 23, 27, 29, 35 &   4 &  17 \\     
2, 4, 10,  &  15, 17, 21, 23, 27, 29, 35 &   2 &  19 \\  
2, 4, 10,  &  15, 19, 21, 23, 27, 29, 35 &   4 &  17 \\  
2, 4, 10,  &  17, 19, 21, 23, 27, 29, 35 & {\bf c} &  11 \\  
2, 4, 10, & 15, 17, 19, 21, 23, 27, 29, 35 & {\bf c} &  11 \\ 
  \hline   
\end{tabular}  
$$
\caption{Non-admissible modules for $(6,8,v)$}
\label{Table6-8-25}
\end{table}
}

\smallskip
Using this table and the list of all $273+25(v-9)/2$,
modules $\Delta$, we checked Conjecture \ref{CONJSUP}
for $t=1$ (i.e. ignoring the dimensions) for
quite a few $v$. Here and in {\em all\,} calculations we performed 
it appeared sufficient to replace the admissibility of $\d$ by 
the admissibility of all $D_i$ (separately), a potentially 
weaker condition.
Generalizing \cite{Pi}, we  
verified here that all subvarieties $J^m_{\r}[\d]$
are diffeomorphic to proper $\mathbb{A}^{\!N}$, so 
we are in the situation of (\ref{conjaff}).  
\smallskip

Let us provide some details. Recall that 
the $D$\~flag of length $m$ originated at $D$ is by definition 
the sequence of $D$\~sets for $\Gamma$\~modules:
\begin{align}\label{d-mods}
\d\!=\!\bigl\{D_0=D, D_1\!=\!D\cup \{g_1\}, \ldots,
D_m\!=\!D\cup \{g_1,g_2,\ldots,g_m\}\bigr\},
\end{align}
where the inequalities $g_1<g_2<\cdots<g_m\in G\setminus D$
are imposed. Let us illustrate numerically
the importance of this very
ordering in our definition of $D$\~flags.
 
As it results from Proposition \ref{NESTED},
each set  $D\cup \{g_i\}$ corresponds to a certain
$\Gamma$\~module for any $1\le i\le m$. However, apart
from torus knots, these conditions (imposed together)
are generally significantly weaker than the conditions
we need. The $\Ga$\~modules $D_i$ from
(\ref{d-mods}) are not always admissible
(come from some $\m$) if all $D\cup \{g_i\}$ do.
  
We set $\ep(D)=1,0$ correspondingly 
for admissible and non-admissible $D$ and 
$\ep(\d)=\prod_{i=1}^m \ep(D_i)$, where 
$D_i=D_{\m_i}$. Then $\ep(\d)=1$ is equivalent to 
the admissibility of $\d$ (in the example under consideration).
We expect this implication to hold in general, but cannot
prove this at the moment. Thus $\ep(\d)=1$ 
implies that $\prod_{i=1}^m \ep(D\cup g_i)=1$, but 
the latter condition is not generally insufficient for the 
former.

Namely, if the admissibility of $D$\~flags were defined 
as $\prod_{i=1}^m \ep(D\cup\{g_i\})=1$ instead of
$\ep(\d)=1$, i.e. separately for each and every 
$D\cup\{g_i\}$, then there 
would  be $14$ extra (wrong) terms (with multiplicities) in 
(\ref{conjaff68}) below. This clearly demonstrates that our 
flags are generally more subtle than using ``marks" for
torus knots in \cite{Gor1} 
and other related works. 

The smallest example is as follows. 
Using Table \ref{Table6-8-25} with $v=9$,
the set of all $g_i$ for 
$D=[10, 17, 19, 23, 27, 29, 35]$ 
such that $D\cup\{g_i\}$ is admissible
is $\{2, 4, 11, 21\}$. 
However 
$D\cup\{2,4\}=[2, 4, 10, 17, 19, 23, 27, 29, 35]$ 
is a non-admissible $D$\~set (it is marked by $b$ in
the table).
\smallskip

Finally, 
we arrive at the following identity:
\begin{align}\label{conjaff68}
\sum_{m=0}^\infty\, \sum_{\d=\{D_0,\ldots,D_m\}}
 \ep(\d)\, q^{|D_0|+m\,} a^m=
\end{align}
\renewcommand{\baselinestretch}{0.5} 
{\small
\(
1+q+2 q^2+3 q^3+5 q^4+7 q^5+10 q^6+12 q^7+16 q^8+19 q^9+22 q^{10}
+24 q^{11}+25 q^{12}+24 q^{13}+22 q^{14}+17 q^{15}+11 q^{16}
+5 q^{17}+q^{18}
+a \bigl(q+2 q^2+4 q^3+7 q^4
+12 q^5+18 q^6+26 q^7+35 q^8+46 q^9+56 q^{10}+66 q^{11}+72 q^{12}
+74 q^{13}+70 q^{14}+59 q^{15}+41 q^{16}+21 q^{17}+5 q^{18}\bigr)
+a^2 \bigl(q^3+2 q^4+5 q^5+9 q^6+16 q^7+24 q^8+36 q^9+48 q^{10}
+62 q^{11}+74 q^{12}+82 q^{13}+83 q^{14}+76 q^{15}+58 q^{16}
+34 q^{17}+10 q^{18}\bigr)
+a^3 \bigl(q^6+2 q^7+5 q^8+9 q^9+15 q^{10}+22 q^{11}+31 q^{12}
+38 q^{13}+44 q^{14}+44 q^{15}+38 q^{16}+26 q^{17}+10 q^{18}\bigr)
+a^4 \bigl(q^{10}+2 q^{11}+4 q^{12}+6 q^{13}+9 q^{14}+11 q^{15}
+11 q^{16}+9 q^{17}+5 q^{18}\bigr)
+a^5 \bigl(q^{15}+q^{16}+q^{17}+q^{18}\bigr).
\)
}
\renewcommand{\baselinestretch}{1.2} 
\smallskip

The latter sum exactly coincides with 
$\h_{\,\vec\rr,\,\vec\ss}\,(\square;q,t=1,a)$
from formula $({\mathbf 4})$ in Section 3.2 of \cite{ChD1},
where 
$$
\vec\rr=\{4,2\},\, \vec\ss=\{3,1\},\, \t= C\!ab(25,2)T(4,3).
$$
Also see formula (3.4) there for the corresponding
DAHA-Betti polynomial, which is 
the $a$\~constant term of (\ref{conjaff68}) upon the
substitution $q\mapsto 1/t$ and multiplication by $t^{18}$:
\renewcommand{\baselinestretch}{0.5} 
{\small
\(
 1+5 t+11 t^2+17 t^3+22 t^4+24 t^5+25 t^6+24 t^7+22 t^8+19 t^9
+16 t^{10}+12 t^{11}
+10 t^{12}+7 t^{13}+5 t^{14}+3 t^{15}+2 t^{16}+t^{17}+t^{18}.
\)
}
\renewcommand{\baselinestretch}{1.2} 
\smallskip
Here we use the super-duality.

We also checked that $J^m_{\r}[\d]$ are always affine spaces,
found formulas for the dimensions of $J^{m=0,1}_{\r}[D]$ and 
correspondingly verified (\ref{conjaff}) for the coefficients 
of superpolynomials of $a^0,a^1$ in the considered case.

\comment{
\begin{table}[ht!]
\centering
\begin{tabular}{|l|l|}
  \hline 
  a & b \\
  c & d \\
  \hline   
\end{tabular}     
\hspace {0.1cm}
\begin{tabular}{|l|l|}
  \hline 
  a & b \\
  c & d \\
  \hline   
\end{tabular}  
\label{Tablex}
\caption{Caption}
\end{table}
}
\smallskip

\subsection{\bf Beyond ``Piontkowski"}
Examples are given in \cite{Pi} where his
approach does not work because the corresponding
cells $J_r[D]$ are non affine spaces. Thus the
count of such cells and knowing their dimensions
is insufficient for obtaining the Betti numbers
and the Euler number of $J_r$. Some 
`negative" examples are provided in the table
after Theorem 13 in \cite{Pi}. We found that
not always the cells in his table are really non-affine,
but the phenomenon he discovered is of course important.
It is unclear whether non-affine cells and non-admissible
$D,\d$ are topological.
This is of obvious interest due to our Main Conjecture
and because of the exciting link to orbital integral in
type $A$. 
\smallskip

We mostly consider a relatively simple singularity with  
Puiseux exponents $(6,9,10)$ and the
cable $C\!ab(19,3)T(3,2)$.
Its ring is $\r=\C[[t^6,t^9+t^{10}]]$, the valuation
semigroup is $\Gamma=\lan 6,9,19\ran$ and 
$\de=|\N\setminus \Ga|=21$ in this case.
See \cite{ChD1}, Section 3.2 for details and the formula
for the DAHA superpolynomial; the DAHA-Betti polynomial
($a=0,q=1$) is formula (3.3) there, which does provide
correct Betti numbers for the corresponding $J_{\r}$.
\smallskip

This example is transitional in a sense;
all cells are still affine spaces, but the justification
of this fact (straightforward, using computers)
becomes more involved. We considered some deformations
of parameters of $\r$ and think that Table \ref{Table6-9-19} 
depends only on the corresponding $\Ga$ (i.e. this table is 
of topological nature) but this is not clear. 

{\tiny
\begin{table}[ht!]
$$
\centering
\begin{tabular}{|ll|}
\hline    
\ \ \ \ {22, 41} \ \ \ \ \ \ \ \ \ \ \ \ \ \ \ \ \ \hfill
{2, 3, 5, 8, 11, 14, 16, 17, 20, 22, 23, 26, 29, 32, 35, 41}&\\
\ \ \ \ {3, 22, 41}\hfill
{1, 3, 4, 7, 10, 13, 16, 17, 20, 22, 23, 26, 29, 32, 35, 41}&\\
\ \ \ \ {22, 35, 41}\hfill
{3, 5, 8, 11, 14, 16, 17, 20, 22, 23, 26, 29, 32, 35, 41}&\\
\ \ \ \ {3, 22, 32, 41}\hfill
{2, 3, 8, 11, 14, 16, 17, 20, 22, 23, 26, 29, 32, 35, 41}&\\
\ \ \ \ {3, 22, 35, 41}\hfill
{2, 3, 5, 8, 11, 14, 17, 20, 22, 23, 26, 29, 32, 35, 41}&\\
\ \ \ \ {16, 22, 35, 41}\hfill
{1, 4, 7, 10, 13, 16, 17, 20, 22, 23, 26, 29, 32, 35, 41}&\\
\ \ \ \ {3, 16, 22, 35, 41}\hfill
{1, 3, 4, 7, 10, 13, 16, 20, 22, 23, 26, 29, 32, 35, 41}&\\
\ \ \ \ {3, 22, 29, 35, 41}\hfill
{3, 8, 11, 14, 16, 17, 20, 22, 23, 26, 29, 32, 35, 41}&\\
\ \ \ \ {3, 22, 32, 35, 41}\hfill
{3, 5, 11, 14, 16, 17, 20, 22, 23, 26, 29, 32, 35, 41}&\\
\ \ \ \ {16, 22, 32, 35, 41}\hfill
{3, 5, 8, 11, 14, 17, 20, 22, 23, 26, 29, 32, 35, 41}&\\
\ \ \ \ {3, 16, 22, 29, 35, 41}\hfill
{3, 4, 7, 10, 13, 16, 20, 22, 23, 26, 29, 32, 35, 41}&\\
\ \ \ \ {3, 16, 22, 32, 35, 41}\hfill
{2, 3, 8, 11, 14, 17, 20, 22, 23, 26, 29, 32, 35, 41}&\\
\ \ \ \ {3, 22, 26, 32, 35, 41}\hfill
{1, 4, 7, 10, 13, 16, 20, 22, 23, 26, 29, 32, 35, 41}&\\
\ \ \ \ {3, 22, 29, 32, 35, 41}\hfill
{4, 7, 10, 13, 16, 20, 22, 23, 26, 29, 32, 35, 41}&\\
\ \ \ \ {13, 16, 22, 32, 35, 41}\hfill
{3, 11, 14, 16, 17, 20, 22, 23, 26, 29, 32, 35, 41}&\\
\ \ \ \ {3, 13, 16, 22, 32, 35, 41}\hfill
{3, 8, 14, 16, 17, 20, 22, 23, 26, 29, 32, 35, 41}&\\
\ \ \ \ {3, 16, 22, 26, 32, 35, 41}\hfill
{3, 8, 11, 14, 17, 20, 22, 23, 26, 29, 32, 35, 41}&\\
\ \ \ \ {3, 16, 22, 29, 32, 35, 41}\hfill
{3, 5, 11, 14, 17, 20, 22, 23, 26, 29, 32, 35, 41}&\\
\ \ \ \ {3, 22, 23, 29, 32, 35, 41}\hfill
{3, 4, 7, 10, 13, 16, 22, 23, 26, 29, 32, 35, 41}&\\
\ \ \ \ {3, 22, 26, 29, 32, 35, 41}\hfill
{4, 7, 10, 13, 16, 22, 23, 26, 29, 32, 35, 41}&\\
\ \ \ \ {13, 16, 22, 29, 32, 35, 41}\hfill
{3, 14, 16, 17, 20, 22, 23, 26, 29, 32, 35, 41}&\\
\ \ \ \ {3, 13, 16, 22, 29, 32, 35, 41}\hfill
{3, 11, 16, 17, 20, 22, 23, 26, 29, 32, 35, 41}&\\
\ \ \ \ {3, 16, 22, 23, 29, 32, 35, 41}\hfill
{3, 11, 14, 17, 20, 22, 23, 26, 29, 32, 35, 41}&\\
\ \ \ \ {3, 16, 22, 26, 29, 32, 35, 41}\hfill
{3, 8, 14, 17, 20, 22, 23, 26, 29, 32, 35, 41}&\\
\ \ \ \ {3, 20, 22, 26, 29, 32, 35, 41}\hfill
{3, 7, 10, 13, 16, 22, 23, 26, 29, 32, 35, 41}&\\
\ \ \ \ {3, 22, 23, 26, 29, 32, 35, 41}\hfill
{7, 10, 13, 16, 22, 23, 26, 29, 32, 35, 41}&\\
\ \ \ \ {10, 13, 16, 22, 29, 32, 35, 41}\hfill
{3, 16, 17, 20, 22, 23, 26, 29, 32, 35, 41}&\\
\ \ \ \ {3, 10, 13, 16, 22, 29, 32, 35, 41}\hfill
{3, 14, 17, 20, 22, 23, 26, 29, 32, 35, 41}&\\
\ \ \ \ {3, 16, 20, 22, 26, 29, 32, 35, 41}\hfill
{3, 14, 16, 20, 22, 23, 26, 29, 32, 35, 41}&\\
\ \ \ \ {3, 16, 22, 23, 26, 29, 32, 35, 41}\hfill
{3, 11, 17, 20, 22, 23, 26, 29, 32, 35, 41}&\\
\ \ \ \ {3, 17, 22, 23, 26, 29, 32, 35, 41}\hfill
{3, 7, 10, 13, 16, 22, 26, 29, 32, 35, 41}&\\
\ \ \ \ {3, 20, 22, 23, 26, 29, 32, 35, 41}\hfill
{7, 10, 13, 16, 22, 26, 29, 32, 35, 41}&\\
\ \ \ \ {10, 13, 16, 22, 26, 29, 32, 35, 41}\hfill
{3, 17, 20, 22, 23, 26, 29, 32, 35, 41}&\\
\ \ \ \ {3, 10, 13, 16, 22, 26, 29, 32, 35, 41}\hfill
{3, 16, 20, 22, 23, 26, 29, 32, 35, 41}&\\
\ \ \ \ {3, 14, 20, 22, 23, 26, 29, 32, 35, 41}\hfill
{3, 16, 17, 22, 23, 26, 29, 32, 35, 41}&\\
\hline   
\end{tabular}     
\hspace {0.1cm}
$$
\caption{Non-admissible $D$ for $\Ga=\lan 6,9,19\ran$}
\label{Table6-9-19}
\end{table}
}

 
The simplest example of non-affine $J_{\r}[D]$ 
in this case given in \cite{Pi} is
$$D=[3, 7, 10, 13, 16, 17, 20, 22, 23, 26, 29, 
32, 35, 41],$$ 
but we found that the corresponding $J_{\r}[D]$
is biregular to $\mathbb{A}^{14}$. The only problem
with this and $2$ other similar sets $D$ (our program
obtained) is that the natural $\la$\~variables 
Piontkowski uses (we too) are 
inconvenient to parametrize 
$J_{\r}[D]$; a certain (linear) change of variables is necessary.
The other two $D$\~sets with a similar behavior (when a straight
elimination of the $\la$-variables is insufficient) are 
\begin{align*}
&[3,10,13,16,20,22,23,26,29,32,35,41], \\ 
&[3,10,13,16,17,20,22,23,26,29,32,35,41].
\end{align*}
All $3$ (and any other cells) are affine spaces 
in this case.  The next example of a non-affine 
cell from \cite{Pi} is for $\C[[t^6,t^9+t^{13}]]$; we
confirm it. See Appendix A in the online version of
this paper.

\smallskip
We note that our program routinely calculates 
$|J_{\r}[D](\mathbb{F}_3)|$ for ``suspicious"
$D$ to double check the direct verification of the 
affineness of the cells (mostly automated). These 
cardinalities must be $3^{dim}$ for affine cells. 
The prime $3$ is a ``place  of good reduction" for this 
$\r$. The reduction is bad modulo $p=2$  
and one needs to switch to the {\em topologically\,}
equivalent ring $\C[[t^6+t^7,t^9]]$ before replacing
$\C\mapsto \mathbb{F}_2$. We omit general analysis
of places of bad reduction in this paper.

We will not discuss much the flags here, but let us mention
that all cells are affine for $J^{m=1}_{\r}[\d]$ 
(their dimensions are all calculated)
in the case of $\r=\C[[t^6,t^9+t^{10}]]$.  
The list of (all) $70$ non-admissible $D$ (i.e. for
$m=0$) will be provided. The total number
of modules $\De$ is $447=377+70$ in this case, 
and the Euler number is 
$e(J_{\r})=377$.  Accordingly,
we checked (numerically) the coincidence 
from (\ref{conjaff}) in the following two cases:
$(a)$ for all $a$ when $t=1$ (for the admissibility of $\d$
understood as the admissibility of all $D_i$ in this flag),
and $(b)$ for the coefficients of
$\h(\square;q,t,a)$  from \cite{ChD1} of $a^{0,1}$.
\smallskip

The calculation in $(a)$ greatly demonstrates the role
of admissibility and the implications of
the ordering $g_1\!<\ldots<\!g_m$, which are quite
non-trivial combinatorially. The corresponding
reduction of the superpolynomial
for $\r=\C[[t^6,t^9+t^{10}]]\,$ is $\ \h(\square;q,t=1,a)=$
\vskip 0.2cm

\renewcommand{\baselinestretch}{0.5} 
{\small
\(
1+q+2 q^2+3 q^3+5 q^4+7 q^5+10 q^6+13 q^7+17 q^8+21 q^9+25 q^{10}
+29 q^{11}+33 q^{12}+36 q^{13}+37 q^{14}+37 q^{15}+34 q^{16}
+28 q^{17}+20 q^{18}+12 q^{19}+5 q^{20}+q^{21}+a^5 \bigl(q^{15}
+q^{16}+2 q^{17}+2 q^{18}+3 q^{19}+2 q^{20}+q^{21}\bigr)
+a^4 \bigl(q^{10}+2 q^{11}+4 q^{12}+7 q^{13}+11 q^{14}+15 q^{15}
+19 q^{16}+22 q^{17}+23 q^{18}+21 q^{19}+13 q^{20}+5 q^{21}\bigr)
+a \bigl(q+2 q^2+4 q^3+7 q^4+12 q^5+18 q^6+27 q^7+37 q^8+50 q^9
+63 q^{10}+78 q^{11}+91 q^{12}+105 q^{13}+113 q^{14}+118 q^{15}
+114 q^{16}+100 q^{17}+76 q^{18}+48 q^{19}+22 q^{20}
+5 q^{21}\bigr)+a^3 \bigl(q^6+2 q^7+5 q^8+9 q^9+16 q^{10}
+24 q^{11}+36 q^{12}+47 q^{13}+61 q^{14}+71 q^{15}+81 q^{16}
+82 q^{17}+76 q^{18}+57 q^{19}+32 q^{20}+10 q^{21}\bigr)
+a^2 \bigl(q^3+2 q^4+5 q^5+9 q^6+16 q^7+25 q^8+38 q^9+53 q^{10}
+71 q^{11}+90 q^{12}+109 q^{13}+126 q^{14}+138 q^{15}
+143 q^{16}+134 q^{17}+111 q^{18}+75 q^{19}+38 q^{20}
+10 q^{21}\bigr).
\)
}
\renewcommand{\baselinestretch}{1.2} 
\smallskip

\setcounter{equation}{0}
\section{\sc Some perspectives}\label{sec:CONCL}
The topics below are mostly open projects, but we believe
that this section can be of interest; the
relation to orbital integrals, affine Springer fibers and
the motivic reformulation (\ref{conjmotxx}) of
(\ref{conjmot}) are the key. We will first address
the absence of colors and links in our Main Conjecture and related
issues.

\smallskip

\subsection{\bf Adding colors}
In contrast to \cite{ChD2} and previous paper, we
restrict this one to algebraic uncolored knots and
only type $A$ is addressed in our conjecture. Let
us briefly comment on this. Adding colors is
expected via the curves in \cite{Ma}, so algebraic 
{\em links\,} are/seem necessary for this. Apart
from  rectangle Young diagrams, the coefficients
of DAHA-superpolynomials are non-positive, which is
an obvious challenge for the geometric interpretation.
This is the same for links (included uncolored ones). 
The Jacobian factors become {\em ind-schemes} for algebraic 
links and we need to follow \cite{KL} (divide by certain 
split tori) to make them {\em proper\,}. 
 
One can try to bypass the non-positivity issues switching to 
some {\em Hilbert schemes} instead of our 
flagged Jacobians, which corresponds to  
the unreduced topological setting. This can be similar to 
\cite{ORS,EGL}; see also Conjecture 5.3 $(ii)$ from \cite{ChD2}.
\smallskip

{\sf Superpolynomials beyond $A_n$\,}.
Our flagged Jacobian factors are related
to the Hitchin and affine Springer fibers. The hope is 
that the DAHA-superpolynomials can be directly and 
geometrically determined by spectral curves for
types $B_n,C_n,D_n$ via the corresponding
flagged Jacobian factors and/or Hilbert schemes. The
rank stabilization was conjectured in \cite{CJ};
the name ``hyperpolynomials" is used instead of
superpolynomials for non-$A$ types. The spectral 
curves are generally non-unibranch, so the passage 
to links is necessary here. A certain confirmation 
is the following observation. 

It is true in all known examples that  
$\h(b\,;\,q,t,u,a\!=\!0)$ for the
hyperpolynomials from Section 4.2 of \cite{CJ}
depend only on the corresponding knot and the color $b$.
This holds for the hyperpolynomials of type $C,D$
calculated there for $T(3,2),T(7,2),T(4,3)$ and 
$b=\om_1,2\om_1,\om_2$. The hyperpolynomials of type $D$ 
are the specialization of those of type $C$ upon $u=1$ 
(see (4.7) there; $u$ is the second $t$ for $B,C$).
The same holds for type $B$ under 
$q^2\mapsto q, u\mapsto t$. 

Furthermore,
the hyperpolynomials $\h_{\rr,\ss}^{\mathfrak{ad}}(q,t,a)$
for the {\em Deligne-Gross series\,} (extending the
series $E_{6,7,8}$) introduced in  
Section 4.2 of \cite{ChEl} coincide at $a=0$ with 
$\h_{\rr,\ss}(\square,\,;\,q,t,a=0)$ of type $A$.
The existence of $\h_{\rr,\ss}^{\mathfrak{ad}}$
was confirmed only partially in \cite{ChEl}
(and only for $T(3,2),T(4,3)$). These 
$\h_{\rr,\ss}^{\mathfrak{ad}}$ are mysterious
from the viewpoint of \cite{CJ}, since there can be no
rank stabilization here. We hope that our present paper can shed
some light on this. 

\subsection{\bf Perfect DAHA modules}
A challenge is an interpretation of DAHA superpolynomials 
for arbitrary algebraic knots, similar to that in \cite{Gor2}
for torus ones via perfect modules of rational DAHA. See also
Theorem 9.5 from \cite{GORS}. There are some obstacles. 

One can use here the classification results from
\cite{Vas,VV,Yun} in terms of $K$\~theory or
homology (for the rational and trigonometric DAHA) of {\em 
iwahoric Springer fibers}. See paper \cite{Vas} for general
theory, including a comprehensive analysis in type $A$.

Following Gorsky, let us consider the root system
$A_{\ss-1}$ for a given torus knot $T(\rr,\ss)$
(so we adjust its rank to the knot)
and take $e\!=\!e_1\!+\!\ldots\!+\!e_{\ss-1}\!+
\!z^{\rr}f_\theta$
in the notations from Section 4.1 from \cite{VV}. Then
$H_*(\y_e)$ for the corresponding Iwahoric Springer fiber
$\y_e$ can be supplied with a natural structure
of the {\em perfect\,} $\HH$\~module for 
$t\!=\!q^{-\frac{\rr}{\ss}}$ 
of dimension $\rr^{\ss-1}$, which is a very explicit
quotient of the polynomial representation of
rational DAHA. This  module becomes simpler 
with $q,t$ (via $K$\~theory), but the grading then will be
missing.

The coefficients of $a^m$ in the superpolynomial are
associated with isotypic components of this $\HH$\~module
for the wedge $m{\hbox{\tiny th}}$ powers of the standard 
$\ss-1$ dimensional representation of $W=\S_{\ss}$.
This construction is quite combinatorial; Gorsky relates $a^m$ 
to $m$\~sets of {\em marks\,}, which are corners in Dyck paths 
below the diagonal in the $\ss\times\rr$\~rectangle. 
In the absence of marks (when $a=0$), the  
space of $W$\~coinvariants is the image of the projection of 
$H_*(\y_e)$ onto  $H_*(\x_e)$ for the affine Springer fiber 
$\x_e$ (see also Section \ref{sec:ORBITAL} below). 
\smallskip

The number of marks corresponds to $m$ in our $J_{\r}^m$.
The latter is a special subvariety of the 
{\em parahoric Springer fiber} of full $m$\~sub-flags 
starting from the top, defined as the
space of flags of $\r$\~modules $\{\m_i, i=0,\ldots,m\}$, where 
$\mathfrak{m}\m_m\subset \m_i$ and dim$_{\C}(\m_{i+1}/\m_i)=1$;
$\mathfrak{m}\subset \r$ is the maximal ideal. They can be
extended to full {\em periodic} flags via the 
forgetful map from $\y_e$ to the parahoric one. 

This construction involves the root system $A_{\ss-1}$, and
the topological $\rr\! \leftrightarrow\! \ss$ symmetry  of
the superpolynomial generally becomes far from obvious.
Our geometric superpolynomials are defined directly in
terms of the singularity and are manifestly
$\rr\! \leftrightarrow\! \ss$\~symmetric.

Furthermore, the (finite-dimensional) spaces $H_*(\y_e)$ for 
arbitrary nil-elliptic (anisotropic) $e$ 
are not generally related to any DAHA\~modules
unless in the torus case.
These spaces are needed for general unibranch
plane curve singularities, but  
their dimensions and other features are 
very different from those of the finite-dimensional modules.
It is not impossible that they can appear as some 
``remarkable" subspaces in infinite-dimensional $\HH$\~modules, 
but this is questionable. 

\subsection{\bf Affine Springer fibers}\label{sec:ORBITAL}
Conjecture \ref{CONJSUP} leads to the following 
connection with the orbital integrals from the 
Fundamental Lemma (in the geometric setting). We will
mainly follow Section 2.4 from \cite{Yu}.  
The example from Section 2.4.2 there is for the spectral
curve $y^{\rr}=x^{\ss}$ \cite{LS}, which can be actually 
generalized to any (germs of) plane curve singularities
$\c$ such that the Jacobian factor 
$J_{\r}$ at $(0,0)$ has the Piontkowski-type decomposition with 
affine cells only. 

The following discussion will be mostly restricted to type $A$
in the anisotropic case. 
Note that fractional ideals of $\r$
of {\em degree zero\,} are considered in \cite{Yu} (and
in the Fundamental Lemma) in the definition of the 
compactified Jacobians. 

Let $1/t=p^{\ell}$ be the cardinality of a finite 
field $\mathbb{F}$.
The choice $t=1/|\mathbb{F}|$ is standard when connecting the
$q,t$\~theory of spherical functions
with the $p$\~adic theory (see e.g.
\cite{C101}), though using $q$ instead of $1/t$
is equally fine due to the
super-duality of the DAHA superpolynomials.
 
Because we stick to the unibranch case, 
$\ga\in \mathfrak{g}(\mathbb{F}((x)))$ is assumed 
{\em nil-elliptic}: no
split tori
over a local field $\mathbb{F}((x))$ in
the stabilizer $G_{\ga}$ ($\ga=e$ was used above). 
The affine Springer fiber $\x_{\ga}$ is then formed by the
classes of $g$ in the {\em affine Grassmannian} 
$G(\mathbb{F}((x)))/G(\mathbb{F}[[x)]])$ over the field
$\mathbb{F}$ such that
$Ad_g^{-1}(\ga)\in \mathfrak{g}(\mathbb{F}[[x]])$. 
Here Lie$(G)=\mathfrak{g}$; see e.g. \cite{Yu}.
Yun denotes the corresponding {\em parahoric,} ones by 
$\x_{P,\ga}$, which contain our flagged Jacobian factors
(though we do not need any $\mathfrak{g},G$). The
Iwahoric Springer fibers are for $P=I$ in his notation.

\vskip 0.2cm
The affine Springer fiber $\x_{\ga}$ can be naturally
identified with the compactified Jacobians for 
rational curves with the local ring $\r$ at its (unique)
singularity; see e.g. \cite{La,ChU}. A general
construction is from 
\cite{ChU}. Let $\mathscr{G}$ be a {\em factorizable} Lie group 
schemes over a smooth projective curve $E$, 
defined by the conditions $H^1(E,Lie(\mathscr{G}))=\{0\}=
H^0(E,Lie(\mathscr{G}))$ for the corresponding sheaf
of Lie algebras $Lie(\mathscr{G})\to E$. For
a scheme subtorus  $\mathscr{T}\subset \mathscr{G}$ such that its
generic fiber is a maximal torus, embeddings 
$f:\mathscr{T}\hookrightarrow \mathscr{G}$ become
conjugations over sufficiently general open $U\subset E$:
$f(\xi)=\phi^{-1} \xi\phi$ for $\xi\in H^0(U,\mathscr{T})$
and meromorphic sections $\phi=\phi_U$ of $\mathscr{G}$
over $U$. $\check{\hbox{C}}$ech cohomology is used.

Assuming that $\phi_i$ exist for open $U_i$ such that 
$E=\cup_i U_i$, 
the map $\{f\}\ni f\mapsto \{\phi_i\phi_j^{-1}\in 
H^0(U_i\cap U_j,\mathscr{G})\}\to H^1(E,\mathscr{T})$
is an isomorphism. Actually, we use this map in the opposite
direction, from $H^1(E,\mathscr{T})$ to $\{f\}$.
Here $H^1(E,\mathscr{T})$ is the generalized 
Jacobian of the cover $F\to E$ if 
$\mathscr{T}=\o_F^*\subset \mathscr{G}$.
Note that $E$ can be only $\mathbb{CP}^1$ or an
elliptic curve; 
generalized Prim varieties appear for the latter.
\smallskip

Let $F$ be rational with exactly one singular point that is
the whole fiber $F_o$ over  some $o\!\in\! E\!=\!\mathbb{CP}^1$;
here $\mathscr{G}_o$ must be $G$. Then $f(\mathscr{T})$ can be 
obtained as $g^{-1}\mathscr{T} g$ for 
{\em rational sections\,} $g$ of $\mathscr{G}$ such that 
$g$ are regular at $E\setminus o$ and 
$g^{-1}\mathscr{T} g$ are regular at $F_o$.
Thus  such $g$ form an open subset in the affine 
Springer fiber $\mathscr{X}_\ga$ at $o$ for a sufficiently general
section $\ga$ of $\hbox{Lie}(\mathscr{T})$ regular
at $o$. For any fiber $F_o$, the map $g\mapsto 
g^{-1}\mathscr{T} g$ goes through the quotient
$L_\ga\!\setminus\! \mathscr{X}_\ga$ for the group $L_\ga$
of rational sections of $F$ regular at $F\setminus F_o$; cf.  
\cite{Yu}, Theorem 2.9 and \cite{KL}.



\vskip 0.1cm
Type $A$ is not necessary here, but we need this in
our paper (and anisotropic $G_{\ga}$). Recall that our 
{\em flagged\,} Jacobian factors  
deviate  from the usual ones. First, we consider only
standard $\r$\~modules. Second,
admissible $D$\~flags $\d\!=\![\,D_0\!\subset\! 
D_1\!=\!D_0\!\cup\{g_1\} 
\ldots\subset D_m\,]$ are subtle; 
$D_i$ must be $D$-sets and the ordering
$g_1< g_2<\ldots <g_m$ is imposed. This ordering and the
admissibility are quite non-trivial geometrically.
\subsection{\bf Motivic approach}\label{sec:MOTIV}
Our construction results in the following generalization
of orbital integrals: 
$t^{\de}\sum_{\d} q^{m+|D_0|}\,a^m
\mid\!J^m_{\r}[\d](\mathbb{F}_{1/t})\!\mid$ in type $A$
(the nil-elliptic case), where the 
summation is over all admissible $D$\~flags $\d$. Their
interpretation as ``natural" orbital integrals is one of
the main challenges triggered by this work. 
It seems doable due to an entirely geometric nature 
of our approach. Let us also mention here potential adding colors
to our construction, another challenge for us and the
specialists in orbital integrals. 

For sufficiently general prime $p$, there is solid 
evidence that such sums  coincide with the DAHA superpolynomials 
$\h_{\c}(\square;q,t,a)$ associated with the singularity $\c$,
which is stated in (\ref{conjmot}). For $q\!=\!1,a\!=\!0$, we 
arrive at {\em orbital integrals\,} $O_{\ga}$. They are 
expected to be $|\mathbb{F}|$\~integral and positive,
which matches our Main Conjecture. 
\smallskip

It is not
impossible that all classes $[J^m_{\r}[d]]$
in $K_0(Sch/\C)$
are sums of classes $[\mathbb{A}^1]^N$ over $\Z_+$, the
strongest possible (motivic) assumption. 
This would match known and conjectured properties of 
$\h_{\c}$. It is possible that  $J^m_{\r}[d]$ are 
always paved by affine spaces (even when $J^m_{\r}[\d]$ are
not all affine); see 
Appendix A in the online version of this paper.

If (\ref{conjmot}) holds,  then DAHA superpolynomials
provide all virtual Hodge numbers of $J_{\r}^{m}[d]$.
This results from Part $(3)$ of Theorem 1 from the Appendix
by N.Katz in \cite{HR}. Following it, let
$E(X; x, y) = \sum_{r,s}e_{r,s}x^ry^s$ for a separated
scheme $X/\C$ of finite type and {\em virtual Hodge numbers\,}
$
e_{r,s} = \sum_i (-1)^i h^{r,s}(gr_W^{r+s}(H_c^i(X^{an},\C))).
$

Then (\ref{conjmot}) gives that $X=J_{\r}^{m}[d]$
is {\em strongly polynomial-count\,} and the following formula 
holds:
\begin{align}\label{conjmotx}
&(xy)^\delta\h_{\,\vec\rr,\,\vec\ss}\,(\square;q,1/(xy),a)=
\sum_{d,m} E(J^m_{\r}[d];x,y) q^{d+m} a^m\notag\\
&=
\sum_{d,m,i,r} (-1)^i\,  q^{d+m} a^m (xy)^r\,\, \hbox{rk} 
(gr_W^{2r}(H_c^i(J^m_{\r}[d],\C)),
\end{align}
which directly links the DAHA superpolynomials to
the weight filtration in $H_c^i(X^{an},\C)$. We use that
$\r$ can be assumed over $\Z$ and then $J^m_{\r}[d]$
are defined over a localization of $\Z$ by finitely many
primes.
Thus $e_{r,s}\!=\!0$ for $r\!\neq\! s$ (i.e. this is a Tate Hodge
structure) and the right-hand side in (\ref{conjmot})
can be replaced by the following $\h^{wt}(q,t,a)$:
\begin{align}\label{conjmotxx}
&\h_{\,\vec\rr,\,\vec\ss}\,(\square;q,t,a)\!=\!
\sum_{d,m,i,r} (-1)^i\,  q^{d+m}\, t^{\de- r}\,
a^m\,  \hbox{rk} 
(gr_W^{2r}(H_c^i(J^m_{\r}[d],\C)).
\end{align}

\smallskip

{\sf Conclusion.}
In spite of some similarity of (\ref{conjmotxx}) and
Theorem 9.5 \cite{GORS}, we 
do not see how they can be connected. First,
$t_{st}^2=q/t$ is used there for the weight 
filtration instead of $t$ in (\ref{conjmotxx}).
Second, our construction
does not require affine Springer
fibers and picking the  corresponding {\em isotypic\,}
components in their homology. Lie groups do not 
appear in our approach; in a sense, this corresponds to
the ``endoscopy part" of Fundamental Lemma.
Third, related to the second, our admissible flags are
new and different from those in \cite{GORS}. 
Fourth, our approach is fully computational and 
we calculate well beyond torus knots and the series
in \cite{Pi}; the examples provided in \cite{ORS}
were only for some simple torus knots. In spite of these 
differences, the ORS conjecture can be still compatible with our
one,  but this is not clear. 
\smallskip

Conceptually, {\em nested Hilbert
schemes} from \cite{ORS} and similar objects in 
related geometry-topology and
physics result in some {\em infinite} Poincar\'e 
series. Generally the ``ultimate" problem is to transform 
them into polynomials or some {\em finite} expressions, 
such that
the resulting coefficients are positive integers. The latter
positivity is generally much more subtle than the positivity
(if any) of the coefficients of the initial series.
In this paper, we provide the conjectural 
geometric interpretation for the DAHA superpolynomials
via the {\em Jacobian factors} instead of Hilbert schemes,
which is therefore ultimate in the sense above. There
is recent progress with the geometric interpretation
of the positivity conjecture for DAHA superpolynomials
colored by symmetric or wedge powers, but so far
not with arbitrary rectangles.

\smallskip

We note that not many formulas are known for the stable
$KhR$\~polynomials (the celebrated Khovanov polynomials for 
$sl_2$ are exceptional). They are mostly for
$T(2\mm +1,2)$ and no formulas are known for iterated
non-torus knots. Thus checking our geometric 
superpolynomials versus the DAHA ones is actually
the only way for iterated torus knots; though see\cite{ChD1} 
for some conditional verifications of our topological
conjecture using the reduction to 
the Khovanov polynomials.
\smallskip

We do not pay any special attention
to torus knots in our work. For such knots, we
generalize the approach from \cite{GM1},\cite{GM2}, 
where  $a=0$. Our usage of {\em standard modules}
is the key; this is not fully understood geometrically,
but we already have some motivic interpretation of
complete geometric superpolynomials (with all $3$
parameters). We note that even for torus knots, our 
{\em flagged Jacobian factors\,} are new.
The conditions on the 
corresponding $\Gamma$\~modules $\De$ in their definition
are not clear by now from the viewpoint of usual flagged 
constructions.


We plan to approach arbitrary 
{\em colored\,} algebraic {\em links\,} and possibly
reach arbitrary root systems our further
papers. It is a must if we want to realize the
potential of DAHA in full and for connections with 
$p$\~adic orbital integrals. Also, there is a realistic 
program of justifying
our conjecture (related to \cite{Ma}). It requires 
knowing the behavior of our geometric
superpolynomials under monoidal-type transformations of 
the corresponding singularities. All of them must be
considered here, not only those for torus knots.

\smallskip
{\sf Acknowledgements.} We thank 
Evgeny Gorsky, Anatoly Libgober, Davesh Maulik and
Andras Szenes for useful 
discussions. The referee's suggestions were very helpful.
Our special thanks to David Kazhdan. The first 
author thanks University of Geneva for the invitation and
hospitality. 

\comment{
Our special thanks go to Vivek Shende for clarifying discussions
on the topics related to \cite{ORS,GSh} and Mikhail
Khovanov for many conversations on his and Rozansky's theory. 
We thank Semen Artamonov for our using his unique software for 
calculating colored HOMFLYPT polynomials, Peter Samuelson for 
sending us his work before it was posted and him and Yuri Berest 
for discussions. The first author thanks Andras Szenes for the 
talks on the algebraic-geometric aspects of our construction and 
the 
University of Geneva for the invitation and hospitality. We
thank the referees for their attention to our paper,
thorough reports and important suggestions.

The paper was mainly written at RIMS. The first author thanks 
Hiraku Nakajima and RIMS for the invitation and hospitality and
the Simons Foundation;
the second author is grateful for the invitation to the 
school \& conference "Geometric Representation Theory" 
at RIMS (July 21- August 1, 2014) and generous RIMS' support.
I.D. acknowledges partial support from.
}

\vskip -4cm
\bibliographystyle{unsrt}

\vfill\eject

\appendix

\centerline{\bf Two appendices}
\vskip 0.5cm 

In these two appendices we will use the notation
and formula numbers from the main body of the
paper. The first one is an important part of the
paper; it provides examples of non-affine Piontkowski
cells, i.e. those non-isomorphic to affine spaces
$\mathbb{A}^N$. Such cells influence  the geometric 
superpolynomial in quite a non-trivial way.
They are differences of some affine cells, so 
their corresponding contributions contain the
$a,q,t$\~monomials with negative coefficients. 
When $t=1$, some such cells contribute $0$, i.e. 
behave as non-admissible ones. 
In all considered examples, we calculated the geometric 
superpolynomials for $t=1$ for all $m$.

The second appendix is a complete list of all
dimensions in the case  of $\r=\C[[z^6,z^9+z^{13}]]$
corresponding to $\Ga\!=\!\lan 6,9,22\ran$.
It is supposed to be used together with the table
of non-affine cells in this case, which we provide.
There are many interesting properties and symmetries
of these dimensions, which we do not systematically
discuss or do not touch at all. We think that this
table somewhat compensates this.

\setcounter{equation}{0}
\section{\sc Non-affine cells}
\label{app:noaff}
\subsection{\bf Puiseux exponents (6,9,13)}
Let us discuss the case of
$\r=\C[[z^6,z^9+z^{13}]]$, the third in the table
after Theorem 13 in \cite{Pi} of the cases beyond the
technique of this paper (including counting the Euler numbers). 
Here $\Ga=\lan 6,9,22\ran, \de=24$, the link is
$C\!ab(22,3)T(3,2)$.

Let $D^\dag$ be the set of entries
in a $D$\~set, which are {\em primitive\,}, i.e. 
the minimal set of generators of the module $\De$ over
$\Ga$; we omit $0$. This is obviously sufficient to recover
$D$. We confirm the claim from \cite{Pi} that
$J_{\r}[D]$ is non-affine for   
$D^\dag\!=\!\{3, 7, 10, 17, 20\}$, i.e. for 
$$
D\!=\![3, 7, 10, 13, 16, 17, 19, 20, 23, 25, 26, 29, 32, 35, 
38, 41, 47].
$$

Table \ref{Table6-9-22} gives {\em all\,} non-affine cells;
we provide the corresponding $D^\dag$,  
$d=|D|$, the (biregular) type of $J_{\r}[D]$, and
the corresponding contribution to the 
(geometric) superpolynomials.
We put $\mathbb{A}^N\vee\mathbb{A}^{N}$ for the amalgam (union)
of two $\mathbb{A}^N$ with the intersection $\mathbb{A}^{N-1}$.

{\small
\begin{table}[ht!]
\centering
\begin{tabular}{|l|c|c|l|}
\hline 
$D^{\dag}$-sets & $|D|$ & types & 
$q,t$-terms \\
\hline    
\{3, 11, 14\} &  14 & 
$\mathbb{A}^{17}\setminus \mathbb{A}^{16}$ & $q^{14}(t^7-t^8)$\\
\{3, 10, 13, 20, 23\} & 15 &
$\mathbb{A}^{16}\vee\mathbb{A}^{16}$ & $q^{15}(2 t^8-t^9)$\\
\{3, 11, 14, 19\} &  15 & 
$\mathbb{A}^{15}\setminus\mathbb{A}^{14}$ & $q^{15}(t^9-t^{10})$\\
\{3, 8, 11, 19\} &  16 & 
$\mathbb{A}^{15}\setminus\mathbb{A}^{14}$ & $q^{16}(t^9-t^{10})$\\
\{3, 10, 13, 17, 20\} & 16 &
$\mathbb{A}^{14}\vee\mathbb{A}^{14}$ & 
$q^{16}(2 t^{10}-t^{11})$\\
\{3, 7, 10, 17, 20\} & 17 &
$\mathbb{A}^{14}\vee\mathbb{A}^{14}$ & 
$q^{17}(2 t^{10}-t^{11})$\\
\hline   
\end{tabular}     
\hspace {0.1cm}
\vspace {0.2cm}
\caption{Non-affine $J_{\r}[D]$ for $\Ga=\lan 6,9,22\ran$}
\label{Table6-9-22}
\end{table}
}

The total number of
modules $\De$ is $605$, $e(J_{\r})=523$ and the total
number of non-admissible modules is $79$. The latter is
{\em not\,} now $605-523=82$, since $3$ modules from
Table \ref{Table6-9-22} of type 
$(\mathbb{A}^{N}\setminus \mathbb{A}^{N-1})$ do not
contribute to the Euler number. The other $3$ contribute
$2-1=1$, as with ordinary affine cells. Actually,
this is not the simplest example;
$\r=\C[[z^6,z^9+z^{11}]]$ has $3$ $D$ and $10$ 
$\d\!=\![D_0,\!D_1]$ of type 
$(\mathbb{A}^N\!\setminus\mathbb{A}^{N\!-\!1}).$

We have a sketch of the justification that 
$J_{\r}[d]$ for $d$ from the table are still paved by affine
spaces due to the ``redistribution" of cells. Namely, 
each subtracted $\mathbb{A}^{N}$ can be compensated
by adding an affine cell $J_{\r}[D']$
to $J_{\r}[D]$ from the table, where $|D'|=|D|$
and $J_{\r}[D']$ is {\em from the
boundary\,} of $J_{\r}[D]$. For instance, 
$
(\mathbb{A}^{N}\vee\mathbb{A}^{N})\cup 
\mathbb{A}^{\!N\!-\!1}
\ =\ \mathbb{A}^{N}\cup \mathbb{A}^{N},
$
where the cell $\mathbb{A}^{\!N\!-\!1}$ 
(of the same $d$) is taken from the boundary of the cell
$(\mathbb{A}^{N}\vee\mathbb{A}^{N})$.
Proposition \ref{DEGJAC} is used here. For instance,
our analysis gives that 2 cells with
$D_1^\dag\!=\!\{4, 7, 10, 14, 20\},\, 
 D_2^\dag\!=\!\{4, 7, 10, 17, 20\}$ belong to the boundary
of that for the last entry $\{3, 7, 10, 17, 20\}$ of the table;
they are with (the same!)
$d\!=\!|D_{1,2}|\!=\!17$ and biregular to $\mathbb{A}^{13}$. 

\smallskip

\subsection{\bf Two-flag cells}
The following is the list of all non-affine cells
for $m=1$, i.e. for $D$\~flags $[D_0,D_1]$. We will 
show the $D^\dag$ from Table \ref{Table6-9-22}
by omitting the second entry. They are of importance
since either $D_0$ or $D_1$ are with non-affine
cells if the flag $[D_0,D_1]$ corresponds to a non-affine
one. This matches Proposition \ref{NESTED}. However, there
is plenty of {\em affine} two-flag cells when one of
the ends is non-affine. Similar to Table \ref{Table6-9-22},
the type is that for $N=dim$ and we provide the contributions
to the corresponding {\em geometric\,} superpolynomial.

{\small
\begin{table}[ht!]
\centering
\begin{tabular}{|l|l|c|c|c|c|}
\hline 
$D_0^{\dag}$-sets & $D_1^{\dag}$-sets & dim &  types & 
$q,t,a$-terms \\                         
\hline
 \{11,14,25\} & \{3,11,14,25\}  & 18 & 
$\mathbb{A}^{N}\!\setminus \mathbb{A}^{N\!-\!1}$ & 
$a q^{14} (t^6\!-\! t^7)$\\ 
 \{3,11,14\} & \ \rule{1cm}{0.4pt} & 17 & 
$\mathbb{A}^{N}\!\setminus \mathbb{A}^{N\!-\!1}$ & 
$q^{14} (t^7\!-\! t^8)$\\ 
 \{3,11,14\} & \{3,11,14,16\}  & 18 & 
$\mathbb{A}^{N}\!\setminus \mathbb{A}^{N\!-\!1}$ & 
$a q^{15} (t^6\!-\! t^7)$\\ 
 \{3,11,14\} & \{3,11,14,19\}  & 17 & 
$\mathbb{A}^{N}\!\setminus \mathbb{A}^{N\!-\!1}$ & 
$a q^{15} (t^7\!-\! t^8)$\\ 
 \{10,13,20,23\} & \{3,10,13,20,23\} & 17 & 
$\mathbb{A}^{N}\!\vee \mathbb{A}^{N}$ & 
$a q^{15} (2t^7\!-\! t^8)$\\ 
 \{11,14,19\} & \{3,11,14,19\} & 16 & 
$\mathbb{A}^{N}\!\setminus \mathbb{A}^{N\!-\!1}$ & 
$a q^{15} (t^8\!-\! t^9)$\\ 
 \{3,10,13,20,23\} & \ \rule{1cm}{0.4pt} & 16 & 
$\mathbb{A}^{N}\!\vee \mathbb{A}^{N}$ & 
$q^{15} (2t^8\!-\! t^9)$\\ 
 \{3,10,13,20,23\} & \{3,10,13,14\} & 17 & 
$\mathbb{A}^{N}\!\vee \mathbb{A}^{N}$ & 
$a q^{16}(2t^7\!-\! t^8)$\\ 
 \{3,10,13,20,23\} & \{3,10,13,17,20\} & 16 & 
$\mathbb{A}^{N}\!\vee \mathbb{A}^{N}$ & 
$a q^{16}(2t^8\!-\! t^9)$\\ 
 \{3,11,14,19\} & \ \rule{1cm}{0.4pt}
   & 15 & 
$\mathbb{A}^{N}\!\setminus \mathbb{A}^{N\!-\!1}$ & 
$q^{15} (t^9\!-\!t^{10})$\\ 
 \{3,11,14,19\} & \{3,8,11,19\} & 17 & 
$\mathbb{A}^{N}\!\setminus \mathbb{A}^{N\!-\!1}$ & 
$a q^{16} (t^7\!-\! t^8)$\\ 
 \{3,11,14,19\} & \{3,11,13,14\} & 16 & 
$\mathbb{A}^{N}\!\setminus \mathbb{A}^{N\!-\!1}$ & 
$a q^{16}(t^8\!-\! t^9)$\\ 
 \{3,11,14,19\} & \{3,11,14,16,19\} & 15 & 
$\mathbb{A}^{N}\!\setminus \mathbb{A}^{N\!-\!1}$ & 
$a q^{16}(t^9\!-\! t^{10})$\\ 
 \{8,11,19\} & \{3,8,11,19\}        & 16 & 
$\mathbb{A}^{N}\!\setminus \mathbb{A}^{N\!-\!1}$ & 
$a q^{16}(t^8\!-\! t^9)$\\ 
 \{10,13,17,20\} & \{3,10,13,17,20\} & 15 & 
$\mathbb{A}^{N}\!\vee \mathbb{A}^{N}$ & 
$a q^{16}(2t^9\!-\! t^{10})$\\ 
\{3,8,11,19\} & \ \rule{1cm}{0.4pt}    & 15 & 
$\mathbb{A}^{N}\!\setminus \mathbb{A}^{N\!-\!1}$ & 
$q^{16} (t^9\!-\! t^{10})$\\     
\{3,8,11,19\} & \{3,8,11,13\}     & 16 & 
$\mathbb{A}^{N}\!\setminus \mathbb{A}^{N\!-\!1}$ & 
$a q^{17} (t^8\!-\! t^9)$\\ 
\{3,8,11,19\} & \{3,8,11,16,19\}     & 15 & 
$\mathbb{A}^{N}\!\setminus \mathbb{A}^{N\!-\!1}$ & 
$a q^{17} (t^9\!-\! t^{10})$\\     
\{3,10,13,17,20\} & \ \rule{1cm}{0.4pt} & 14 & 
$\mathbb{A}^{N}\!\vee \mathbb{A}^{N}$ & 
$q^{16} (2t^{10}\!-\! t^{11})$\\ 
\{3,10,13,17,20\} & \{3,7,10,17,20\} & 16 & 
$\mathbb{A}^{N}\!\vee \mathbb{A}^{N}$ & 
$a q^{17} (2t^8\!-\! t^9)$\\ 
\{3,10,13,17,20\} & \{3,10,11,13\} & 15 & 
$\mathbb{A}^{N}\!\vee \mathbb{A}^{N}$ & 
$a q^{17} (2t^9\!-\! t^{10})$\\ 
\{3,10,13,17,20\} & \{3,10,13,14,17\} & 14 & 
$\mathbb{A}^{N}\!\vee \mathbb{A}^{N}$ & 
$a q^{17} (2t^{10}\!-\! t^{11})$\\ 
\{7,10,17,20\} & \{3,7,10,17,20\}        & 15 & 
$\mathbb{A}^{N}\!\vee \mathbb{A}^{N}$ & 
$a q^{17} (2t^9\!-\! t^{10})$\\ 
\{3,7,10,17,20\} & \ \rule{1cm}{0.4pt}     & 14 & 
$\mathbb{A}^{N}\!\vee \mathbb{A}^{N}$ & 
$2q^{17} (t^{10}\!-\! t^{11})$\\ 
\{3,7,10,17,20\} & \{3,7,10,11\}   & 15 & 
$\mathbb{A}^{N}\!\vee \mathbb{A}^{N}$ & 
$a q^{18} (2t^9\!-\! t^{10})$\\ 
\{3,7,10,17,20\} & \{3,7,10,14,17\}   & 14 & 
$\mathbb{A}^{N}\!\vee \mathbb{A}^{N}$ & 
$2a q^{18} (t^{10}\!-\! t^{11})$\\ 
\hline   
\end{tabular}     
\hspace {0.1cm}
\vspace {0.2cm}
\caption{Non-affine $J^{m=1}_{\r}[D_0,\!D_1]$ 
for $\Ga\!=\!\lan 6,9,22\ran$}
\label{Table6-9-22-1}
\end{table}
}

Since there are non-affine 
Piontkowski cells,
(\ref{conjaff}) is not applicable. We checked 
(\ref{conjmot}) and the coincidence
from (\ref{conjcoh}) 
for the coefficients of $a^{0,1}$ and any $q,t$ and 
also under $t=1$ for any powers of $a$ and any $q$,
understanding the admissibility of $\d$
as that of all $D_i$, which is potentially weaker than
the actual admissibility of $\d$
but sufficient for the match with the DAHA superpolynomial.
I.e. the geometric superpolynomials $\h^{h\!om}$ and
$\h^{mod}$ coincide with
$\h_{\,\vec\rr,\,\vec\ss}\,(\square;q,t,a)$ for such
$q,t,a$, quite a confirmation of Conjecture \ref{CONJSUP}.

Recall that (\ref{conjmot}) is
equivalent to (\ref{conjmotxx}) for $\h^{wt}(q,t,a)$
in terms of the weight
filtration. An obvious advantage of modular (\ref{conjmot})
vs. (\ref{conjcoh}) is that we do not need to understand the 
geometry of $J_{\r}[d]$. A possible passage
to obtaining (\ref{conjcoh}) is as follows.  One can try to 
``recombine" the Piontkowski
cells of $J_{\r}[d]$ (within a fixed $d$) to obtain
{\em affine\,} cells; the boundaries of $J_{\r}[D]$
must be known for this. We have a sketch of the corresponding
theory, but it is not finished. In the considered case,
this is doable and 
the resulting complex provides (\ref{conjcoh}).

The modular approach requires only counting  
points of $J_{\r}[\d]$ over $\mathbb{F}_{1/t}$; 
the types from Table \ref{Table6-9-22} are sufficient for
this (knowing the boundaries is unnecessary). Our programs 
determine such types automatically (algebraically)
for ``almost" all cells. We always combine this with 
counting $\mathbb{F}_{1/t}$\~points (for all cells). 
Note that $p=2$ is the only place of {\em bad reduction\,}
here due to $\nu((z^9+z^{13})^2-(z^6)^3)=22.$

A straightforward elimination of the standard $\la$\~parameters
of our modules $M$ 
(see \cite{Pi} and above) gives that only $9$ $D$\~cells can be
potentially non-affine; $6$ of them  
from Table \ref{Table6-9-22} are non-affine indeed.

\subsection{\bf Superpolynomial}\label{sec:sup6-9-22}
Let us provide the corresponding DAHA superpolynomial,
which is for $\vec\rr=\{3,3\},\,
\vec\ss=\{2,4\}$. Recall that
the DAHA construction requires the matrices from
$PSL_2(\Z)$ with the first columns $(\rr_i,\ss_i)^{tr}$,
namely:
{\small $
\begin{pmatrix}3,*\\2,*\end{pmatrix}\,=\,
\tau_+\tau_-^2\,,\ \ 
 \begin{pmatrix}3,*\\4,*\end{pmatrix}\,=\,
\tau_-\tau_+^2\tau_-\,.
$} Then:
$$
\vec\rr=\{3,3\},\,
\vec\ss=\{2,4\},\,
\t= C\!ab(22,3)T(3,2);\ 
\h_{\,\vec\rr,\,\vec\ss}\,(\square\,;\,\,q,t,a)=
$$
\renewcommand{\baselinestretch}{0.5} 
{\small
\(
1+q t+q^2 t+q^3 t+q^4 t+q^5 t+q^2 t^2+q^3 t^2+2 q^4 t^2+2 q^5 t^2
+3 q^6 t^2+2 q^7 t^2+q^8 t^2+q^3 t^3+q^4 t^3+2 q^5 t^3+3 q^6 t^3
+4 q^7 t^3+4 q^8 t^3+4 q^9 t^3+q^{10} t^3+q^4 t^4+q^5 t^4
+2 q^6 t^4+3 q^7 t^4+5 q^8 t^4+5 q^9 t^4+6 q^{10} t^4
+4 q^{11} t^4+q^{12} t^4+q^5 t^5+q^6 t^5+2 q^7 t^5+3 q^8 t^5
+5 q^9 t^5+6 q^{10} t^5+7 q^{11} t^5+6 q^{12} t^5+4 q^{13} t^5
+q^{14} t^5+q^6 t^6+q^7 t^6+2 q^8 t^6+3 q^9 t^6+5 q^{10} t^6
+6 q^{11} t^6+8 q^{12} t^6+7 q^{13} t^6+6 q^{14} t^6+3 q^{15} t^6
+q^7 t^7+q^8 t^7+2 q^9 t^7+3 q^{10} t^7+5 q^{11} t^7+6 q^{12} t^7
+8 q^{13} t^7+8 q^{14} t^7+7 q^{15} t^7+3 q^{16} t^7+q^{17} t^7
+q^8 t^8+q^9 t^8+2 q^{10} t^8+3 q^{11} t^8+5 q^{12} t^8
+6 q^{13} t^8+8 q^{14} t^8+8 q^{15} t^8+8 q^{16} t^8+3 q^{17} t^8
+q^9 t^9+q^{10} t^9+2 q^{11} t^9+3 q^{12} t^9+5 q^{13} t^9
+6 q^{14} t^9+8 q^{15} t^9+8 q^{16} t^9+7 q^{17} t^9+3 q^{18} t^9
+q^{10} t^{10}+q^{11} t^{10}+2 q^{12} t^{10}+3 q^{13} t^{10}
+5 q^{14} t^{10}+6 q^{15} t^{10}+8 q^{16} t^{10}+8 q^{17} t^{10}
+6 q^{18} t^{10}+q^{19} t^{10}+q^{11} t^{11}+q^{12} t^{11}
+2 q^{13} t^{11}+3 q^{14} t^{11}+5 q^{15} t^{11}+6 q^{16} t^{11}
+8 q^{17} t^{11}+7 q^{18} t^{11}+4 q^{19} t^{11}+q^{12} t^{12}
+q^{13} t^{12}+2 q^{14} t^{12}+3 q^{15} t^{12}+5 q^{16} t^{12}
+6 q^{17} t^{12}+8 q^{18} t^{12}+6 q^{19} t^{12}+q^{20} t^{12}
+q^{13} t^{13}+q^{14} t^{13}+2 q^{15} t^{13}+3 q^{16} t^{13}
+5 q^{17} t^{13}+6 q^{18} t^{13}+7 q^{19} t^{13}+4 q^{20} t^{13}
+q^{14} t^{14}+q^{15} t^{14}+2 q^{16} t^{14}+3 q^{17} t^{14}
+5 q^{18} t^{14}+6 q^{19} t^{14}+6 q^{20} t^{14}+q^{21} t^{14}
+q^{15} t^{15}+q^{16} t^{15}+2 q^{17} t^{15}+3 q^{18} t^{15}
+5 q^{19} t^{15}+5 q^{20} t^{15}+4 q^{21} t^{15}+q^{16} t^{16}
+q^{17} t^{16}+2 q^{18} t^{16}+3 q^{19} t^{16}+5 q^{20} t^{16}
+4 q^{21} t^{16}+q^{22} t^{16}+q^{17} t^{17}+q^{18} t^{17}
+2 q^{19} t^{17}+3 q^{20} t^{17}+4 q^{21} t^{17}+2 q^{22} t^{17}
+q^{18} t^{18}+q^{19} t^{18}+2 q^{20} t^{18}+3 q^{21} t^{18}
+3 q^{22} t^{18}+q^{19} t^{19}+q^{20} t^{19}+2 q^{21} t^{19}
+2 q^{22} t^{19}+q^{23} t^{19}+q^{20} t^{20}+q^{21} t^{20}
+2 q^{22} t^{20}+q^{23} t^{20}+q^{21} t^{21}+q^{22} t^{21}
+q^{23} t^{21}+q^{22} t^{22}+q^{23} t^{22}+q^{23} t^{23}
+q^{24} t^{24}
\)

\vfil\smallskip
\noindent
\(
+a^5 \bigl(q^{15}+q^{16} t+q^{17} t+q^{17} t^2
+q^{18} t^2+q^{19} t^2+q^{18} t^3+q^{19} t^3+q^{20} t^3
+q^{21} t^3+q^{19} t^4+q^{20} t^4+q^{21} t^4+q^{20} t^5
+q^{21} t^5+q^{22} t^5+q^{21} t^6+q^{22} t^6+q^{22} t^7
+q^{23} t^7+q^{23} t^8+q^{24} t^9\bigr)
\)

\noindent
\(
+a^4 \bigl(q^{10}
+q^{11}+q^{12}+q^{13}+q^{14}+q^{11} t+2 q^{12} t+3 q^{13} t
+3 q^{14} t+3 q^{15} t+q^{16} t+q^{12} t^2+2 q^{13} t^2
+4 q^{14} t^2+5 q^{15} t^2+5 q^{16} t^2+3 q^{17} t^2+q^{18} t^2
+q^{13} t^3+2 q^{14} t^3+4 q^{15} t^3+6 q^{16} t^3+7 q^{17} t^3
+5 q^{18} t^3+3 q^{19} t^3+q^{20} t^3+q^{14} t^4+2 q^{15} t^4
+4 q^{16} t^4+6 q^{17} t^4+8 q^{18} t^4+6 q^{19} t^4+3 q^{20} t^4
+q^{21} t^4+q^{15} t^5+2 q^{16} t^5+4 q^{17} t^5+6 q^{18} t^5
+8 q^{19} t^5+6 q^{20} t^5+3 q^{21} t^5+q^{16} t^6+2 q^{17} t^6
+4 q^{18} t^6+6 q^{19} t^6+8 q^{20} t^6+5 q^{21} t^6+q^{22} t^6
+q^{17} t^7+2 q^{18} t^7+4 q^{19} t^7+6 q^{20} t^7+7 q^{21} t^7
+3 q^{22} t^7+q^{18} t^8+2 q^{19} t^8+4 q^{20} t^8+6 q^{21} t^8
+5 q^{22} t^8+q^{23} t^8+q^{19} t^9+2 q^{20} t^9+4 q^{21} t^9
+5 q^{22} t^9+3 q^{23} t^9+q^{20} t^{10}+2 q^{21} t^{10}
+4 q^{22} t^{10}+3 q^{23} t^{10}+q^{24} t^{10}+q^{21} t^{11}
+2 q^{22} t^{11}+3 q^{23} t^{11}+q^{24} t^{11}+q^{22} t^{12}
+2 q^{23} t^{12}+q^{24} t^{12}+q^{23} t^{13}+q^{24} t^{13}
+q^{24} t^{14}\bigr)
\)

\vfil
\noindent
\(
+a^3 \bigl(q^6+q^7+2 q^8+2 q^9+2 q^{10}
+q^{11}+q^{12}+q^7 t+2 q^8 t+4 q^9 t+6 q^{10} t+7 q^{11} t
+6 q^{12} t+4 q^{13} t+2 q^{14} t+q^8 t^2+2 q^9 t^2+5 q^{10} t^2
+8 q^{11} t^2+12 q^{12} t^2+12 q^{13} t^2+10 q^{14} t^2
+5 q^{15} t^2+2 q^{16} t^2+q^9 t^3+2 q^{10} t^3+5 q^{11} t^3
+9 q^{12} t^3+14 q^{13} t^3+17 q^{14} t^3+16 q^{15} t^3
+11 q^{16} t^3+5 q^{17} t^3+2 q^{18} t^3+q^{10} t^4+2 q^{11} t^4
+5 q^{12} t^4+9 q^{13} t^4+15 q^{14} t^4+19 q^{15} t^4
+21 q^{16} t^4+16 q^{17} t^4+9 q^{18} t^4+3 q^{19} t^4+q^{20} t^4
+q^{11} t^5+2 q^{12} t^5+5 q^{13} t^5+9 q^{14} t^5+15 q^{15} t^5
+20 q^{16} t^5+23 q^{17} t^5+19 q^{18} t^5+10 q^{19} t^5
+3 q^{20} t^5+q^{12} t^6+2 q^{13} t^6+5 q^{14} t^6+9 q^{15} t^6
+15 q^{16} t^6+20 q^{17} t^6+24 q^{18} t^6+19 q^{19} t^6
+9 q^{20} t^6+2 q^{21} t^6+q^{13} t^7+2 q^{14} t^7+5 q^{15} t^7
+9 q^{16} t^7+15 q^{17} t^7+20 q^{18} t^7+23 q^{19} t^7
+16 q^{20} t^7+5 q^{21} t^7+q^{14} t^8+2 q^{15} t^8+5 q^{16} t^8
+9 q^{17} t^8+15 q^{18} t^8+20 q^{19} t^8+21 q^{20} t^8
+11 q^{21} t^8+2 q^{22} t^8+q^{15} t^9+2 q^{16} t^9+5 q^{17} t^9
+9 q^{18} t^9+15 q^{19} t^9+19 q^{20} t^9+16 q^{21} t^9
+5 q^{22} t^9+q^{16} t^{10}+2 q^{17} t^{10}+5 q^{18} t^{10}
+9 q^{19} t^{10}+15 q^{20} t^{10}+17 q^{21} t^{10}
+10 q^{22} t^{10}+2 q^{23} t^{10}+q^{17} t^{11}+2 q^{18} t^{11}
+5 q^{19} t^{11}+9 q^{20} t^{11}+14 q^{21} t^{11}
+12 q^{22} t^{11}+4 q^{23} t^{11}+q^{18} t^{12}+2 q^{19} t^{12}
+5 q^{20} t^{12}+9 q^{21} t^{12}+12 q^{22} t^{12}+6 q^{23} t^{12}
+q^{24} t^{12}+q^{19} t^{13}+2 q^{20} t^{13}+5 q^{21} t^{13}
+8 q^{22} t^{13}+7 q^{23} t^{13}+q^{24} t^{13}+q^{20} t^{14}
+2 q^{21} t^{14}+5 q^{22} t^{14}+6 q^{23} t^{14}+2 q^{24} t^{14}
+q^{21} t^{15}+2 q^{22} t^{15}+4 q^{23} t^{15}+2 q^{24} t^{15}
+q^{22} t^{16}+2 q^{23} t^{16}+2 q^{24} t^{16}+q^{23} t^{17}
+q^{24} t^{17}+q^{24} t^{18}\bigr)
\)

\vfil
\noindent
\(
+a^2 \bigl(q^3+q^4+2 q^5
+2 q^6+2 q^7+q^8+q^9+q^4 t+2 q^5 t+4 q^6 t+6 q^7 t+8 q^8 t
+7 q^9 t+6 q^{10} t+3 q^{11} t+q^{12} t+q^5 t^2+2 q^6 t^2
+5 q^7 t^2+8 q^8 t^2+13 q^9 t^2+15 q^{10} t^2+15 q^{11} t^2
+10 q^{12} t^2+5 q^{13} t^2+q^{14} t^2+q^6 t^3+2 q^7 t^3
+5 q^8 t^3+9 q^9 t^3+15 q^{10} t^3+20 q^{11} t^3+24 q^{12} t^3
+19 q^{13} t^3+12 q^{14} t^3+5 q^{15} t^3+q^{16} t^3+q^7 t^4
+2 q^8 t^4+5 q^9 t^4+9 q^{10} t^4+16 q^{11} t^4+22 q^{12} t^4
+29 q^{13} t^4+28 q^{14} t^4+21 q^{15} t^4+11 q^{16} t^4
+4 q^{17} t^4+q^{18} t^4+q^8 t^5+2 q^9 t^5+5 q^{10} t^5
+9 q^{11} t^5+16 q^{12} t^5+23 q^{13} t^5+31 q^{14} t^5
+33 q^{15} t^5+29 q^{16} t^5+16 q^{17} t^5+6 q^{18} t^5
+q^{19} t^5+q^9 t^6+2 q^{10} t^6+5 q^{11} t^6+9 q^{12} t^6
+16 q^{13} t^6+23 q^{14} t^6+32 q^{15} t^6+35 q^{16} t^6
+32 q^{17} t^6+19 q^{18} t^6+6 q^{19} t^6+q^{20} t^6+q^{10} t^7
+2 q^{11} t^7+5 q^{12} t^7+9 q^{13} t^7+16 q^{14} t^7
+23 q^{15} t^7+32 q^{16} t^7+36 q^{17} t^7+32 q^{18} t^7
+16 q^{19} t^7+4 q^{20} t^7+q^{11} t^8+2 q^{12} t^8+5 q^{13} t^8
+9 q^{14} t^8+16 q^{15} t^8+23 q^{16} t^8+32 q^{17} t^8
+35 q^{18} t^8+29 q^{19} t^8+11 q^{20} t^8+q^{21} t^8+q^{12} t^9
+2 q^{13} t^9+5 q^{14} t^9+9 q^{15} t^9+16 q^{16} t^9
+23 q^{17} t^9+32 q^{18} t^9+33 q^{19} t^9+21 q^{20} t^9
+5 q^{21} t^9+q^{13} t^{10}+2 q^{14} t^{10}+5 q^{15} t^{10}
+9 q^{16} t^{10}+16 q^{17} t^{10}+23 q^{18} t^{10}
+31 q^{19} t^{10}+28 q^{20} t^{10}+12 q^{21} t^{10}+q^{22} t^{10}
+q^{14} t^{11}+2 q^{15} t^{11}+5 q^{16} t^{11}+9 q^{17} t^{11}
+16 q^{18} t^{11}+23 q^{19} t^{11}+29 q^{20} t^{11}
+19 q^{21} t^{11}+5 q^{22} t^{11}+q^{15} t^{12}+2 q^{16} t^{12}
+5 q^{17} t^{12}+9 q^{18} t^{12}+16 q^{19} t^{12}
+22 q^{20} t^{12}+24 q^{21} t^{12}+10 q^{22} t^{12}
+q^{23} t^{12}+q^{16} t^{13}+2 q^{17} t^{13}+5 q^{18} t^{13}
+9 q^{19} t^{13}+16 q^{20} t^{13}+20 q^{21} t^{13}
+15 q^{22} t^{13}+3 q^{23} t^{13}+q^{17} t^{14}+2 q^{18} t^{14}
+5 q^{19} t^{14}+9 q^{20} t^{14}+15 q^{21} t^{14}
+15 q^{22} t^{14}+6 q^{23} t^{14}+q^{18} t^{15}+2 q^{19} t^{15}
+5 q^{20} t^{15}+9 q^{21} t^{15}+13 q^{22} t^{15}+7 q^{23} t^{15}
+q^{24} t^{15}+q^{19} t^{16}+2 q^{20} t^{16}+5 q^{21} t^{16}
+8 q^{22} t^{16}+8 q^{23} t^{16}+q^{24} t^{16}+q^{20} t^{17}
+2 q^{21} t^{17}+5 q^{22} t^{17}+6 q^{23} t^{17}+2 q^{24} t^{17}
+q^{21} t^{18}+2 q^{22} t^{18}+4 q^{23} t^{18}+2 q^{24} t^{18}
+q^{22} t^{19}+2 q^{23} t^{19}+2 q^{24} t^{19}+q^{23} t^{20}
+q^{24} t^{20}+q^{24} t^{21}\bigr)
\)

\vskip 0.1cm
\noindent
\(+a \bigl(q+q^2+q^3+q^4+q^5
+q^2 t+2 q^3 t+3 q^4 t+4 q^5 t+5 q^6 t+4 q^7 t+2 q^8 t+q^9 t
+q^3 t^2+2 q^4 t^2+4 q^5 t^2+6 q^6 t^2+9 q^7 t^2+10 q^8 t^2
+9 q^9 t^2+5 q^{10} t^2+2 q^{11} t^2+q^4 t^3+2 q^5 t^3+4 q^6 t^3
+7 q^7 t^3+11 q^8 t^3+14 q^9 t^3+16 q^{10} t^3+13 q^{11} t^3
+6 q^{12} t^3+2 q^{13} t^3+q^5 t^4+2 q^6 t^4+4 q^7 t^4+7 q^8 t^4
+12 q^9 t^4+16 q^{10} t^4+20 q^{11} t^4+20 q^{12} t^4
+14 q^{13} t^4+6 q^{14} t^4+2 q^{15} t^4+q^6 t^5+2 q^7 t^5
+4 q^8 t^5+7 q^9 t^5+12 q^{10} t^5+17 q^{11} t^5+22 q^{12} t^5
+24 q^{13} t^5+21 q^{14} t^5+13 q^{15} t^5+4 q^{16} t^5
+q^{17} t^5+q^7 t^6+2 q^8 t^6+4 q^9 t^6+7 q^{10} t^6
+12 q^{11} t^6+17 q^{12} t^6+23 q^{13} t^6+26 q^{14} t^6
+25 q^{15} t^6+17 q^{16} t^6+7 q^{17} t^6+q^{18} t^6+q^8 t^7
+2 q^9 t^7+4 q^{10} t^7+7 q^{11} t^7+12 q^{12} t^7+17 q^{13} t^7
+23 q^{14} t^7+27 q^{15} t^7+27 q^{16} t^7+18 q^{17} t^7
+7 q^{18} t^7+q^{19} t^7+q^9 t^8+2 q^{10} t^8+4 q^{11} t^8
+7 q^{12} t^8+12 q^{13} t^8+17 q^{14} t^8+23 q^{15} t^8
+27 q^{16} t^8+27 q^{17} t^8+17 q^{18} t^8+4 q^{19} t^8
+q^{10} t^9+2 q^{11} t^9+4 q^{12} t^9+7 q^{13} t^9+12 q^{14} t^9
+17 q^{15} t^9+23 q^{16} t^9+27 q^{17} t^9+25 q^{18} t^9
+13 q^{19} t^9+2 q^{20} t^9+q^{11} t^{10}+2 q^{12} t^{10}
+4 q^{13} t^{10}+7 q^{14} t^{10}+12 q^{15} t^{10}
+17 q^{16} t^{10}+23 q^{17} t^{10}+26 q^{18} t^{10}
+21 q^{19} t^{10}+6 q^{20} t^{10}+q^{12} t^{11}+2 q^{13} t^{11}
+4 q^{14} t^{11}+7 q^{15} t^{11}+12 q^{16} t^{11}
+17 q^{17} t^{11}+23 q^{18} t^{11}+24 q^{19} t^{11}
+14 q^{20} t^{11}+2 q^{21} t^{11}+q^{13} t^{12}+2 q^{14} t^{12}
+4 q^{15} t^{12}+7 q^{16} t^{12}+12 q^{17} t^{12}
+17 q^{18} t^{12}+22 q^{19} t^{12}+20 q^{20} t^{12}
+6 q^{21} t^{12}+q^{14} t^{13}+2 q^{15} t^{13}+4 q^{16} t^{13}
+7 q^{17} t^{13}+12 q^{18} t^{13}+17 q^{19} t^{13}
+20 q^{20} t^{13}+13 q^{21} t^{13}+2 q^{22} t^{13}+q^{15} t^{14}
+2 q^{16} t^{14}+4 q^{17} t^{14}+7 q^{18} t^{14}+12 q^{19} t^{14}
+16 q^{20} t^{14}+16 q^{21} t^{14}+5 q^{22} t^{14}+q^{16} t^{15}
+2 q^{17} t^{15}+4 q^{18} t^{15}+7 q^{19} t^{15}+12 q^{20} t^{15}
+14 q^{21} t^{15}+9 q^{22} t^{15}+q^{23} t^{15}+q^{17} t^{16}
+2 q^{18} t^{16}+4 q^{19} t^{16}+7 q^{20} t^{16}+11 q^{21} t^{16}
+10 q^{22} t^{16}+2 q^{23} t^{16}+q^{18} t^{17}+2 q^{19} t^{17}
+4 q^{20} t^{17}+7 q^{21} t^{17}+9 q^{22} t^{17}+4 q^{23} t^{17}
+q^{19} t^{18}+2 q^{20} t^{18}+4 q^{21} t^{18}+6 q^{22} t^{18}
+5 q^{23} t^{18}+q^{20} t^{19}+2 q^{21} t^{19}+4 q^{22} t^{19}
+4 q^{23} t^{19}+q^{24} t^{19}+q^{21} t^{20}+2 q^{22} t^{20}
+3 q^{23} t^{20}+q^{24} t^{20}+q^{22} t^{21}+2 q^{23} t^{21}
+q^{24} t^{21}+q^{23} t^{22}+q^{24} t^{22}+q^{24} t^{23}\bigr).
\)
}
\renewcommand{\baselinestretch}{1.2} 
\vskip 0.2cm

Concerning practical aspects, the production of this DAHA
superpolynomial requires a couple of minutes. About the same
time is needed to calculate all dimensions of $J_{\r}[D]$,
including the list of non-admissible modules and 
potentially non-affine cells. 
Such a calculation with $J_{\r}^{m=1}[\d]$ takes about 
10 minutes.

\subsection {\bf Exponents (6,9,14),(6,9,16)} 
An example of a cell which type is different from those in Table 
\ref{Table6-9-22} is for $\r=\C[[z^6,z^9+z^{14}]]$, where
$\Ga=\lan 6,9,23\ran,\de= 25$ and the 
corresponding link is $C\!ab(23,3)T(3,2)$.
Namely, there is exactly one cell $J_{\r}[D]$ for 
$D^\dag=\{3,10,13,20\}$ that is biregular to 
$(\mathbb{A}^{15}\setminus \mathbb{A}^{14})\cup
(\mathbb{A}^{15}\setminus \mathbb{A}^{14})$;
here $|D|=16$, the union is disjoint and 
the contribution to the superpolynomial
is $2a q^{16}(t^9-t^{10})$. This is the
simplest {\em disconnected\,} $D$-cell we found.
For $2$-flags, new types of cells are of the
same kind $(\mathbb{A}^{N}\setminus \mathbb{A}^{N\!-\!1})\cup
(\mathbb{A}^{N}\setminus \mathbb{A}^{N\!-\!1})$;
 they are:
{\small
\begin{align*}
&D_0^\dag=\{10,13,20\},&D_1^\dag=\{3,10,13,20\},\ &dim\!=\!17,
&2a q^{16} (t^8-t^9),\\
&D_0^\dag=\{3,10,13,20\},&D_1^\dag=\{3,10,13,14\},\ &dim\!=\!17,
&2a q^{17} (t^8-t^9),\\
&D_0^\dag=\{3,10,13,20\},&D_1^\dag=\{3,10,13,17,20\},\ 
&dim\!=\!16,
&2a q^{17} (t^9\!-\!t^{10}).
\end{align*}
}
Also, there are
$4+5$ non-affine "old" $D$\~cells of types 
$\mathbb{A}^{N}\!\setminus \mathbb{A}^{N\!-\!1}$ and
$\mathbb{A}^{N}\!\vee \mathbb{A}^{N}$. 
The total number of non-affine cells for
$m=0,1$ (contributing to  
$a^{0,1}$) is $44$. 
It takes
about 5 min  to calculate all dim$J_{\r}[D]$ and about
30 min for obtaining all dim$J_{\r}^{m=1}[D_0,D_1]$.
The computer program almost always finds the types of 
non-affine cells (reporting questionable cases).  The match with 
the DAHA-superpolynomial for the coefficients of $a^{0,1}$ or 
under $t=1$ (any $q,a$) is perfect.
\smallskip

In the case $\r=\C[[z^6,z^9+z^{16}]]$, where
$\Ga=\lan 6,9,25\ran,\de= 27$ and the 
corresponding link is $C\!ab(25,3)T(3,2)$, 
(exactly) one new type of cells $J_{\r}[D]$ appears.
There is a unique cell of the following kind:
\begin{align*}
&D^\dag=\{3,10,13,20,23\}, \
 |D|\!=\!18,\, dim\!=\!16, \  J_{\r}[D]\cong A_1\cup A_2 
\cup A_3,\\
&A_i\cong\mathbb{A}^{16}, A_1\cap A_2\cong \mathbb{A}^{15}\cong
A_2\cap A_3, A_1\cap A_2\cap A_3=A_1\cap A_3\cong \mathbb{A}^{14}.
\end{align*}
Accordingly,  its contribution to the geometric
superpolynomial equals $a q^{18}(3 t^{11}-2 t^{12})$;\, 
$11=\de-16=
27-16$. The new types for $2$\~flags are of the
same kind; they (and their contributions)
are as follows:
{\small
\begin{align*}
&D_0^\dag\!=\!\{10,13,20,23\},\!\!&D_1^\dag\!=\!\{3,10,13,20,23\},
 &\ dim\!=\!17,
&a q^{18} (3t^{10}\!-\!2t^{11}),\\
&D_0^\dag\!=\!\{3,10,13,20,23\},\!\!&D_1^\dag\!=\!\{3,10,13,14\},
 &\ dim\!=\!17,
&a q^{19} (3t^{10}\!-\!2t^{11}),\\
&D_0^\dag\!=\!\{3,10,13,20,23\},\!\!&D_1^\dag\!=
\!\{3,10,13,17,20\}, 
&\ dim\!=\!16, &a q^{19} (3t^{11}\!-\!2t^{12}). 
\end{align*}
} 
The number of all non-affine cells is $\,62\,$ ($m=0,1$) in this 
case; the match with the DAHA-superpolynomial at 
$a^{0,1}$ is perfect for this $\r$.
The similarity with the previous $1+3$ ``new cells" 
for $\C[[z^6,z^9+z^{14}]]$
of type $2(\mathbb{A}^{N}\setminus \mathbb{A}^{N\!-\!1})$
is hardly accidental; indeed, we just add $23$ to $D_{0,1}$.

We performed the same calculation (and the check vs. DAHA) for 
$\C[[z^6,z^9+z^{17}]]$ (which took about 300 min).
The total number of non-affine cells becomes $\,85\,$ ($m=0,1$)
and no ``new types" appear vs. the previous $\r$.
The total number of $\d$\~flags for $m=0,1$ is $3102$. Due
to our extensive numerical experiments, we expect that 
non-affine Piontkowski cells can be only as described above 
($4$ types) for the whole family
$\r= \C[[z^6,z^9+z^{3p\pm 1}]]$ with any $p\ge 4$.

\setcounter{equation}{0}
\section{\sc Dimensions for (6,9,13)}
\label{app:dimens}
In the case of $\r=\C[[z^6,z^9+z^{13}]]$
with $\Ga\!=\!\lan 6,9,22\ran$,
we will provide all dimensions dim$J_{\r}^{m=0,1}[\d]$.
The corresponding exponents are $(6,9,13)$ in this case. 
Importantly, the cells $J_{\r}^{m=0,1}[\d]$
are not always affine; use Table \ref{Table6-9-22-1} 
for $\Ga\!=\!\lan 6,9,22\ran$ for the list of non-affine cells 
in this case. The types and the ``generic" dimension
are sufficient to determine the corresponding contribution
to the geometric superpolynomial.

The dimension tables for such $m=0,1$ (and the corresponding 
full DAHA superpolynomials) are also available upon request for 
$\r=\C[[z^6,z^9+z^{14,16}]]$. 
We also constantly calculate the geometric superpolynomials
for any $m$ under $q=1$ in these and all considered
examples, which will not be provided.

Only $D_0^\dag$ or $D_1^\dag$ are shown in the table followed
by the corresponding dimensions and  
$|D_0|$ (after dim$J_{\r}[D_0]$) for the lines without "$+$". 
We note that $|D_0|$ is needed to be shown since it is 
not immediate to calculate it in terms of primitive $D_0^\dag$,
though the latter of course uniquely determines $D_0$. 
For the lines with "$+$", we put $D_1^\dag$ followed by  
dim$J_{\r}^{m=1}[D_0,D_1]$; recall that $|D_1|=|D_0|+1$.
The corresponding $D_0^\dag$ must be taken from 
the {\em closest previous entry\,} without "$+$".

We fill the first row, then the second and so on; by {\em na\,}
we mean {\em non-admissible\,} $D$\~flags (which do not
contribute to the superpolynomial). 

Let us mention that Tables \ref{Table6-9-22},
\ref{Table6-9-22-1} are actually special cases of the
the below table. For instance,  dim$J_{\r}[D]=16$
for $D$ corresponding to $D^\dag=\{3,10,13,20,23\}$ 
(find below this $D^\dag$ in lines without $+$),
which contribute $q^{15}(2t^8-t^9)$ according to 
Tables \ref{Table6-9-22}, \ref{Table6-9-22-1}. If this cell
were affine, it would result in pure 
$q^{15} t^8$, where $8=\delta-\dim J_{\r}[D]=24-16$. 
\smallskip

The match with the corresponding coefficients
of DAHA superpolynomial is perfect, as well as in the
case of $\Ga\!=\!\lan 6,9,23\ran$ (which will not be
discussed here). The dimensions are
provided below. Note that there are exactly $5$ 
pairs $\{D_0,D_1\}$ of maximal cell-dimension $\de=24$:
\vskip -0.5cm

{\small
\begin{align*}
&\{[\emptyset],[47]\}\ ,\
\{[47], [38, 47]\}\ ,\ \{[38, 47], [25, 38, 47]\},\\ 
&\{[25, 38, 47], [16, 25, 38, 47]\}\ ,\  
\{[16, 25, 38, 47], [3, 16, 25, 38, 47]\},
\end{align*}
}
\vskip -0.5cm

\noindent
which contribute $a(q+q^2+q^3+q^4+q^5)$ to the superpolynomial
from Section \ref{sec:sup6-9-22}. In terms of the primitive
$D^\dag$ in the flags (used
in the table below instead of the ``complete" $D$) they
are: 

{\small
\begin{align*}
&\bigl\{\{\,\},\{47\}\bigr\},
\bigl\{\{47\},\{38\}\bigr\}, \bigl\{\{38\}, \{25,38\}\bigr\}, 
\bigl\{\{25, 38\}, \{16\}\bigr\}, \bigl\{\{16\},\{3,16\}\bigr\}.
\end{align*}
}

The table of dimensions for $\Ga\!=\!\lan 6,9,22\ran$
is as follows:

\comment{
{\tiny
\begin{table}[ht!]
\centering
\begin{longtable}{|c|c|c|c|}
\hline 
$D_0^{\dag}$-sets & $D_1^{\dag}$-sets & dim &  types \\     
\hline                
\{ 3,11,14\}\ :\ 17 & 
+\{ 3,5,11,14\}:\ na & 
+\{ 3,8,11,14\}:\ na & 
+\{ 3,11,14,16\}:\ 18\\ 
\hline   
\end{longtable}     
\hspace {0.1cm}
\vspace {0.2cm}
\caption{Test $J^{m=1}_{\r}[D_0,\!D_1]$ 
for $\Ga\!=\!\lan 6,9,22\ran$}
\label{Table6-9-22-x}
\end{table}
}}
\input D6-9-22.tex

\end{document}

%% file: macro.tex
\renewcommand{\tilde}{\widetilde}
\renewcommand{\hat}{\widehat}

\newcommand{\BR}{{\mathbb R}}
\newcommand{\BQ}{{\mathbb Q}}
\newcommand{\BC}{{\mathbb C}}
\newcommand{\BP}{{\mathbb P}}
\newcommand{\BZ}{{\mathbb Z}}
\newcommand{\BN}{{\mathbb N}}
\newcommand{\BS}{{\mathbb S}}

\newcommand{\cH}{{\mathcal H}}
\newcommand{\cA}{{\mathcal A}}
\newcommand{\cB}{{\mathcal B}}
\newcommand{\ccF}{{\mathfrak F}}
\newcommand{\cD}{{\mathcal D}}
\newcommand{\cL}{{\mathcal L}}
\newcommand{\cF}{{\mathcal F}}
\newcommand{\cP}{{\mathcal P}}
\newcommand{\cX}{{\mathcal X}}
\newcommand{\cY}{{\mathcal Y}}
\newcommand{\cS}{{\mathcal S}}
\newcommand{\cSol}{\hbox{$\mathcal Sol$}}
\newcommand{\cT}{\hbox{$\mathcal T$}}

\newcommand{\Z}{{\mathbb Z}}
\newcommand{\Q}{{\mathbb Q}}
\newcommand{\N}{{\mathbb N}}
\newcommand{\C}{{\mathbb C}}
\newcommand{\R}{{\mathbb R}}
\newcommand{\X}{{\mathbb X}}
\newcommand{\Y}{{\mathbb Y}}

\newcommand{\CH}{{\mathcal H}}
\newcommand{\CA}{{\mathcal A}}

\def\HH{\mbox{${\mathcal H}$\kern-5.2pt${\mathcal H}$}}

\newcommand{\binomial}[2]{\genfrac{(}{)}{0pt}{}{ #1 }{ #2 }}
\newcommand{\qbinomial}[2]{\genfrac{[}{]}{0pt}{}{ #1 }{ #2 }_q }
\newcommand{\qbinom}[3]{\genfrac{[}{]}{0pt}{}{ #1 }{ #2 }_{ #3 } }


\def\der{\partial}
\def\tensor{\otimes}
\def\gam{\gamma} \def\Gam{\Gamma}
\def\del{\delta} \def\Del{\Delta}
\def\kap{\kappa}
\def\lam{\lambda} \def\Lam{\Lambda}
\def\Comp{{\mathbb C}}
\def\sM{{\mathcal M}}

\newtheorem{theorem}{Theorem}[section]
\newtheorem{maintheorem}[theorem]{Main Theorem}
\newtheorem{proposition}[theorem]{Proposition}
\newtheorem{definition}[theorem]{Definition}
\newtheorem{lemma}[theorem]{Lemma}
\newtheorem{corollary}[theorem]{Corollary}
\newtheorem{notation}[theorem]{Notation}
\newtheorem{remark}[theorem]{Remark}
\newtheorem{example}[theorem]{Example}

\newtheorem{theorem }{Theorem}[section]
\newtheorem{maintheorem }[theorem]{Main Theorem}
\newtheorem{proposition }[theorem]{Proposition}
\newtheorem{definition }[theorem]{Definition}
\newtheorem{lemma }[theorem]{Lemma}
\newtheorem{corollary }[theorem]{Corollary}
\newtheorem{notation }[theorem]{Notation}
\newtheorem{remark }[theorem]{Remark}
\newtheorem{example }[theorem]{Example}

\newtheorem{ maintheorem }[theorem]{Main Theorem}
\newtheorem{ theorem}{Theorem}[section]
\newtheorem{ proposition}[theorem]{Proposition}
\newtheorem{ definition}[theorem]{Definition}
\newtheorem{ lemma}[theorem]{Lemma}
\newtheorem{ corollary}[theorem]{Corollary}
\newtheorem{ notation}[theorem]{Notation}
\newtheorem{ remark}[theorem]{Remark}
\newtheorem{ example}[theorem]{Example}

\newtheorem{thm}{Theorem}[section]
\newtheorem{prop}[thm]{Proposition}
\newtheorem{lem}[thm]{Lemma}
\newtheorem{cor}[thm]{Corollary}
\newtheorem{conj}[thm]{Conjecture}
\newtheorem{con}[thm]{Conjecture}
\newtheorem{dfn}[thm]{Definition}
\newtheorem{df}[thm]{Definition}
 \newcommand{\rem}{{\bf Comment.\ }}
 \newcommand{\rmk}{{\bf Comment.\ }}
 \newcommand{\exmp}{{\bf Example.\ }}
 \newcommand{\ex}{{\bf Example.\ }}
 \newcommand{\prob}{{\bf Problem.\ }}

\newtheorem{note}{Note} 
\renewcommand{\thenote}{}
\newtheorem*{acka}{Acknowledgments}
\newtheorem{ack}{Acknowledgments}
\renewcommand{\theack}{}
\renewcommand{\appendixname}{\bf Appendix}
\renewcommand{\proof}{{\em Proof.\ }}

\hyphenation{
ap-pen-dix as-ymp-tot-ic at-trib-uted at-trib-ut-able
Bry-li-n-sky com-mu-ta-tion de-ge-ne-rate
de-riv-a-tive dis-trib-ute equi-vari-ant ex-tra-or-di-nary  
geo-met-ric griev-ance griev-ous grad-ed ho-lo-no-my ho-mo-thetic
in-fin-ite-ly in-fin-i-tes-i-mal Ha-rish Cha-n-dra mul-ti-plic-able 
non-euclid-ean non-iso-mor-phic non-smooth par-a-digm 
par-a-bol-ic pa-rab-o-loid pa-ram-e-trize phe-nom-e-non 
post-script pseu-do-dif-fer-en-tial pseu-do-fi-nite 
qua-drat-ics quad-ra-ture Han-kel rec-tan-gle semi-def-i-nite 
set-up wide-spread Euler-ian Feb-ru-ary Gauss-ian Grothen-dieck 
Hamil-ton-ian Her-mi-t-ian her-mi-t-ian Jan-u-ary 
Japan-ese Ka-shi-wa-ra Kor-te-weg Le-gendre No-vem-ber Rie-mann-ian 
Sep-tem-ber Za-mo-lo-d-chi-kov Kni-zh-nik quan-tum Op-dam
Mac-do-nald Ca-lo-ge-ro Su-ther-land Mo-ser 
Ol-sha-net-sky  Pe-re-lo-mov in-de-pen-dent ope-ra-tors 
cy-clo-to-mic ra-tio-nal de-gen-er-a-tion 
in-ter-est-ing de-for-ma-tions de-for-ma-tion pro-ce-dure 
fol-lows ope-ra-tors  pre-serve suf-fices ap-proach 
for-mu-las con-sider its com-ple-tion cor-re-spond-ing 
au-to-mor-phism be-cause pro-por-tional fi-nal-ly let-ting 
equi-v-a-lence ge-n-er-al-ized Mac-do-nald iden-ti-ties 
cor-re-s-pond sub-dia-grams par-ti-tion na-t-u-ral-ly 
or-dered stan-dard de-for-ma-tion ar-gu-ment com-bined 
sphe-r-i-cal rep-re-sen-ta-tions tri-go-no-me-t-ric
ge-n-er-al-ly speak-ing pri-m-it-ive ir-re-du-cible 
sum-ma-tion  rep-re-sen-ta-tives pro-por-ti-o-na-li-ty
ultra-sphe-ri-cal Ro-gers}

\def\ffor{\quad\hbox{ for }\quad}
\def\wwhen{\quad\hbox{ when }\quad}
\def\wwhere{\quad\hbox{ where }\quad}
\def\aand{\quad\hbox{ and }\quad}
\def\for{\  \hbox{ for } \ }
\def\iif{ \ \hbox{ if } \ }
\def\when{ \ \hbox{ when } \ }
\def\where{\  \hbox{ where } \ }
\def\and{\  \hbox{ and } \ }
\def\and{\  \hbox{ and } \ }
\def\oor{\  \hbox{ or } \ }
\def\proof{{\em Proof. \  }}

\def\equal{\stackrel{\,\mathbf{def}}{= \kern-3pt =}}

\def\la{\lambda}
\def\La{\Lambda}
\def\om{\omega}
\def\Om{\Omega}
\def\Th{\Theta}
\def\th{\theta}
\def\al{\alpha}
\def\be{\beta}
\def\ga{\gamma}
\def\ep{\epsilon}
\def\up{\upsilon}
\def\Up{\Upsilon}
\def\de{\delta}
\def\De{\Delta}
\def\ka{\kappa}
\def\kapp{\hbox{\bf \ae}}
\def\si{\sigma}
\def\Si{\Sigma}
\def\Ga{\Gamma}
\def\ze{\zeta}
\def\io{\iota}
\def\bio{b^\iota}
\def\aio{a^\iota}
\def\twio{\tilde{w}^\iota}
\def\hwio{\hat{w}^\iota}
\def\gio{\g^\iota}
\def\Bio{B^\iota}

\def\del{\delta}
\def\pa{\partial}
\def\vp{\varphi}
\def\ve{\varepsilon}
\def\inf{\infty}

\def\vph{\varphi}
\def\vps{\varpsi}
\def\vPh{\varPhi}
\def\vep{\varepsilon}
\def\vpi{{\varpi}}
\def\vth{{\vartheta}}
\def\vsi{{\varsigma}}
\def\vrh{{\varrho}}

\def\bph{\bar{\phi}}
\def\bsi{\bar{\si}}
\def\bvp{\bar{\varphi}}

\newcommand{\bS}{{\mathbf S}}
\newcommand{\bH}{{\mathbf H}}
\newcommand{\bF}{{\mathbf F}}
\newcommand{\bE}{{\mathbf E}}

\def\tal{\tilde{\alpha}}
\def\tbe{\tilde{\beta}}
\def\tde{\tilde{\delta}}
\def\tpi{\tilde{\pi}}
\def\txi{\tilde{\xi}}
\def\tPi{\tilde{\Pi}}
\def\tPhi{\tilde{\Phi}}
\def\tV{\tilde{V}}
\def\tJ{\tilde{J}}
\def\tla{\tilde{\lambda}}
\def\tga{\tilde{\gamma}}
\def\tGa{\tilde{\Gamma}}
\def\tvs{\tilde{{\varsigma}}}
\def\tu{\tilde{u}}
\def\tU{\tilde{U}}
\def\tw{\widetilde w}
\def\tW{\widetilde W}
\def\tB{\tilde B}
\def\tv{\tilde v}
\def\tV{\tilde V}
\def\tz{\tilde z}
\def\tb{\tilde b}
\def\ta{\tilde a}
\def\tih{\tilde h}
\def\trh{\tilde {\rho}}
\def\tx{\tilde x}
\def\tf{\tilde f}
\def\tg{\tilde g}
\def\tG{\tilde G}
\def\tk{\tilde k}
\def\tl{\tilde l}
\def\tL{\tilde L}
\def\tD{\tilde D}
\def\tR{\tilde R}
\def\tP{\tilde P}
\def\tH{\tilde H}
\def\tp{\tilde p}

\def\hH{\hat{H}}
\def\hh{\hat{h}}
\def\hR{\hat{R}}
\def\hY{\hat{Y}}
\def\hX{\hat{X}}
\def\hP{\hat{P}}
\def\hT{\hat{T}}
\def\hV{\hat{V}}
\def\hG{\hat{G}}
\def\hF{\hat{F}}
\def\hw{\widehat{w}}
\def\hW{\widehat{W}}
\def\hu{\hat{u}}
\def\hs{\hat{s}}
\def\hv{\hat{v}}
\def\hb{\hat{b}}
\def\hB{\widehat{B}}
\def\hze{\hat{\zeta}}
\def\hsi{\hat{\sigma}}
\def\hrh{\hat{\rho}}
\def\hth{\hat{\theta}}
\def\hy{\hat{y}}
\def\hx{\hat{x}}
\def\hz{\hat{z}}
\def\hg{\hat{g}}
\def\he{\hat{e}}
\def\hE{\widehat{E}}

\def\B{\mathbf{B}}
\def\I{\mathbf{I}}
\def\P{\mathbf{P}}
\def\G{\mathbf{G}}
\def\S{\mathbf{S}}
\def\F{\mathbf{F}}
\def\one{\mathbf{1}}
\def\Sn{\mathbf{S}_n}
\def\0{\mathbf{0}}
\def\H{\mathbf{H}}
\def\V{\mathbf{V}}

\def\f{\mathcal{F}}
\def\çF{\mathcal{F}}
\def\o{\mathcal{O}}
\def\t{\mathcal{T}}
\def\r{\mathcal{R}}
\def\l{\mathcal{L}}
\def\m{\mathcal{M}}
\def\k{\mathcal{K}}
\def\n{\mathcal{N}}
\def\d{\mathcal{D}}
\def\p{\mathcal{P}}
\def\cP{\mathcal{P}}
\def\a{\mathcal{A}}
\def\h{\mathcal{H}}
\def\c{\mathcal{C}}
\def\y{\mathcal{Y}}
\def\e{\mathcal{E}}
\def\v{\mathcal{V}}
\def\z{\mathcal{Z}}
\def\x{\mathcal{X}}
\def\s{\mathcal{S}}
\def\g{\mathcal{G}}
\def\u{\mathcal{U}}
\def\w{\mathcal{W}}
\def\i{\mathcal{I}}
\def\j{\mathcal{J}}
\def\b{\mathcal{B}}

\def\lan{\langle}
\def\llb{(\!(}
\def\ran{\rangle}
\def\rrb{)\!)}
 \def\dim{{\hbox{\rm dim}}_{\mathbb C}\,}
\def\lng{\hbox{\rm{\tiny lng}}}
\def\sht{\hbox{\rm{\tiny sht}}}
\def\sph{\hbox{\rm{\tiny sph}}}
\def\inv{\hbox{\rm{\tiny inv}}}

\def\br#1{\langle #1 \rangle}

\def\rank{\hbox{rank}}
\def\gl{\mathfrak{gl}_N}

\newcommand{\Aut}{\operatorname{Aut}}
\newcommand{\Hom}{\operatorname{Hom}}
\newcommand{\End}{\operatorname{End}}
\newcommand{\Ind}{\operatorname{Ind}}
\newcommand{\ad}{\operatorname{ad}}
\newcommand{\pr}{\operatorname{pr}}
\newcommand{\aweyl}{\tilde{\mathbb S}_n}
\newcommand{\hec}{{\mathcal H}^t_n}
\newcommand{\Func}{{\mathcal F}({\mathbb C}^n,{\mathcal H}^t_n)}
\newcommand{\tr}{\operatorname{tr}}
\newcommand{\Out}{\operatorname{Out}}
\newcommand{\Rad}{\operatorname{Rad}}
\newcommand{\Spec}{\operatorname{Spec}}
\newcommand{\id}{\operatorname{id}}
\newcommand{\Int}{\operatorname{Int}}
\newcommand{\ct} {\operatorname{ct}}

\newcommand{\rat}{{\mathbb Q}}
\newcommand{\real}{{\mathbb R}}
\newcommand{\cplx}{{\mathbb C}}
\newcommand{\zint}{{\mathbb Z}}

\newcommand{\sq}{\phantom{1}\hfill$\qed$}
\newcommand{\Rea}{\Re}
\newcommand{\Ima}{\Im}

\newcommand{\st}{\bowtie}
\newcommand{\modd}{\mbox{\,mod\,}}
\newcommand{\lr}{\langle}
\newcommand{\rr}{\rangle}
\newcommand{\eps}{\varepsilon}
\newcommand{\phk}{\phi^{(k)}}
\newcommand{\psk}{\psi^{(k)}}
\newcommand{\Res}{\mbox{Res}\;}
\newcommand{\sgn}{\mbox{sgn}}
\newcommand{\mn} {\left\{ \begin{array}{c}m\\
n\end{array}\right\}}

\def\sX{\mathscr{X}}
\def\sH{\mathscr{H}}
\def\sY{\mathscr{Y}}
\def\TT{\mathfrak{T}}
\def\JJ{\mathfrak{J}}
\def\HH{\mathfrak{H}}
\def\FF{\mathfrak{F}}
\def\GG{\mathfrak{G}}
\def\CC{\mathfrak{C}}
\def\LL{\mathfrak{L}}

\def\BB{\mathfrak{B}}
\def\AA{\mathfrak{A}}
\def\ZZ{\mathfrak{Z}}
\def\HH{\hbox{${\mathcal H}$\kern-5.2pt${\mathcal H}$}}
\def\HHH{\hbox{${\mathbb H}$\kern-4.2pt${\mathbb H}$}}
\def\tHH{\widetilde{\HH\ }}

\font\smm=msbm10 at 12pt 
\def\symbol#1{\hbox{\smm #1}}
\def\lsmash{{\symbol n}}
\def\rsmash{{\symbol o}}
\def\#{\sharp}

\font\tenbf=cmbx10
\font\tenrm=cmr10
\font\tenit=cmti10
\font\ninebf=cmbx9
\font\ninerm=cmr9
\font\nineit=cmti9
\font\eightbf=cmbx8
\font\eightrm=cmr8
\font\eightit=cmti8
\font\sevenrm=cmr7
\font\sevenbf=cmbx7


%% file: D6-9-22.tex
{\tiny
\begin{longtable}{|c|c|c|}
\hline
\{\ \}\ :\ 24,0 & 
+\{ 47\}:\ 24 & 
\{ 47\}\ :\ 23,1\\ 
+\{ 25\}:\ na & 
+\{ 38\}:\ 24 & 
+\{ 41\}:\ 23\\ 
\{ 25\}\ :\ na,2 & 
\{ 38\}\ :\ 23,2 & 
+\{ 25, 38\}:\ 24\\ 
+\{ 38, 41\}:\ 23 & 
\{ 41\}\ :\ 22,2 & 
+\{ 25, 41\}:\ na\\ 
+\{ 35\}:\ 23 & 
+\{ 38, 41\}:\ 22 & 
\{ 3\}\ :\ na,3\\ 
\{ 25, 38\}\ :\ 23,3 & 
+\{ 3, 38\}:\ na & 
+\{ 16\}:\ 24\\ 
+\{ 25, 38, 41\}:\ 23 & 
\{ 25, 41\}\ :\ na,3 & 
\{ 35\}\ :\ 22,3\\ 
+\{ 25, 35\}:\ 23 & 
+\{ 35, 38\}:\ 22 & 
\{ 38, 41\}\ :\ 21,3\\
\end{longtable}
}